\documentclass[12pt,reqno]{amsart}

\usepackage[colorlinks,linkcolor={blue},
citecolor={blue},urlcolor={blue}]{hyperref}

\usepackage{mathrsfs}
\usepackage{bbm}
\usepackage{amsfonts}
\usepackage{amsfonts}
\usepackage{amsfonts}
\usepackage{amsfonts}
\usepackage{amsfonts}
\usepackage{amsfonts}
\usepackage{amsfonts}
\usepackage{amsfonts}
\usepackage{amsfonts}
\usepackage{amsfonts}
\usepackage{amsfonts}
\usepackage{amsfonts}
\usepackage{amsfonts}
\usepackage{amsfonts}
\usepackage{amsfonts}
\usepackage{amsfonts}
\usepackage{amsfonts}
\usepackage{amsfonts}
\usepackage{amsfonts}
\usepackage{}
\usepackage{pifont}
\usepackage[shortlabels]{enumitem}
\usepackage{empheq}
\usepackage{amsfonts}
\usepackage{amsmath, amsrefs, amsthm, amssymb}
\usepackage[active]{srcltx}		
\usepackage{pdfsync}					
\usepackage{graphicx}
\usepackage [latin1]{inputenc}
\usepackage{bbm}
\usepackage{amsfonts}
\usepackage{amsthm}
\usepackage{graphicx}
\usepackage{framed}
\usepackage{multicol,multienum}
\usepackage{graphicx} 
\usepackage{booktabs} 
\usepackage{mathrsfs}
\usepackage{tcolorbox}
\usepackage{color}
\usepackage{xcolor}

\newtheorem{theorem}{Theorem}[section]
\newtheorem{proposition}[theorem]{Proposition}
\newtheorem{lemma}[theorem]{Lemma}
\newtheorem{corollary}[theorem]{Corollary}

\theoremstyle{definition}
\newtheorem{definition}[theorem]{Definition}

\numberwithin{equation}{section}

\begin{document}

\title [Well-posedness and invariant measures for the SLLBar equation]{Well-posedness and invariant measures for the stochastically perturbed Landau-Lifshitz-Baryakhtar equation}
\author{Fan Xu}
\address{School of Mathematics and Statistics, Hubei Key Laboratory of Engineering Modeling  and Scientific Computing, Huazhong University of Science and Technology,  Wuhan 430074, Hubei, P.R. China.}
\email{d202280019@hust.edu.cn (F. Xu)}

\author{Lei Zhang}
\address{School of Mathematics and Statistics, Hubei Key Laboratory of Engineering Modeling  and Scientific Computing, Huazhong University of Science and Technology,  Wuhan 430074, Hubei, P.R. China.}
\email{lei\_zhang@hust.edu.cn (L. Zhang)}

\author{Bin Liu}
\address{School of Mathematics and Statistics, Hubei Key Laboratory of Engineering Modeling  and Scientific Computing, Huazhong University of Science and Technology,  Wuhan 430074, Hubei, P.R. China.}
\email{binliu@mail.hust.edu.cn (B. Liu)}



\begin{abstract}
In this paper, we study the initial-boundary value problem for the stochastic Landau-Lifshitz-Baryakhtar (SLLBar) equation with  Stratonovich-type noise in   bounded domains $\mathcal{O}\subset\mathbb{R}^d$, $d=1,2,3$. Our main results can be briefly described as follows: (1) for $d=1,2,3$ and any $\mathbf{u}_0\in\mathbb{H}^1$, the SLLBar equation admits a unique local-in-time pathwise weak solution; (2) for $d=1$ and small-data $\mathbf{u}_0\in\mathbb{H}^1$, the SLLBar equation has a unique global-in-time pathwise weak solution and at least one invariant measure; (3) for $d=1,2$ and small-data $\mathbf{u}_0\in\mathbb{L}^2$, the SLLBar equation possesses a unique global-in-time pathwise very weak solution and at least one invariant measure, while for $d=3$ only the existence of martingale solution is obtained due to the loss of pathwise uniqueness.
\end{abstract}

\maketitle
\section{Introduction}\label{sec1}
The study of the theory of ferromagnetism was initiated by Weiss \cite{ad1}. In 1935, Landau and Lifshitz developed the dispersive theory of magnetization of ferromagnets and introduced the well-known Landau-Lifshitz (LL) equation of ferromagnetic spin chains \cite{30}. In 1955,  Gilbert \cite{23} further developed the theory of ferromagnetism and proposed the Landau-Lifshitz-Gilbert (LLG) equation to describe the evolution of the spin magnetic moment in magnetic systems, particularly the precession and dissipation behavior under the influence of an external magnetic field. This LLG equation can be described by the following partial differential  equation:
\begin{equation}\label{1sys}
\begin{aligned}
&\frac{\partial \mathbf{u}}{\partial t}=-\lambda_1\mathbf{u}\times \mathbf{H}_{\textrm{eff}}-\lambda_2\mathbf{u}\times(\mathbf{u}\times \mathbf{H}_{\textrm{eff}}),
\end{aligned}
\end{equation}
where the unknown quantity $\mathbf{u}(t,x)\in\mathbb{R}^3$ denotes the magnetization vector of a magnetic body $\mathcal{O}\subset\mathbb{R}^d$, $d=1,2,3$. $\mathbf{H}_{\textrm{eff}}$ denotes effective field which consists of the external magnetic field, the demagnetizing field (magnetic field due to the magnetization) and some quantum mechanical effects, etc. Here, the effective field is taken to be $\mathbf{H}_{\textrm{eff}}=\Delta \mathbf{u}\in\mathbb{R}^3$. The vector cross product of $\textbf{a}$ and $\textbf{b}$ in $\mathbb{R}^3$ is given by $\textbf{a} \times \textbf{b}$. The real numbers $\lambda_1>0$ and $\lambda_2>0$ are the gyromagnetic ratio and a phenomenological damping parameter, respectively.  It is well known that the LLG equation \eqref{1sys} is valid for temperatures below the the critical (so-called Curie) temperature $T_c$. And it has been widely applied in the study of magnetic materials, particularly in fields such as magnetic storage technology and magnetic nanomaterials \cite{ad3,ad2}. From the mathematical point of view, the existence, uniqueness, and regularity of solutions to the LLG equation have been extensively discussed, see for example \cite{1,15,20,26} and the references cited therein. It is worth noting that, to restore a more realistic physical background, Brze\'{z}niak et al. \cite{8,9} first introduced Gaussian-type noise into the LLG equation and studied the well-posedness problem for associated stochastic partial differential equation. Subsequently, in 2019, Brze\'{z}niak et al. \cite{11} further discussed the more general stochastic LLG equation driven by L\'{e}vy noise. Very recently, they also established  the existence of global solutions to coupled systems of stochastic LLG equations and Maxwell's equations \cite{14}.

To effectively handle the high-temperature situation that is invalid in the LLG equation, Garanin \cite{22} introduced a thermodynamically consistent approach and derived the the following Landau-Lifshitz-Bloch (LLB) equation for ferromagnets:
\begin{equation}\label{121}
\begin{aligned}
& \frac{\partial \mathbf{u}}{\partial t}=\gamma\mathbf{u}\times \mathbf{H}_{\textrm{eff}}+L_1\frac{1}{|\mathbf{u}|^2}(\mathbf{u}\cdot\mathbf{H}_{\textrm{eff}})\mathbf{u}-L_2\frac{1}{|\mathbf{u}|^2}\mathbf{u}\times(\mathbf{u}\times\mathbf{H}_{\textrm{eff}}),
\end{aligned}
\end{equation}
where $\gamma>0$ is the gyromagnetic ratio, and $L_1$ and $L_2$ are the longitudinal and transverse damping parameters, respectively. The effective field considered in \eqref{121} is formulated by
$\mathbf{H}_{\text {eff }}=\Delta \mathbf{u}-\frac{1}{\chi_{\|}} (1+\frac{3}{5} \frac{T}{T-T_c}|\mathbf{u}|^2 ) \mathbf{u}$, where $\chi_{\|}$ is the longitudinal susceptibility.
The LLB equation essentially interpolates between the LLG equation at low temperatures and the Ginzburg-Landau theory of phase transitions.  The introduction of the LLB equation has enriched our understanding of the behavior of magnetic materials, particularly in the study of magnetic nanoparticles and single-molecule magnets \cite{ad4}. In mathematics, the existence and regularity properties for LLB equation have been studied in \cite{31,32}. Naturally, similar to the LLG equation, the LLB equation affected by random noises arising from the environment has also been investigated by several authors. For instance, Jiang et al. \cite{28} demonstrated the existence of at least one martingale solution to a stochastic LLB equation within a three-dimensional bounded domain. Building upon this result, Brze\'{z}niak et al. \cite{13} subsequently proved the existence of at least one invariant measure for the stochastic LLB equation within one-dimensional and two-dimensional bounded domains.

Nevertheless, neither the LLG nor LLB equations fail to explain the certain experimental data and microscopic calculations, such as the nonlocal damping observed in magnetic metals and crystals \cite{17,43}, or the higher-than-expected spin wave decrement for short-wave magnons \cite{4}. To overcome this problem, based on Onsager's relations, Baryakhtar \cite{2,3,4} extended the LLG and LLB equations so that he introduced the so-called Landau-Lifshitz-Baryakhtar (LLBar) equation \cite{17,18,42}. The LLBar equation in its most general form \cite{4,42} reads
$$
\frac{\partial \mathbf{u}}{\partial t}=-\lambda_1 \mathbf{u} \times \mathbf{H}_{\mathrm{eff}}+\mathbf{\Lambda}_r \cdot \mathbf{H}_{\mathrm{eff}}-\mathbf{\Lambda}_{e, i j} \frac{\partial^2 \mathbf{H}_{\mathrm{eff}}}{\partial \mathbf{x}_i \partial \mathbf{x}_j},
$$
where $\mathbf{u}$ represents the magnetisation vector, $\mathbf{\Lambda}_r$ and $\mathbf{\Lambda}_e$ denote the relaxation tensor and the exchange tensor, respectively. Since for a polycrystalline, amorphous soft magnetic materials and magnetic metals at moderate temperature, where nonlocal damping and longitudinal relaxation are significant, the simplified form of the LLBar equation is as follows \cite{17,42}:
\begin{equation}\label{2sys}
\begin{aligned}
&\frac{\partial \mathbf{u}}{\partial t}=-\lambda_1\mathbf{u}\times \mathbf{H}_{\textrm{eff}}+\lambda_r\mathbf{H}_{\textrm{eff}}-\lambda_e\Delta\mathbf{H}_{\textrm{eff}},
\end{aligned}
\end{equation}
where the positive constants $\lambda_1$, $\lambda_r$, and $\lambda_e$ are the electron gyromagnetic ratio, relativistic damping constant, and exchange damping constant, respectively. The effective field $\mathbf{H}_{\textrm{eff}}$ in \eqref{2sys} is given by
\begin{equation}\label{141}
\begin{aligned}
\mathbf{H}_{\textrm{eff}}=\Delta \mathbf{u}+\frac{1}{2\chi}(1-|\mathbf{u}|^2)\mathbf{u},
\end{aligned}
\end{equation}
where $\chi>0$ is the magnetic susceptibility of the material.  It is worth mentioning that various micromagnetic simulations have demonstrated that the LLBar equation  agrees with some of the observed experimental phenomena in micromagnetics, especially regarding ultrafast magnetization at high temperatures, see for instance \cite{17,42,43}. To our knowledge, there are relatively few mathematical analysis results available for the LLBar equation \eqref{2sys} except the recent work \cite{39}. If the exchange interaction is dominant as the case for ordinary ferromagnetic material,  Soenjaya and Tran \cite{39} demonstrated that \eqref{2sys} possesses a unique global weak/strong solution within bounded domains in one, two, and three dimensions.

Physically speaking, in many practical applications the effective magnetic fields  $\mathbf{H}_{\textrm{eff}}$ in \eqref{2sys} are inevitably influenced by random factors from surroundings, such as the thermal fluctuations, magnetic field fluctuations, and external sources of noise. These noise sources can introduce randomness into the effective field, which in turn affect the evolution of the spin magnetic moment. On the other hand, in the theory of ferromagnetism, describing the phase transitions between different equilibrium states induced by thermal fluctuations of the field $\mathbf{H}_{\textrm{eff}}$ is an important problem. Being inspired by these reasons, it is necessary to appropriately modify model \eqref{2sys} by incorporating the random fluctuations of the effective field $\mathbf{H}_{\textrm{eff}}$ into the dynamics of magnetization $\mathbf{u}$, so as to describe the phase transitions between ferromagnetic equilibrium states induced by noise. As a matter of fact, the initiative to analyze noise-induced transitions was started by N\'{e}el \cite{34}, and subsequent advancements were made in \cite{6,29} and others. Taking the ideas from \cite{8,9,28}, one of the effective ways is to perturb the effective field $\mathbf{H}_{\textrm{eff}}$ by adding a Gaussian-type stochastic external forcing, that is, to make the substitution
$
 \mathbf{H}_{\textrm{eff}}  ~~\longmapsto ~~\mathbf{H}_{\textrm{eff}}+\xi,
$
where $\xi$ denotes the white noise with respect to time variable. Then the perturbed LLBar equation \eqref{2sys} can be formulated as
\begin{equation}\label{3sys}
\begin{aligned}
&\frac{\partial \mathbf{u}}{\partial t}=-\lambda_1\mathbf{u}\times (\mathbf{H}_{\textrm{eff}}+\xi)+\lambda_r(\mathbf{H}_{\textrm{eff}}+\xi)-\lambda_e\Delta(\mathbf{H}_{\textrm{eff}}+\xi).
\end{aligned}
\end{equation}
For mathematical qualitative analysis, let us make some assumptions on the random noise $\xi$. It is well known \cite{37} in the theory of SPDEs that a rigorous interpretation of $\xi$ is via the relationship $\xi=\dot{W}$, where $(W(t),t\geq0)$ is a Wiener process defined in a probability space. In particular, the Gaussian noise we considered in this paper takes the following form
\begin{equation}\label{131}
\begin{aligned}
\xi(t)=\sum_{j=1}^{\infty}\mathbf{h}_j\circ\,\frac{dW_j(t)}{dt},
\end{aligned}
\end{equation}
where $\{W_{j},j\geq 1\}$ is a family of independent real-valued Wiener processes, and $\{\mathbf{h}_j,j\geq 1\}$ are space-dependent coefficients satisfying suitable regularity conditions. The term $\circ\mathrm{d}W_j$ in \eqref{131} need to be understand in the Stratonovich sense. In the sequel, we assume that
\begin{equation}\label{condh}
\begin{aligned}
\sum_{j=1}^{\infty}\|\mathbf{h}_j\|_{\mathbb{W}^{3,2}(\mathcal{O})}^2\leq C_{\mathbf{h}}<\infty.
\end{aligned}
\end{equation}

If the exchange interaction is dominant as the ordinary ferromagnetic material, by plugging the formulation \eqref{141} and \eqref{131} into \eqref{3sys}, we obtain the stochastic version of the LLBar equation considered in this work:
\begin{equation}\label{sys1}
\left\{
\begin{aligned}
&\mathrm{d}\mathbf{u}= \left[\beta_1\Delta \mathbf{u}-\beta_2\Delta^2\mathbf{u}+\beta_3(1-|\mathbf{u}|^2)\mathbf{u}-\beta_4\mathbf{u}\times\Delta\mathbf{u}+\beta_5\Delta(|\mathbf{u}|^2\mathbf{u})\right]\,\mathrm{d}t\\
&\quad \quad +\sum_{j=1}^{\infty}(-\mathbf{u}\times \mathbf{h}_j+\mathbf{h}_j-\Delta\mathbf{h}_j)\circ\mathrm{d}W_j(t),&&\textrm{in}~\mathbb{R}_+\times\mathcal{O},\\
&\frac{\partial\mathbf{u}}{\partial\mathbf{n}}=\frac{\partial\Delta\mathbf{u}}{\partial\mathbf{n}}=0,&&\textrm{on}~\mathbb{R}_+\times\partial\mathcal{O},\\
&\mathbf{u}(0)=\mathbf{u}_0,&&\textrm{in}~\mathcal{O},
\end{aligned}
\right.
\end{equation}
where $\mathcal{O}\subset\mathbb{R}^d$, $d=1,2,3$, is a bounded smooth domain, and $\mathbf{n}$ denotes the exterior unit normal vector of the boundary $\partial\mathcal{O}$. $\beta_1=\lambda_r-\frac{\lambda_e}{2\chi}$ is a real constant (may be positive or negative), and $\beta_2,...,\beta_5$ are positive constants.

To our best knowledge, there seems to be no results concerning the mathematical analysis for the LLBar by considering the random noises arising from environment. The main objective of this paper is to provide a ground for the well-posedness for the stochastic LLBar equation under proper assumptions. Meanwhile, when the equation permits a global pathwise solution, we also investigate the existence of invariance measure to the initial-boundary problem.
Let us mention that there have some works considering the stochastic LLG and stochastic LLB equations associated to \eqref{1sys} and \eqref{121}, respectively, see for example \cite{8,9,11,13,28}. However, due to the appearance of the important nonlocal damping and longitudinal relaxation in \eqref{sys1}, the LLBar equation becomes a forth-order parabolic SPDEs that involves more complicated nonlinear structure, which makes the derivation of several key estimations (especially in stochastic setting) more subtle, and we shall overcome these problems by virtue of some techniques from the stochastic analysis.

To state the main results, let us define some spaces which will be used frequently in the sequel.

Let $\mathcal{O}\subset\mathbb{R}^d$, $d=1,2,3$, be an open bounded domain with smooth boundary. The function space $\mathbb{L}^p:=\mathbb{L}^p(\mathcal{O};\mathbb{R}^3)$ denotes the space of p-th integrable functions taking values in $\mathbb{R}^3$ and $\mathbb{W}^{k,p}:=\mathbb{W}^{k,p}(\mathcal{O};\mathbb{R}^3)$ denotes the Sobolev space of functions on $\mathcal{O}$ taking values in $\mathbb{R}^3$. In particular, let $\mathbb{H}^p:=\mathbb{W}^{2,p}$. Let $X$ and $Y$ be two Banach spaces. The symbol $\langle\cdot,\cdot\rangle_{X^*,X}$ stands for the standard duality pairing, where $X^*:=\mathcal{L}(X;\mathbb{R})$ is the dual space of $X$. In particular, if $X$ is a Hilbert space, then the symbol $(\cdot,\cdot)_X$ denotes the scalar product. Let $Q_w$ be a   Hilbert space $Q$ endowed with the weak topology, and $C([0,T];Q_w)$ be the space of weakly continuous functions $f:[0,T]\rightarrow Q$ with the weakest topology. Let $\mathbb{B}_w^k(R)$ be the ball $\mathbb{B}^k(R):=\{f\in\mathbb{H}^k:\|f\|_{\mathbb{H}^k}\leq R\}$ endowed with the weak topology. Then $\mathbb{B}_w^k(R)$ is metrizable \cite{5-1}. Let
$
C([0,T];\mathbb{B}_w^k(R)):=\{f\in C([0,T];\mathbb{H}_w^k):\sup_{t\in[0,T]}\|f\|_{\mathbb{H}^k}\leq R\}.
$
The space $\left(C([0,T];\mathbb{B}_w^k(R)),\rho\right)$ is a complete metric space with
$
\rho(f,g)=\sup_{t\in[0,T]}q(f(t),g(t)),
$
where $q$ is the metric compatible with the weak topology on $\mathbb{B}^k$.

Now, let us provide the definitions of the solutions to the SLLBar equation \eqref{sys1}. In the first one, we assume that the initial data take values in $\mathbb{H}^1$.

\begin{definition}\label{def1}
Suppose that $\mathbf{u}_0\in \mathbb{H}^1$. Fix a stochastic basis $\left(\Omega,\mathcal{F},\mathbb{F}:=\{\mathcal{F}_t\}_{t\geq0},\mathbb{P},W\right)$.
\begin{enumerate}
\item [(1)]  A \textsf{local pathwise weak solution} of \eqref{sys1} is a pair $(\mathbf{u},\tau)$, where $\tau$ is strictly positive stopping time relative to $\mathbb{F}$, and
$
\mathbf{u}(\cdot\wedge\tau)\in L^{\infty}(0,T;\mathbb{H}^1)\cap L^2(0,T;\mathbb{H}^3)
$, $\mathbb{P}$-a.s. Moreover there holds $\mathbb{P}$-a.s.,
\begin{equation}\label{1def1}
\begin{split}
&(\mathbf{u}(t\wedge\tau),\phi)_{\mathbb{L}^2}=(\mathbf{u}_0,\phi)_{\mathbb{L}^2}-\beta_1\int_0^{t\wedge\tau}\left(\nabla\mathbf{u},\nabla\phi\right)_{\mathbb{L}^2}\,\mathrm{d}s+\beta_2\int_0^{t\wedge\tau}\left(\nabla\Delta\mathbf{u},\nabla\phi\right)_{\mathbb{L}^2}\,\mathrm{d}s\\
&\quad +\beta_3\int_0^{t\wedge\tau}\left((1-|\mathbf{u}|^2)\mathbf{u},\phi\right)_{\mathbb{L}^2}\,\mathrm{d}s+\beta_4\int_0^{t\wedge\tau}\left(\mathbf{u}\times\nabla\mathbf{u},\nabla\phi\right)_{\mathbb{L}^2}\,\mathrm{d}s\\
&\quad-\beta_5\int_0^{t\wedge\tau}\left(\nabla\left(|\mathbf{u}|^2\mathbf{u}\right),\nabla\phi\right)_{\mathbb{L}^2}\,\mathrm{d}s+\sum_{j=1}^{\infty}\int_0^{t\wedge\tau}\left(-\mathbf{u}\times \mathbf{h}_j+\mathbf{h}_j-\Delta\mathbf{h}_j,\phi\right)_{\mathbb{L}^2}\circ\mathrm{d}W_j(s),
 \end{split}
\end{equation}
for every $t\in[0,T]$ and $\phi\in \mathbb{H}^1$.

\item [(2)] The local pathwise weak solutions are said to be \textsf{unique}, if given any two pair of local pathwise weak solutions $(\mathbf{u}_1,\tau_1)$ and $(\mathbf{u}_2,\tau_2)$ with the same initial value, then
\begin{equation*}
\begin{split}
\mathbb{P}\{\mathbf{u}_1(t,x)=\mathbf{u}_2(t,x),~\forall(t,x)\in[0,\tau_1\wedge\tau_2]\times\mathcal{O}\}=1.
 \end{split}
\end{equation*}
\item [(3)] A \textsf{maximal pathwise weak solution} is a triple $\left(\mathbf{u},\tau,\{\tau^n\}_{n\geq1}\right)$, if each pair $\left(\mathbf{u},\{\tau^n\}_{n\geq1}\right)$ is a local pathwise weak solution, $\tau^n$ increases with $\lim_{n\rightarrow\infty}\tau^n=\tau$ such that
\begin{equation*}
\begin{split}
\lim_{n\rightarrow\infty}\sup_{t\in[0,\tau^n]}\|\nabla\mathbf{u}(t)\|_{\mathbb{L}^2}=\infty~\textrm{on the set}~\{\tau<\infty\}.
 \end{split}
\end{equation*}
\item [(4)] $\mathbf{u} $ is said to be a \textsf{global pathwise weak solution} of \eqref{sys1} if for every $T>0$, $\mathbf{u}(\cdot)\in L^{\infty}(0,T;\mathbb{H}^1)\cap L^2(0,T;\mathbb{H}^3)$ and satisfies the identity \eqref{1def1} $\mathbb{P}$-a.s.
\end{enumerate}
\end{definition}

We also consider the initial-boundary value problem with initial data in $\mathbb{L}^2$. In present case, the solution naturally has lower regularity than the ones in Definition \ref{def1}, so we introduce the following weaker definition.

\begin{definition}\label{def2}
Suppose that $\mathbf{u}_0\in \mathbb{L}^2$. Fix a stochastic basis $\left(\Omega,\mathcal{F},\mathbb{F},\mathbb{P},W\right)$. A $ \mathbb{F}$-predictable process $\mathbf{u} $ is said to be a \textsf{global pathwise very weak solution} of \eqref{sys1} if for every $T>0$,
$
\mathbf{u}(\cdot)\in L^{\infty}(0,T;\mathbb{L}^2)\cap L^2(0,T;\mathbb{H}^2),
$
and for every $t\in[0,T]$, $\phi\in \mathbb{H}^2$ and $\mathbb{P}$-a.s.
\begin{equation}\label{1def2}
\begin{split}
(\mathbf{u}(t),\phi)_{\mathbb{L}^2}=&(\mathbf{u}_0,\phi)_{\mathbb{L}^2}-\beta_1\int_0^{t}\left(\nabla\mathbf{u},\nabla\phi\right)_{\mathbb{L}^2}\,\mathrm{d}s-\beta_2\int_0^{t}\left(\Delta\mathbf{u},\Delta\phi\right)_{\mathbb{L}^2}\,\mathrm{d}s\\
&+\beta_3\int_0^{t}\left((1-|\mathbf{u}|^2)\mathbf{u},\phi\right)_{\mathbb{L}^2}\,\mathrm{d}s+\beta_4\int_0^{t}\left(\mathbf{u}\times\nabla\mathbf{u},\nabla\phi\right)_{\mathbb{L}^2}\,\mathrm{d}s\\
&+\beta_5\int_0^{t}\left(|\mathbf{u}|^2\mathbf{u},\Delta\phi\right)_{\mathbb{L}^2}\,\mathrm{d}s+\sum_{j=1}^{\infty}\int_0^{t}\left(-\mathbf{u}\times \mathbf{h}_j+\mathbf{h}_j-\Delta\mathbf{h}_j,\phi\right)_{\mathbb{L}^2}\circ\mathrm{d}W_j(s).
\end{split}
\end{equation}
\end{definition}

Now we can state the main results of this paper. The first main result states that for all $d=1,2,3$ and any initial data in $\mathbb{H}^1$ without size limit,   the SLLBar equation \eqref{sys1} admits a unique probabilistically strong solution up to a almost-surely finite stopping time.

\begin{theorem}[\textsf{Local large-data solution in $\mathbb{H}^1$}]\label{the1}
Let $\mathcal{O}\subset\mathbb{R}^d$, $d=1,2,3$, be a bounded domain with $C^{2,1}$-boundary. Assume that $\mathbf{u}_0\in\mathbb{H}^1$. Then there exists a unique local maximal pathwise weak solution $\left(\mathbf{u},\tau,\{\tau^n\}_{n\geq1}\right)$ to \eqref{sys1} in the sense of Definition  \ref{def1} such that
$
\mathbf{u}(\cdot\wedge\tau)\in L^p\left(\Omega;L^{\infty}(0,T;\mathbb{H}^1)\cap L^2(0,T;\mathbb{H}^3)\right),
$
for every $p\geq1$.
\end{theorem}

Our second result concerns the global-in-time solvability for \eqref{sys1} in $\mathbb{H}^1$. It is shown that small initial data stimulate to the existence of global solutions. Meanwhile, we also prove that the associated equation possesses an invariance measure.

\begin{theorem}[\textsf{Global small-data solution in $\mathbb{H}^1$}]\label{the2} Let $\mathcal{O}\subset\mathbb{R}^d$, $d=1$, be a bounded domain with $C^{2,1}$-boundary. Assume that $
\|\mathbf{u}_0\|_{\mathbb{H}^1}\leq C_0$ for some $C_0>0$.
Then

(1) there exists a unique global pathwise weak solution to \eqref{sys1} in the sense of Definition \ref{def1} such that
$
\mathbf{u}(\cdot)\in L^p\left(\Omega;L^{\infty}(0,T;\mathbb{H}^1)\cap L^2(0,T;\mathbb{H}^3)\right),
$
and
\begin{equation}\label{1then}
\begin{split}
\mathbb{E}\|\mathbf{u}\|_{W^{\alpha,p}(0,T;(\mathbb{H}^{1})^*)}^q+\mathbb{E}\|\mathbf{u}\|_{L^{\infty}(0,T;\mathbb{H}^1)\cap L^{2}(0,T;\mathbb{H}^3)}^p\leq C,
\end{split}
\end{equation}
for every $p$, $q\geq1$ and $\alpha\in(0,\frac{1}{2})$;

(2) there exists at least one invariant measure for equation \eqref{sys1}.
\end{theorem}

In our third main result, we prove that \eqref{sys1} admits a unique global pathwise solution with   initial data $\mathbf{u}_0\in\mathbb{L}^2$ in dimension one and two. In both of the cases, we show that the associated equation has an invariant measure. In dimension three, we prove that \eqref{sys1} admits a global probabilistically weak solution while  leaving the uniqueness part to be open.

\begin{theorem} [\textsf{Global small-data solution in $\mathbb{L}^2$}] \label{the3}
Assume that there is a constant $C_0'>0$ such that $\|\mathbf{u}_0\|_{\mathbb{L}^2}\leq C_0'$. Then we have

(1) if $d=3$, there exists at least one global martingale very weak solution of \eqref{sys1};

(2) if $d=1,2$, there exists a unique global pathwise very weak solution of \eqref{sys1} in the sense of Definition \ref{def2} such that $\mathbf{u}(\cdot)\in L^p\left(\Omega;L^{\infty}(0,T;\mathbb{L}^2)\cap L^2(0,T;\mathbb{H}^2)\right)$, and
\begin{equation}\label{1the2}
\begin{split}
\mathbb{E}\|\mathbf{u}\|_{W^{\alpha,p}(0,T;(\mathbb{H}^{2})^*)}^q+\mathbb{E}\|\mathbf{u}\|_{L^{\infty}(0,T;\mathbb{L}^2)\cap L^{2}(0,T;\mathbb{H}^2)}^p\leq C,
\end{split}
\end{equation}
for every $p$, $q\geq1$ and $\alpha\in(0,\frac{1}{2})$;

(3) if $d=1,2$, there exists at least one invariant measure for equation \eqref{sys1}.
\end{theorem}

The main framework for proving the main theorems of this paper relies on the classical Faedo-Galerkin approximation combined with the stochastic compactness method. This framework was first used in the work of Flandoli and Gatarek \cite{21}, and has since been effectively applied in discussions of stochastic LLG and stochastic LLB equations \cite{8,9,13,28}. To obtain weak solutions of equation \eqref{sys1}, we first derive uniformly bounded estimates for the solutions of the Galerkin approximation equation. Subsequently, by compactness arguments, Skorohod representation theorem \cite{27,37}, and convergence of approximate solutions pointwise, we obtain weak solutions of the original equation in the sense of probability. To prove that the weak solution in the sense of probability is also a strong solution, it suffices to show that the solution is pathwise unique, which is a direct result of the Yamada-Watanabe theorem \cite{35,44}. It is worth noting that to obtain solutions in the sense of Definition \ref{def1}, we cannot directly derive a priori uniformly bounded estimates for the Galerkin approximation equation of \eqref{sys1}. Fortunately, drawing inspiration from the literature \cite{16,24}, we can introduce a appropriate truncation function to correct the original equation, allowing us to derive uniformly bounded estimates for the modified Galerkin equation \eqref{sys2} and obtain global solutions for the modified equation \eqref{sys1-2}. To transform the modified equation back to the original equation, we simply need to introduce a stopping time to remove the truncation function, thereby ultimately obtaining local the solution of equation \eqref{sys1}. In particular, if considering only the one-dimensional case or aiming solely to obtain weak solutions under Definition \ref{def2}, it is unnecessary to introduce the truncation function. Finally, we prove the existence of invariant measures, which is based on the utilization of the Maslowski-Seidler theorem \cite{33}.

This paper is organized as follows. In section \ref{sec2}, we construct the solutions to an approximate scheme by
the Faedo-Galerkin approximations and prove for them some uniform bounds in various norms.
Section \ref{sec3} is devoted to obtain the tightness of approximate solutions which allows us to use the Skorohod theorem. The main results on well-posedness and invariant measure are proved in sections \ref{sec4}, \ref{sec5} and \ref{sec7}. Some auxiliary lemmas are given in Appendix.

\section{Faedo-Galerkin apprpximation}\label{sec2}
The main aim of this section is first to introduce the approximation equation with solutions in finite-dimensional Hilbert spaces, and then derive some necessary uniform  a priori estimates for the approximation solutions.

Let $\{\mathbf{e}_i\}_{i=1}^{\infty}$ denote an orthonormal basis of $\mathbb{L}^2$ consisting of eigenvectors for the Neumann Laplacian $A=-\Delta$ such that
$
A\mathbf{e}_i=\lambda_i\mathbf{e}_i~~\textrm{in}~~\mathcal{O},~~\textrm{and}~~\frac{\partial\mathbf{e}_i}{\partial\mathbf{n}}=0~~\textrm{on}~~\partial\mathcal{O},
$
where $\lambda_i>0$ are the eigenvalues of $A$, associated with $\mathbf{e}_i$. According to elliptic regularity results, $\mathbf{e}_i$ is smooth up to the boundary, and we have
$
A^2\mathbf{e}_i=\lambda_i^2\mathbf{e}_i~~\textrm{in}~~\mathcal{O},~~\textrm{and}~~\frac{\partial\mathbf{e}_i}{\partial\mathbf{n}}=\frac{\partial\Delta\mathbf{e}_i}{\partial\mathbf{n}}=0~~\textrm{on}~~\partial\mathcal{O}.
$
Let $S_n:=\textrm{span}\{\mathbf{e}_1...,\mathbf{e}_n\}$ and $\Pi:\mathbb{L}^2\rightarrow S_n$ be the orthogonal projection defined by
$
(\Pi f,g)_{\mathbb{L}^2}=(f,g)_{\mathbb{L}^2},~g\in S_n,~f\in\mathbb{L}^2.
$
We note that $\Pi_n$ is self-adjoint and satisfies
$
\|\Pi f\|_{\mathbb{L}^2}\leq\|f\|_{\mathbb{L}^2},~f\in\mathbb{L}^2.
$

To prove the existence of a local martingale weak solution to \eqref{sys1}, we will use the Faedo-Galerkin method and introduce a truncation function $\theta_R(\cdot)$. Fix $R>0$ to be determined, choosing a $C^{\infty}$-smooth nonincreasing truncation function $\theta_R:[0,\infty)\mapsto[0,1]$ such that
\begin{equation*}
\theta_R(x)=\left\{
\begin{aligned}
&1,~\textrm{for}~|x|<R,\\
&0,~\textrm{for}~|x|>2R.
\end{aligned}
\right.
\end{equation*}
We consider the following Galerkin approximation scheme for \eqref{sys1}
\begin{equation}\label{sys2}
\left\{
\begin{aligned}
&\mathrm{d}\mathbf{u}_n=\bigl[\beta_1\Delta \mathbf{u}_n-\beta_2\Delta^2\mathbf{u}_n+\beta_3\Pi_n\left((1-|\mathbf{u}_n|^2)\mathbf{u}_n\right)-\beta_4\Pi_n(\mathbf{u}_n\times\Delta\mathbf{u}_n)\\
&~~+\beta_5\theta_R(\|\nabla\mathbf{u}_n\|_{\mathbb{L}^2})\Pi_n\Delta(|\mathbf{u}_n|^2\mathbf{u}_n)\bigl]\,\mathrm{d}t+\sum_{j=1}^{n}\Pi_n(-\mathbf{u}_n\times \mathbf{h}_j+\mathbf{h}_j-\Delta\mathbf{h}_j)\circ\mathrm{d}W_j(t),&&\textrm{in}~(0,\infty)\times\mathcal{O},\\
&\mathbf{u}_n(0)=\Pi_n\mathbf{u}_0,&&\textrm{in}~\mathcal{O}.
\end{aligned}
\right.
\end{equation}
The existence of a local solution to the SDE \eqref{sys2} is a consequence of the following lemma, whose proof is standard.
\begin{lemma}\label{lem31} For $n\in\mathbb{N}$, define the maps:
\begin{equation*}
\begin{split}
&F_n^1:S_n\ni f\mapsto\Delta f\in S_n,\\
&F_n^2:S_n\ni f\mapsto\Delta^2 f\in S_n,\\
&F_n^3:S_n\ni f\mapsto\Pi_n\left(|f|^2f\right)\in S_n,\\
&F_n^4:S_n\ni f\mapsto\Pi_n\left(f\times\Delta f\right)\in S_n,\\
&F_n^5:S_n\ni f\mapsto\theta_R(\|\nabla f\|_{\mathbb{L}^2})\Pi_n\Delta\left(|f|^2f\right)\in S_n,\\
&G_{nj}:S_n\ni f\mapsto\Pi_n\left(-f\times\mathbf{h}_j+\mathbf{h}_j-\Delta\mathbf{h}_j\right)\in S_n.
\end{split}
\end{equation*}
\end{lemma}
Then $F_n^1$, $F_n^2$ and $G_{nj}$ are globally Lipschitz while $F_n^3$, $F_n^4$ and $F_n^5$ are locally Lipschitz.

Let us recall the relation between the Stratonovich and It\^{o} differentials: if $W_j$ is an $\mathbb{R}$-valued standard Wiener process defined on a certain filtered probability space $\left(\Omega,\mathcal{F},\mathbb{F},\mathbb{P}\right)$ then
\begin{equation*}
\begin{split}
G_{nj}(f)\circ\mathrm{d}W_j(t)&=\frac{1}{2}G_{nj}'(f)\left[G_{nj}(f)\right]+G_{nj}(f)\mathrm{d}W_j(t)\\
&=-\frac{1}{2}\Pi_n\left(G_{nj}(f)\times\mathbf{h}_j\right)+G_{nj}(f)\mathrm{d}W_j(t).
\end{split}
\end{equation*}

We now proceed to prove uniform bounds for the approximate solutions of \eqref{sys2}.
\begin{lemma}\label{lem32} Let $\mathcal{O}\subset\mathbb{R}^d$,~$d=1,2,3$, be a bounded domain with $C^{2,1}$-boundary. Then for any $p\geq1$, $n\in\mathbb{N}$ and every $t\in[0,T]$, there exists a positive constant $C=C(\|\mathbf{u}_0\|_{\mathbb{L}^2},p,\mathbf{h},T)$ independent of $n$ such that
\begin{equation}\label{1lem32}
\begin{split}
&\mathbb{E}\sup_{s\in[0,t]}\|\mathbf{u}_n(s)\|_{\mathbb{L}^2}^{2p}+\mathbb{E}\left(\int_0^t\|\mathbf{u}_n(s)\|_{\mathbb{H}^2}^2\,\mathrm{d}s\right)^p+\mathbb{E}\left(\int_0^t\|\mathbf{u}_n(s)\|_{\mathbb{L}^4}^4\,\mathrm{d}s\right)^p\leq C.
\end{split}
\end{equation}
\end{lemma}
\noindent\textbf{Proof.} Applying the It\^{o} lemma to $\|\mathbf{u}_n\|_{\mathbb{L}^2}^2$, we have
\begin{equation}\label{2lem32}
\begin{split}
&\frac{1}{2}\|\mathbf{u}_n(t)\|_{\mathbb{L}^2}^2+\beta_1\int_0^t\|\nabla \mathbf{u}_n(s)\|_{\mathbb{L}^2}^2\,\mathrm{d}s+\beta_2\int_0^t\|\Delta \mathbf{u}_n(s)\|_{\mathbb{L}^2}^2\,\mathrm{d}s+\beta_3\int_0^t\| \mathbf{u}_n(s)\|_{\mathbb{L}^4}^4\,\mathrm{d}s\\
&=\frac{1}{2}\|\mathbf{u}_n(0)\|_{\mathbb{L}^2}^2+\beta_3\int_0^t\| \mathbf{u}_n(s)\|_{\mathbb{L}^2}^2\,\mathrm{d}s+\beta_5\int_0^t\left(\theta_R(\|\nabla \mathbf{u}_n(s)\|_{\mathbb{L}^2})\Delta\left(|\mathbf{u}_n(s)|\mathbf{u}_n(s)\right),\mathbf{u}_n(s)\right)_{\mathbb{L}^2}\,\mathrm{d}s\\
&-\frac{1}{2}\sum_{j=1}^n\int_0^t\left(G_{nj}(\mathbf{u}_n(s))\times\mathbf{h}_j,\mathbf{u}_n(s)\right)_{\mathbb{L}^2}\,\mathrm{d}s+\frac{1}{2}\sum_{j=1}^n\int_0^t\|G_{nj}(\mathbf{u}_n(s))\|_{\mathbb{L}^2}^2\,\mathrm{d}s\\
&+\sum_{j=1}^n\int_0^t\left(G_{nj}(\mathbf{u}_n(s)),\mathbf{u}_n(s)\right)_{\mathbb{L}^2}\,\mathrm{d}W_j(s).
\end{split}
\end{equation}
Here for the above equality we have used the fact that $\left(\mathbf{u}_n\times\Delta \mathbf{u}_n,\mathbf{u}_n\right)_{\mathbb{L}^2}=0$. Since $\frac{\partial \mathbf{u}_n}{\partial \mathbf{n}}=0$ on $\partial \mathcal{O}$, we use integration by parts to obtain that
\begin{equation}\label{3lem32}
\begin{split}
&\left(\theta_R(\|\nabla \mathbf{u}_n\|_{\mathbb{L}^2})\Delta\left(|\mathbf{u}_n|\mathbf{u}_n\right),\mathbf{u}_n\right)_{\mathbb{L}^2}=-\theta_R(\|\nabla \mathbf{u}_n\|_{\mathbb{L}^2})\left(\nabla\left(|\mathbf{u}_n|\mathbf{u}_n\right),\nabla\mathbf{u}_n\right)_{\mathbb{L}^2}\\
&=-2\theta_R(\|\nabla \mathbf{u}_n\|_{\mathbb{L}^2})\|\mathbf{u}_n\cdot\nabla\mathbf{u}_n\|_{\mathbb{L}^2}^2-\theta_R(\|\nabla \mathbf{u}_n\|_{\mathbb{L}^2})\||\mathbf{u}_n||\nabla\mathbf{u}_n|\|_{\mathbb{L}^2}^2.
\end{split}
\end{equation}
By using the H\"{o}lder inequality and Young's inequality, it follows from \eqref{condh} that
\begin{equation}\label{4lem32}
\begin{split}
&\sum_{j=1}^n\left|\left(G_{nj}(\mathbf{u}_n)\times\mathbf{h}_j,\mathbf{u}_n\right)_{\mathbb{L}^2}\right|\leq\sum_{j=1}^n\|\mathbf{h}_j\|_{\mathbb{L}^{\infty}}\|G_{nj}(\mathbf{u}_n)\|_{\mathbb{L}^{2}}\|\mathbf{u}_n\|_{\mathbb{L}^{2}}\\
&\leq\sum_{j=1}^n\|\mathbf{h}_j\|_{\mathbb{L}^{\infty}}\left(\|\mathbf{h}_j\|_{\mathbb{L}^{\infty}}\|\mathbf{u}_n\|_{\mathbb{L}^{2}}+\|\mathbf{h}_j\|_{\mathbb{L}^{2}}+\|\Delta\mathbf{h}_j\|_{\mathbb{L}^{2}}\right)\|\mathbf{u}_n\|_{\mathbb{L}^{2}}\\
&\leq C_{\mathbf{h}}+C_{\mathbf{h}}\|\mathbf{u}_n\|_{\mathbb{L}^{2}}^2.
\end{split}
\end{equation}
Similarly, we have
\begin{equation}\label{5lem32}
\begin{split}
&\sum_{j=1}^n\|G_{nj}(\mathbf{u}_n)\|_{\mathbb{L}^2}^2\leq\sum_{j=1}^n4\left(\|\mathbf{u}_n\times\mathbf{h}_j\|_{\mathbb{L}^2}^2+\|\mathbf{h}_j\|_{\mathbb{L}^2}^2+\|\Delta\mathbf{h}_j\|_{\mathbb{L}^2}^2\right)\leq C_{\mathbf{h}}+C_{\mathbf{h}}\|\mathbf{u}_n\|_{\mathbb{L}^{2}}^2.
\end{split}
\end{equation}
In addition, using integration by parts, the H\"{o}lder inequality and Young's inequality, we have
\begin{equation}\label{6lem32}
\begin{split}
|\beta_1|\|\nabla \mathbf{u}_n\|_{\mathbb{L}^2}^2=-|\beta_1|(\mathbf{u}_n,\Delta\mathbf{u}_n)_{\mathbb{L}^2}^2\leq\varepsilon\|\Delta\mathbf{u}_n\|_{\mathbb{L}^2}^2+C_{\varepsilon}\|\mathbf{u}_n\|_{\mathbb{L}^2}^2.
\end{split}
\end{equation}
Plugging \eqref{3lem32}-\eqref{6lem32} into \eqref{2lem32} and choosing $\varepsilon$ small enough, we infer that
\begin{equation}\label{7lem32}
\begin{split}
&\|\mathbf{u}_n(t)\|_{\mathbb{L}^2}^2+\int_0^t\|\Delta \mathbf{u}_n\|_{\mathbb{L}^2}^2\,\mathrm{d}s+\int_0^t\| \mathbf{u}_n\|_{\mathbb{L}^4}^4\,\mathrm{d}s+\int_0^t\theta_R(\|\nabla \mathbf{u}_n\|_{\mathbb{L}^2})\|\mathbf{u}_n\cdot\nabla\mathbf{u}_n\|_{\mathbb{L}^2}^2\,\mathrm{d}s\\
&+\int_0^t\theta_R(\|\nabla \mathbf{u}_n\|_{\mathbb{L}^2})\| |\mathbf{u}_n||\nabla \mathbf{u}_n|\|_{\mathbb{L}^2}^2\,\mathrm{d}s\\
&\leq C\|\mathbf{u}_n(0)\|_{\mathbb{L}^2}^2+C_{\mathbf{h}}+C_{\mathbf{h}}\int_0^t\|\mathbf{u}_n\|_{\mathbb{L}^2}^2\,\mathrm{d}s+C\sum_{j=1}^n\int_0^t\left(G_{nj}(\mathbf{u}_n),\mathbf{u}_n\right)_{\mathbb{L}^2}\,\mathrm{d}W_j(s).
\end{split}
\end{equation}
Using the BDG inequality and the H\"{o}lder inequality, we see that for any $p\geq1$
\begin{equation}\label{8lem32}
\begin{split}
&\mathbb{E}\sup_{s\in[0,t]}\left|\sum_{j=1}^n\int_0^s\left(G_{nj}(\mathbf{u}_n(s')),\mathbf{u}_n(s')\right)_{\mathbb{L}^2}\,\mathrm{d}W_j(s')\right|^p\\
&=\mathbb{E}\sup_{s\in[0,s]}\left|\sum_{j=1}^n\int_0^t\left(\mathbf{h}_j+\Delta\mathbf{ h}_j,\mathbf{u}_n\right)_{\mathbb{L}^2}\,\mathrm{d}W_j(s)\right|^p\\
&\leq C\mathbb{E}\left|\sum_{j=1}^n\int_0^t\left(\|\mathbf{h}_j\|_{L^2}^2+\|\Delta\mathbf{ h}_j\|_{L^2}^2\right)\|\mathbf{u}_n\|_{\mathbb{L}^2}^2\,\mathrm{d}s\right|^\frac{p}{2}\\
&\leq C\left(\sum_{j=1}^n\|\mathbf{h}_j\|_{L^2}^2+\|\Delta\mathbf{h}_j\|_{L^2}^2\right)^{\frac{p}{2}}\mathbb{E}\left(\left|\int_0^t\|\mathbf{u}_n\|_{\mathbb{L}^2}^2\,\mathrm{d}s\right|^\frac{p}{2}\right)\leq C_{\mathbf{h}}+C_{\mathbf{h}}\int_0^t\mathbb{E}\|\mathbf{u}_n(s)\|_{\mathbb{L}^2}^{2p}\,\mathrm{d}s.
\end{split}
\end{equation}
Thus it follows from \eqref{7lem32}, \eqref{8lem32} and Jensen's inequality that
\begin{equation*}
\begin{split}
&\mathbb{E}\sup_{s\in[0,t]}\|\mathbf{u}_n(s)\|_{\mathbb{L}^2}^{2p}+\mathbb{E}\left(\int_0^t\|\nabla \mathbf{u}_n(s)\|_{\mathbb{L}^2}^2\,\mathrm{d}s\right)^p+\mathbb{E}\left(\int_0^t\|\Delta \mathbf{u}_n(s)\|_{\mathbb{L}^2}^2\,\mathrm{d}s\right)^p\\
&+\mathbb{E}\left(\int_0^t\|\mathbf{u}_n(s)\|_{\mathbb{L}^4}^4\,\mathrm{d}s\right)^p\leq C+C\int_0^t\mathbb{E}\|\mathbf{u}_n(s)\|_{\mathbb{L}^2}^{2p}\,\mathrm{d}s.
\end{split}
\end{equation*}
The estimate \eqref{1lem32} then follows from the Gronwall lemma, which completes the proof.~~$\Box$

\begin{lemma}\label{lem33} Under the same assumptions as in Lemma \ref{lem32}. For any $p\geq1$, $n\in\mathbb{N}$ and every $t\in[0,T]$, there exists a positive constant $C=C(\|\nabla\mathbf{u}_0\|_{\mathbb{L}^2},p,\mathbf{h},T,R)$ independent of $n$ such that
\begin{equation}\label{1lem33}
\begin{split}
\mathbb{E}\sup_{s\in[0,t]}\|\mathbf{u}_n(s)\|_{\mathbb{H}^1}^{2p}+\mathbb{E}\left(\int_0^t\|\mathbf{u}_n(s)\|_{\mathbb{H}^3}^2\,\mathrm{d}s\right)^p\leq C.
\end{split}
\end{equation}
\end{lemma}
\noindent\textbf{Proof.} Applying the It\^{o} Lemma to $\|\nabla\mathbf{u}_n\|_{\mathbb{L}^2}^2$, we see that
\begin{equation}\label{2lem33}
\begin{split}
&\frac{1}{2}\|\nabla\mathbf{u}_n(t)\|_{\mathbb{L}^2}^2+\beta_1\int_0^t\|\Delta \mathbf{u}_n\|_{\mathbb{L}^2}^2\,\mathrm{d}s+\beta_2\int_0^t\|\nabla\Delta \mathbf{u}_n\|_{\mathbb{L}^2}^2\,\mathrm{d}s+\beta_3\int_0^t\left(\nabla(|\mathbf{u}_n|^2\mathbf{u}_n),\nabla\mathbf{u}_n\right)_{\mathbb{L}^2}\,\mathrm{d}s\\
&=\frac{1}{2}\|\nabla\mathbf{u}_n(0)\|_{\mathbb{L}^2}^2+\beta_3\int_0^t\|\nabla\mathbf{u}_n\|_{\mathbb{L}^2}^2\,\mathrm{d}s+\beta_5\int_0^t\left(\theta_R(\|\nabla \mathbf{u}_n\|_{\mathbb{L}^2})\nabla\Delta\left(|\mathbf{u}_n|\mathbf{u}_n\right),\nabla\mathbf{u}_n\right)_{\mathbb{L}^2}\,\mathrm{d}s\\
&-\frac{1}{2}\sum_{j=1}^n\int_0^t\left(\nabla\left(G_{nj}(\mathbf{u}_n)\times\mathbf{h}_j\right),\nabla\mathbf{u}_n\right)_{\mathbb{L}^2}\,\mathrm{d}s+\frac{1}{2}\sum_{j=1}^n\int_0^t\|\nabla G_{nj}(\mathbf{u}_n)\|_{\mathbb{L}^2}^2\,\mathrm{d}s\\
&+\sum_{j=1}^n\int_0^t\left(\nabla G_{nj}(\mathbf{u}_n),\nabla\mathbf{u}_n\right)_{\mathbb{L}^2}\,\mathrm{d}W_j(s).
\end{split}
\end{equation}
Through direct calculation, we have
\begin{equation}\label{3lem33}
\begin{split}
\left(\nabla(|\mathbf{u}_n|^2\mathbf{u}_n),\nabla\mathbf{u}_n\right)_{\mathbb{L}^2}=2\|\mathbf{u}_n\cdot\nabla\mathbf{u}_n\|_{\mathbb{L}^2}^2+\||\mathbf{u}_n||\nabla\mathbf{u}_n|\|_{\mathbb{L}^2}^2,
\end{split}
\end{equation}
and
\begin{equation}\label{4lem33}
\begin{split}
&\left(\theta_R(\|\nabla \mathbf{u}_n\|_{\mathbb{L}^2})\nabla\Delta\left(|\mathbf{u}_n|\mathbf{u}_n\right),\nabla\mathbf{u}_n\right)_{\mathbb{L}^2}=-\theta_R(\|\nabla \mathbf{u}_n\|_{\mathbb{L}^2})\left(\Delta\left(|\mathbf{u}_n|\mathbf{u}_n\right),\Delta\mathbf{u}_n\right)_{\mathbb{L}^2}\\
&=-\theta_R(\|\nabla \mathbf{u}_n\|_{\mathbb{L}^2})\Bigl[2\left(|\nabla\mathbf{u}_n|^2\mathbf{u}_n,\Delta\mathbf{u}_n\right)_{\mathbb{L}^2}+2\|\mathbf{u}_n\cdot\Delta\mathbf{u}_n\|_{\mathbb{L}^2}^2+\||\mathbf{u}_n||\Delta \mathbf{u}_n|\|_{\mathbb{L}^2}^2\\
&+4\left(\nabla\mathbf{u}_n(\mathbf{u}_n\cdot\nabla\mathbf{u}_n)^{\top},\Delta\mathbf{u}_n\right)_{\mathbb{L}^2}\Bigl].
\end{split}
\end{equation}
Using the H\"{o}lder inequality and Young's inequality as well as the condition \eqref{condh}, we have
\begin{equation}\label{5lem33}
\begin{split}
&\sum_{j=1}^n\left|\left(\nabla\left(G_{nj}(\mathbf{u}_n)\times\mathbf{h}_j\right),\nabla\mathbf{u}_n\right)_{\mathbb{L}^2}\right|=\sum_{j=1}^n\left|\left(\left(G_{nj}(\mathbf{u}_n)\times\mathbf{h}_j\right),\Delta\mathbf{u}_n\right)_{\mathbb{L}^2}\right|\\
&\leq\sum_{j=1}^n\|\mathbf{h}_j\|_{\mathbb{L}^{\infty}}\|G_{nj}(\mathbf{u}_n)\|_{\mathbb{L}^2}\|\Delta\mathbf{u}_n\|_{\mathbb{L}^2}\\
&\leq\sum_{j=1}^n\|\mathbf{h}_j\|_{\mathbb{L}^{\infty}}\left(\|\mathbf{h}_j\|_{\mathbb{L}^{\infty}}\|\mathbf{u}_n\|_{\mathbb{L}^2}+\|\mathbf{h}_j\|_{\mathbb{L}^{2}}+\|\Delta\mathbf{h}_j\|_{\mathbb{L}^{2}}\right)\|\Delta\mathbf{u}_n\|_{\mathbb{L}^2}\\
&\leq\varepsilon\|\Delta\mathbf{u}_n\|_{\mathbb{L}^2}^2+C_{\mathbf{h}}\|\mathbf{u}_n\|_{\mathbb{L}^2}^2+C_{\mathbf{h},\varepsilon}.
\end{split}
\end{equation}
Similarly, it follows that
\begin{equation}\label{6lem33}
\begin{split}
&\sum_{j=1}^n\|\nabla G_{nj}(\mathbf{u}_n)\|_{\mathbb{L}^2}^2\leq\sum_{j=1}^n4\left(\|\nabla\left(\mathbf{u}_n\times\mathbf{h}_j\right)\|_{\mathbb{L}^2}^2+\|\nabla\mathbf{h}_j\|_{\mathbb{L}^2}^2+\|\nabla\Delta\mathbf{h}_j\|_{\mathbb{L}^2}^2\right)\\
&\leq C_{\mathbf{h}}+C_{\mathbf{h}}\|\nabla\mathbf{u}_n\|_{\mathbb{L}^2}^2+C_{\mathbf{h}}\|\mathbf{u}_n\|_{\mathbb{L}^2}^2.
\end{split}
\end{equation}
In addition, by the standard elliptic regularity result with Neumann boundary data \cite{25}, it follows that
\begin{equation}\label{7lem33}
\begin{split}
\|\Delta \mathbf{u}_n\|_{\mathbb{L}^2}^2\leq\|\nabla \mathbf{u}_n\|_{\mathbb{L}^2}\|\nabla\Delta \mathbf{u}_n\|_{\mathbb{L}^2}\leq\varepsilon\|\nabla\Delta \mathbf{u}_n\|_{\mathbb{L}^2}^2+C_{\varepsilon}\|\nabla\mathbf{u}_n\|_{\mathbb{L}^2}^2.
\end{split}
\end{equation}
Then plugging \eqref{3lem33}-\eqref{7lem33} into \eqref{2lem33} and choosing $\varepsilon$ small enough, we have
\begin{equation}\label{8lem33}
\begin{split}
&\|\nabla\mathbf{u}_n(t)\|_{\mathbb{L}^2}^2+\int_0^t\|\nabla\Delta \mathbf{u}_n\|_{\mathbb{L}^2}^2\,\mathrm{d}s+\int_0^t\|\mathbf{u}_n\cdot\nabla\mathbf{u}_n\|_{\mathbb{L}^2}^2\,\mathrm{d}s+\int_0^t\||\mathbf{u}_n||\nabla\mathbf{u}_n|\|_{\mathbb{L}^2}^2\,\mathrm{d}s\\
&+C_1\int_0^t\theta_R(\|\nabla \mathbf{u}_n\|_{\mathbb{L}^2})\||\mathbf{u}_n||\Delta\mathbf{u}_n|\|_{\mathbb{L}^2}^2\,\mathrm{d}s\\
&\leq C\|\nabla\mathbf{u}_n(0)\|_{\mathbb{L}^2}^2+C\int_0^t\|\nabla\mathbf{u}_n\|_{\mathbb{L}^2}^2\,\mathrm{d}s+C\int_0^t\|\mathbf{u}_n\|_{\mathbb{L}^2}^2\,\mathrm{d}s\\
&+C_2\int_0^t\theta_R(\|\nabla \mathbf{u}_n\|_{\mathbb{L}^2})\left(|\mathbf{u}_n||\nabla\mathbf{u}_n|^2,|\Delta\mathbf{u}_n|\right)_{\mathbb{L}^2}\,\mathrm{d}s+\sum_{j=1}^n\int_0^t\left(\nabla G_{nj}(\mathbf{u}_n),\nabla\mathbf{u}_n\right)_{\mathbb{L}^2}\,\mathrm{d}W_j(s).
\end{split}
\end{equation}
Using the H\"{o}lder inequality and Young's inequality, we have
\begin{equation*}
\begin{split}
&C_2\theta_R(\|\nabla \mathbf{u}_n\|_{\mathbb{L}^2})\left(|\mathbf{u}_n||\nabla\mathbf{u}_n|^2,|\Delta\mathbf{u}_n|\right)_{\mathbb{L}^2}\\
&\leq C_2\theta_R(\|\nabla \mathbf{u}_n\|_{\mathbb{L}^2})\||\mathbf{u}_n||\Delta\mathbf{u}_n|\|_{\mathbb{L}^2}\|\nabla\mathbf{u}_n\|_{\mathbb{L}^4}^2\\
&\leq\frac{1}{2}C_1\theta_R(\|\nabla \mathbf{u}_n\|_{\mathbb{L}^2})\||\mathbf{u}_n||\Delta\mathbf{u}_n|\|_{\mathbb{L}^2}^2+C\theta_R(\|\nabla \mathbf{u}_n\|_{\mathbb{L}^2})\|\nabla\mathbf{u}_n\|_{\mathbb{L}^4}^4.
\end{split}
\end{equation*}
This together with \eqref{8lem33} yields
\begin{equation}\label{9lem33}
\begin{split}
&\|\nabla\mathbf{u}_n(t)\|_{\mathbb{L}^2}^2+\int_0^t\|\nabla\Delta \mathbf{u}_n\|_{\mathbb{L}^2}^2\,\mathrm{d}s+\int_0^t\||\mathbf{u}_n||\nabla\mathbf{u}_n|\|_{\mathbb{L}^2}^2\,\mathrm{d}s\\\
&\leq C\|\nabla\mathbf{u}_n(0)\|_{\mathbb{L}^2}^2+C\int_0^t\|\nabla\mathbf{u}_n\|_{\mathbb{L}^2}^2\,\mathrm{d}s+C\int_0^t\|\mathbf{u}_n\|_{\mathbb{L}^2}^2\,\mathrm{d}s+C_3\int_0^t\theta_R(\|\nabla \mathbf{u}_n\|_{\mathbb{L}^2})\|\nabla\mathbf{u}_n\|_{\mathbb{L}^4}^4\,\mathrm{d}s\\
&+\sum_{j=1}^n\int_0^t\left(\nabla G_{nj}(\mathbf{u}_n),\nabla\mathbf{u}_n\right)_{\mathbb{L}^2}\,\mathrm{d}W_j(s).
\end{split}
\end{equation}
Now we estimate the term $\|\nabla\mathbf{u}_n\|_{\mathbb{L}^4}^4$ by invoking the Gagliardo-Nirenberg (GN) inequality \cite{5}. In the case of $d=1$, we have
\begin{equation}\label{ad1lem33}
\begin{split}
&\|\nabla\mathbf{u}_n\|_{\mathbb{L}^4}^4\leq C\|\mathbf{u}_n\|_{\mathbb{L}^2}^{\frac{7}{3}}\|\mathbf{u}_n\|_{\mathbb{H}^3}^{\frac{5}{3}}\leq\varepsilon\|\mathbf{u}_n\|_{\mathbb{H}^3}^2+C_{\varepsilon}\|\mathbf{u}_n\|_{\mathbb{L}^2}^{14}.
\end{split}
\end{equation}
By the standard elliptic regularity result with Neumann boundary data, it follows that
\begin{equation}\label{ad2lem33}
\begin{split}
&\|\mathbf{u}_n\|_{\mathbb{H}^3}^2\leq C\left(\|\mathbf{u}_n\|_{\mathbb{L}^2}^2+\|\nabla\mathbf{u}_n\|_{\mathbb{L}^2}^2+\|\nabla\Delta\mathbf{u}_n\|_{\mathbb{L}^2}^2\right),
\end{split}
\end{equation}
which together with \eqref{ad1lem33} implies that
\begin{equation}\label{ad3lem33}
\begin{split}
&\|\nabla\mathbf{u}_n\|_{\mathbb{L}^4}^4\leq \varepsilon\|\nabla\Delta\mathbf{u}_n\|_{\mathbb{L}^2}^2+\varepsilon\|\nabla\mathbf{u}_n\|_{\mathbb{L}^2}^2+C_{\varepsilon}\left(1+\|\mathbf{u}_n\|_{\mathbb{L}^2}^{14}\right).
\end{split}
\end{equation}
In the case of $d=2$, applying the GN inequality, Young's inequality as well as inequality \eqref{7lem33}, it follows that
\begin{equation}\label{10lem33}
\begin{split}
&\|\nabla\mathbf{u}_n\|_{\mathbb{L}^4}^4\leq C\|\nabla\mathbf{u}_n\|_{\mathbb{L}^2}^2\|\nabla\mathbf{u}_n\|_{\mathbb{H}^1}^2\leq C\left(1+\|\Delta\mathbf{u}_n\|_{\mathbb{L}^2}^2\right)\|\nabla\mathbf{u}_n\|_{\mathbb{L}^2}^2\\
&\leq C\left(1+\|\nabla \mathbf{u}_n\|_{\mathbb{L}^2}\|\nabla\Delta \mathbf{u}_n\|_{\mathbb{L}^2}\right)\|\nabla\mathbf{u}_n\|_{\mathbb{L}^2}^2\\
&\leq C\|\nabla\mathbf{u}_n\|_{\mathbb{L}^2}^2+\varepsilon\|\nabla\Delta \mathbf{u}_n\|_{\mathbb{L}^2}^2+C_{\varepsilon}\|\nabla\mathbf{u}_n\|_{\mathbb{L}^2}^6.
\end{split}
\end{equation}
Choosing $\varepsilon$ small enough such that $\varepsilon C_3\leq\frac{1}{2}$, then we see from \eqref{10lem33} that
\begin{equation}\label{11lem33}
\begin{split}
&C_3\int_0^t\theta_R(\|\nabla \mathbf{u}_n\|_{\mathbb{L}^2})\|\nabla\mathbf{u}_n\|_{\mathbb{L}^4}^4\,\mathrm{d}s\\
&\leq\frac{1}{2}\int_0^t\|\nabla\Delta \mathbf{u}_n\|_{\mathbb{L}^2}^2\,\mathrm{d}s+C\left(1+R^4\right)\int_0^t\|\nabla\mathbf{u}_n\|_{\mathbb{L}^2}^2\,\mathrm{d}s.
\end{split}
\end{equation}
Similarly, in the case of $d=3$, we have
\begin{equation}\label{12lem33}
\begin{split}
&\|\nabla\mathbf{u}_n\|_{\mathbb{L}^4}^4\leq C\|\nabla\mathbf{u}_n\|_{\mathbb{L}^2}^{\frac{5}{2}}\|\nabla\mathbf{u}_n\|_{\mathbb{H}^2}^{\frac{3}{2}}\leq C\|\nabla\mathbf{u}_n\|_{\mathbb{L}^2}^{\frac{5}{2}}\left(1+\|\nabla\mathbf{u}_n\|_{\mathbb{L}^2}^2+\|\nabla\Delta\mathbf{u}_n\|_{\mathbb{L}^2}^2\right)^{\frac{3}{4}}\\
&\leq \varepsilon\|\nabla\Delta \mathbf{u}_n\|_{\mathbb{L}^2}^2+C+C_{\varepsilon}\|\nabla\mathbf{u}_n\|_{\mathbb{L}^2}^{10}.
\end{split}
\end{equation}
Choosing $\varepsilon$ small enough such that $\varepsilon C_3\leq\frac{1}{2}$, it follows from \eqref{12lem33} that
\begin{equation}\label{13lem33}
\begin{split}
&C_3\int_0^t\theta_R(\|\nabla \mathbf{u}_n\|_{\mathbb{L}^2})\|\nabla\mathbf{u}_n\|_{\mathbb{L}^4}^4\,\mathrm{d}s\\
&\leq\frac{1}{2}\int_0^t\|\nabla\Delta \mathbf{u}_n\|_{\mathbb{L}^2}^2\,\mathrm{d}s+Ct+CR^8\int_0^t\|\nabla\mathbf{u}_n\|_{\mathbb{L}^2}^2\,\mathrm{d}s.
\end{split}
\end{equation}
Thus for $d=1,2,3$, we see from \eqref{ad3lem33}, \eqref{11lem33}, \eqref{13lem33} and \eqref{9lem33} that
\begin{equation}\label{14lem33}
\begin{split}
&\|\nabla\mathbf{u}_n(t)\|_{\mathbb{L}^2}^2+\frac{1}{2}\int_0^t\|\nabla\Delta \mathbf{u}_n\|_{\mathbb{L}^2}^2\,\mathrm{d}s+\int_0^t\||\mathbf{u}_n||\nabla\mathbf{u}_n|\|_{\mathbb{L}^2}^2\,\mathrm{d}s\\\
&\leq C+C\int_0^t1+\|\mathbf{u}_n\|_{\mathbb{L}^2}^{14}\,\mathrm{d}s+C_{R}\int_0^t\|\nabla\mathbf{u}_n\|_{\mathbb{L}^2}^2\,\mathrm{d}s+\sum_{j=1}^n\int_0^t\left(\nabla G_{nj}(\mathbf{u}_n),\nabla\mathbf{u}_n\right)_{\mathbb{L}^2}\,\mathrm{d}W_j(s).
\end{split}
\end{equation}
Using the BDG inequality, the H\"{o}lder inequality and Young's inequality, it follows that for any $p\geq1$
\begin{equation}\label{15lem33}
\begin{split}
&\mathbb{E}\sup_{s\in[0,t]}\left|\sum_{j=1}^n\int_0^s\left(\nabla G_{nj}(\mathbf{u}_n),\nabla\mathbf{u}_n\right)_{\mathbb{L}^2}\,\mathrm{d}W_j(s')\right|^p\\
&\leq C\mathbb{E}\left|\sum_{j=1}^n\int_0^t\left(\nabla G_{nj}(\mathbf{u}_n),\nabla\mathbf{u}_n\right)_{\mathbb{L}^2}^2\,\mathrm{d}s\right|^\frac{p}{2}\\
&\leq C\mathbb{E}\left|\sum_{j=1}^n\int_0^t\left(\mathbf{u}_n\times\nabla\mathbf{h}_j+\nabla\mathbf{h}_j+\nabla\Delta\mathbf{h}_j,\nabla\mathbf{u}_n\right)_{\mathbb{L}^2}^2\,\mathrm{d}s\right|^\frac{p}{2}\\
&\leq C\mathbb{E}\left|\sum_{j=1}^n\int_0^t\left(\|\nabla\mathbf{h}_j\|_{\mathbb{L}^{\infty}}\||\mathbf{u}_n||\nabla\mathbf{u}_n|\|_{\mathbb{L}^1}+\|\nabla\mathbf{h}_j\|_{\mathbb{L}^{2}}\|\nabla\mathbf{u}_n\|_{\mathbb{L}^2}+\|\nabla\Delta\mathbf{h}_j\|_{\mathbb{L}^{2}}\|\nabla\mathbf{u}_n\|_{\mathbb{L}^2}\right)^2\,\mathrm{d}s\right|^\frac{p}{2}\\
&\leq C_{\mathbf{h}}\mathbb{E}\left|\int_0^t\||\mathbf{u}_n||\nabla\mathbf{u}_n|\|_{\mathbb{L}^2}^2+\|\nabla\mathbf{u}_n\|_{\mathbb{L}^2}^2\,\mathrm{d}s\right|^\frac{p}{2}\\
&\leq\varepsilon\mathbb{E}\left(\int_0^t\||\mathbf{u}_n||\nabla\mathbf{u}_n|\|_{\mathbb{L}^2}^2\,\mathrm{d}s\right)^p+C\mathbb{E}\left(\int_0^t\|\nabla\mathbf{u}_n\|_{\mathbb{L}^2}^{2p}\,\mathrm{d}s\right)+C.
\end{split}
\end{equation}
Choosing $\varepsilon$ small enough, it follows from \eqref{1lem32}, \eqref{14lem33}, \eqref{15lem33} and Jensen's inequality that
\begin{equation*}
\begin{split}
&\mathbb{E}\sup_{s\in[0,t]}\|\nabla\mathbf{u}_n(s)\|_{\mathbb{L}^2}^{2p}+\mathbb{E}\left(\int_0^t\|\nabla\Delta\mathbf{u}_n(s)\|_{\mathbb{L}^2}^2\,\mathrm{d}s\right)^p+\mathbb{E}\left(\int_0^t\||\mathbf{u}_n||\nabla\mathbf{u}_n|\|_{\mathbb{L}^2}^2\,\mathrm{d}s\right)^p\\
&\leq C+C_{R}\int_0^t\mathbb{E}\|\nabla\mathbf{u}_n(s)\|_{\mathbb{L}^2}^{2p}\,\mathrm{d}s.
\end{split}
\end{equation*}
Thus by applying the Gronwall lemma and using the inequality \eqref{1lem32}, we obtain the estimate \eqref{1lem33}.~~$\Box$

\begin{corollary}\label{cor32} Under the same assumptions as in Lemma \ref{lem32}. Let $q\geq1$, $p>2$ and $\alpha\in(0,\frac{1}{2})$ with $p\alpha>1$. Then for any $n\in\mathbb{N}$ and every $t\in[0,T]$, there exists a positive constant $C$ independent of $n$ such that
\begin{equation}\label{2cor32}
\begin{split}
\mathbb{E}\|\mathbf{u}_n\|_{W^{\alpha,p}(0,T;(\mathbb{H}^{1})^*)}^{q}\leq C.
\end{split}
\end{equation}
\end{corollary}
\noindent\textbf{Proof.} Equation \eqref{sys2} can be written as follows:
\begin{equation}\label{3cor32}
\begin{split}
\mathbf{u}_n(t)&=\mathbf{u}_n(0)+\beta_1\int_0^tF_n^1(\mathbf{u}_n)\,\mathrm{d}s-\beta_2\int_0^tF_n^2(\mathbf{u}_n)\,\mathrm{d}s+\beta_3\int_0^t\Pi_n\mathbf{u}_n\,\mathrm{d}s-\beta_3\int_0^tF_n^3(\mathbf{u}_n)\,\mathrm{d}s\\
&-\beta_4\int_0^tF_n^4(\mathbf{u}_n)\,\mathrm{d}s+\beta_5\int_0^tF_n^5(\mathbf{u}_n)\,\mathrm{d}s-\frac{1}{2}\sum_{j=1}^n\int_0^t\Pi_n\left(G_{nj}(\mathbf{u}_n)\times\mathbf{h}_j\right)\,\mathrm{d}s\\
&+\sum_{j=1}^n\int_0^tG_{nj}(\mathbf{u}_n)\mathrm{d}W_j(t)\\
&:=\mathbf{u}_n(0)+\sum_{k=1}^7\mathbf{B}_{n,k}(\mathbf{u}_n)(t)+\mathbf{B}_{n,8}(\mathbf{u}_n,W_n)(t),~t\in[0,T].
\end{split}
\end{equation}
Let $\phi\in \mathbb{H}^1$. Then by Holder's inequality and Sobolev embedding $\mathbb{H}^1\hookrightarrow\mathbb{L}^6$ \cite{19}, we have
\begin{equation*}
\begin{split}
|\left(F_n^1(\mathbf{u}_n),\phi\right)_{\mathbb{L}^2}|&=|\left(\nabla\mathbf{u}_n,\nabla\phi\right)_{\mathbb{L}^2}|\leq\|\nabla\mathbf{u}_n\|_{\mathbb{L}^2}\|\phi\|_{\mathbb{H}^1},\\
|\left(F_n^2(\mathbf{u}_n),\phi\right)_{\mathbb{L}^2}|&=|\left(\nabla\Delta\mathbf{u}_n,\nabla\phi\right)_{\mathbb{L}^2}|\leq\|\nabla\Delta\mathbf{u}_n\|_{\mathbb{L}^2}\|\phi\|_{\mathbb{H}^1},\\
|\left(F_n^3(\mathbf{u}_n),\phi\right)_{\mathbb{L}^2}|&\leq\|\mathbf{u}_n\|_{\mathbb{L}^2}\|\mathbf{u}_n\|_{\mathbb{L}^6}^2\|\phi\|_{\mathbb{L}^6}\leq C\|\mathbf{u}_n\|_{\mathbb{L}^2}\|\mathbf{u}_n\|_{\mathbb{H}^1}^2\|\phi\|_{\mathbb{H}^1},\\
|\left(F_n^4(\mathbf{u}_n),\phi\right)_{\mathbb{L}^2}|&=|\left(\mathbf{u}_n\times\nabla\mathbf{u}_n,\nabla\phi\right)_{\mathbb{L}^2}|\leq\|\mathbf{u}_n\|_{\mathbb{L}^4}\|\nabla\mathbf{u}_n\|_{\mathbb{L}^4}\|\nabla\phi\|_{\mathbb{L}^2}\leq C\|\mathbf{u}_n\|_{\mathbb{H}^1}\|\nabla\mathbf{u}_n\|_{\mathbb{H}^1}\|\nabla\phi\|_{\mathbb{L}^2},\\
|\left(F_n^5(\mathbf{u}_n),\phi\right)_{\mathbb{L}^2}|&=|\theta_R(\|\nabla \mathbf{u}_n\|_{\mathbb{L}^2})\left(\nabla(|\mathbf{u}_n|^2\mathbf{u}_n),\nabla\phi\right)_{\mathbb{L}^2}|\leq\|\nabla(|\mathbf{u}_n|^2\mathbf{u}_n)\|_{\mathbb{L}^2}\|\nabla\phi\|_{\mathbb{L}^2}\\
&\leq3\||\mathbf{u}_n|^2|\nabla\mathbf{u}_n|\|_{\mathbb{L}^2}\|\nabla\phi\|_{\mathbb{L}^2}\leq3\|\mathbf{u}_n\|_{\mathbb{L}^6}^2\|\nabla\mathbf{u}_n\|_{\mathbb{L}^6}\|\nabla\phi\|_{\mathbb{L}^2}\leq C\|\mathbf{u}_n\|_{\mathbb{H}^1}^2\|\nabla\mathbf{u}_n\|_{\mathbb{H}^1}\|\phi\|_{\mathbb{H}^1}.
\end{split}
\end{equation*}
Similarly, through direct calculation, it follows from \eqref{condh} that
\begin{equation*}\label{9cor32}
\begin{split}
\left\|\sum_{j=1}^n\Pi_n\left(G_{nj}(\mathbf{u}_n)\times\mathbf{h}_j\right)\right\|_{(\mathbb{H}^{1})^*}\leq C_{\mathbf{h}}+C_{\mathbf{h}}\|\mathbf{u}_n\|_{\mathbb{L}^2}.
\end{split}
\end{equation*}
Thus we derive that for all $q\geq1$
\begin{equation}\label{10cor32}
\begin{split}
&\sum_{k=1}^7\mathbb{E}\|\mathbf{B}_{n,k}(\mathbf{u}_n)\|_{W^{1,2}(0,T;(\mathbb{H}^{1})^*)}^q\\
&\leq C+C\mathbb{E}\left(\int_0^T\|\mathbf{u}_n\|_{\mathbb{H}^{1}}^2\,\mathrm{d}s\right)^{\frac{q}{2}}+C\mathbb{E}\left(\int_0^T\|\nabla\Delta\mathbf{u}_n\|_{\mathbb{L}^{2}}^2\,\mathrm{d}s\right)^{\frac{q}{2}}+C\mathbb{E}\left(\int_0^T\|\mathbf{u}_n\|_{\mathbb{L}^{2}}^2\|\mathbf{u}_n\|_{\mathbb{H}^{1}}^4\,\mathrm{d}s\right)^{\frac{q}{2}}\\
&+C\mathbb{E}\left(\int_0^T\|\mathbf{u}_n\|_{\mathbb{H}^{1}}^2\|\nabla\mathbf{u}_n\|_{\mathbb{H}^{1}}^2\,\mathrm{d}s\right)^{\frac{q}{2}}+C\mathbb{E}\left(\int_0^T\|\mathbf{u}_n\|_{\mathbb{H}^{1}}^4\|\nabla\mathbf{u}_n\|_{\mathbb{H}^{1}}^2\,\mathrm{d}s\right)^{\frac{q}{2}}\\
&\leq C+C\mathbb{E}\left[\sup_{t\in[0,T]}\|\mathbf{u}_n(t)\|_{\mathbb{H}^{1}}^{2q}\left(\int_0^T\|\mathbf{u}_n\|_{\mathbb{H}^{2}}^2\,\mathrm{d}s\right)^{\frac{q}{2}}\right]\\
&\leq C+C\mathbb{E}\sup_{t\in[0,T]}\|\mathbf{u}_n(t)\|_{\mathbb{H}^{1}}^{4q}+C\mathbb{E}\left(\int_0^T\|\mathbf{u}_n\|_{\mathbb{H}^{2}}^2\,\mathrm{d}s\right)^{q}\leq C.
\end{split}
\end{equation}
To estimate the last term on the right-hand side of \eqref{3cor32}, we use Lemma \ref{lem92} to derive that
\begin{equation}\label{11cor32}
\begin{split}
&\mathbb{E}\|\mathbf{B}_{n,8}(\mathbf{u}_n,W_n)\|_{W^{\alpha,p}(0,T;(\mathbb{H}^{1})^*)}^q\leq C\mathbb{E}\left[\int_0^T\left(\sum_{j=1}^n\|G_{nj}(\mathbf{u}_n)\|_{(\mathbb{H}^{1})^*}^2\right)^{\frac{q}{2}}\,\mathrm{d}s\right]\\
&\leq C_{\mathbf{h}}+C_{\mathbf{h}}\mathbb{E}\left(\int_0^T\|\mathbf{u}_n\|_{\mathbb{L}^2}^q\,\mathrm{d}s\right)\leq C.
\end{split}
\end{equation}
Since $
W^{1,2}(0,T;(\mathbb{H}^{1})^*)\hookrightarrow W^{\alpha,p}(0,T;(\mathbb{H}^{1})^*),~\textrm{if}~\frac{1}{2}+\frac{1}{p}>\alpha,
$
the estimate \eqref{2cor32} then follows from \eqref{10cor32} and \eqref{11cor32}.~~$\Box$

\section{Tightness result}\label{sec3}
\begin{lemma}\label{lem41}
If $\beta>1$ and $p>1$, then the measures $\{\mathcal{L}(\mathbf{u}_n)\}_{n\in\mathbb{N}}$ on $\mathcal{X}:=C([0,T];(\mathbb{H}^{\beta})^*)\cap L^2(0,T;\mathbb{H}^2)\cap L^p(0,T;\mathbb{L}^4)$ are tight.
\end{lemma}
\noindent\textbf{Proof.} According to the uniform bounds provided by Lemma \ref{lem33} and Corollary \ref{cor32}, we have
\begin{equation*}
\begin{split}
\mathbb{E}\|\mathbf{u}_n\|_{W^{\alpha,p}(0,T;(\mathbb{H}^{1})^*)\cap L^p(0,T;\mathbb{H}^1)\cap L^2(0,T;\mathbb{H}^3)}\leq C.
\end{split}
\end{equation*}
Thanks to Lemma \ref{lem93} and Lemma \ref{lem94}, we obtain the following compact embeddings
\begin{equation*}
\begin{split}
W^{\alpha,p}(0,T;(\mathbb{H}^{1})^*)\cap L^p(0,T;\mathbb{H}^1)\hookrightarrow C([0,T];(\mathbb{H}^{\beta})^*)\cap L^p(0,T;\mathbb{L}^4),\\
W^{\alpha,p}(0,T;(\mathbb{H}^{1})^*)\cap L^2(0,T;\mathbb{H}^3)\hookrightarrow C([0,T];(\mathbb{H}^{\beta})^*)\cap L^2(0,T;\mathbb{H}^2),
\end{split}
\end{equation*}
which implies the tightness of $\{\mathcal{L}(\mathbf{u}_n)\}_{n\in\mathbb{N}}$.~~$\Box$

By Lemma \ref{lem41} and Prokhorov theorem, the Borel subsets of $C([0,T];S_n)$ are Borel subsets of $\mathcal{X}$. Moreover, noticing that $\mathcal{X}$ is a separable metric space, we can apply the Skorohod theorem \cite{27} to obtain the following results.
\begin{proposition}\label{pro41}
There exist
\begin{enumerate}
\item [(1)] a probability space $(\Omega',\mathcal{F}',\mathbb{P}')$,
\item [(2)] a sequence $\{(\mathbf{u}_n',W_n')\}$ of random variables defined on $(\Omega',\mathcal{F}',\mathbb{P}')$ and taking values in the space $\mathcal{X}\times C([0,T];\mathbb{R}^{\infty})$,
\item [(3)] a random variable $(\mathbf{u}',W')$ defined on $(\Omega',\mathcal{F}',\mathbb{P}')$ and taking values in the space $\mathcal{X}\times C([0,T];\mathbb{R}^{\infty})$,
\end{enumerate}
such that in the space $\mathcal{X}\times C([0,T];\mathbb{R}^{\infty})$ there hold
\begin{enumerate}
\item [(a)] $\mathcal{L}(\mathbf{u}_n,W_n)=\mathcal{L}(\mathbf{u}_n',W_n')$,
\item [(b)] $(\mathbf{u}_n',W_n')$ converges $\mathbb{P}'$-almost surely to $(\mathbf{u}',W')$ in the topology of $\mathcal{X}$.
\end{enumerate}
\end{proposition}

Since the laws on $C([0,T];S_n)$ of $\mathbf{u}_n$ and $\mathbf{u}_n'$ are equal, we have the following estimate.
\begin{corollary}\label{cor41}
For every $p\geq1$ and $T>0$,
\begin{equation}\label{1cor41}
\begin{split}
\sup_{n\in\mathbb{N}}\mathbb{E}\sup_{t\in[0,T]}\|\mathbf{u}_n'(t)\|_{\mathbb{H}^1}^{2p}+\mathbb{E}\left(\int_0^T\|\mathbf{u}_n'(t)\|_{\mathbb{H}^3}^2\,\mathrm{d}t\right)^p<\infty.
\end{split}
\end{equation}
\end{corollary}

Moreover, we have the following weak convergence result.
\begin{lemma}\label{lem54} For every $p\geq1$ there holds
\begin{equation*}
\begin{split}
\mathbf{u}_n'\rightarrow\mathbf{u}'~\textrm{weakly~in}~L^{2p}(\Omega';L^{\infty}(0,T;\mathbb{H}^1)\cap L^2(0,T;\mathbb{H}^3)).
\end{split}
\end{equation*}
Moreover, the process $\mathbf{u}'\in L^{2p}(\Omega';C([0,T];\mathbb{H}^1_w))$.
\end{lemma}
\noindent\textbf{Proof.} We shall fist prove that
\begin{equation}\label{2lem5-5}
\begin{split}
\mathbf{u}_n'\rightarrow\mathbf{u}'~\textrm{weakly in}~L^{\frac{4}{3}}(\Omega';L^{4}(0,T;\mathbb{L}^2)).
\end{split}
\end{equation}
Since $\mathbf{u}_n'\rightarrow\mathbf{u}'~\textrm{weakly in}~L^{4}(0,T;\mathbb{L}^4)\hookrightarrow L^{4}(0,T;\mathbb{L}^2)$, $\mathbb{P}'$-a.s., for any $\phi\in L^{4}(\Omega';L^{\frac{4}{3}}(0,T;\mathbb{L}^2))$ there holds
$
\int_0^T(\mathbf{u}_n',\phi)_{\mathbb{L}^2}\,\mathrm{d}t\rightarrow\int_0^T(\mathbf{u}',\phi)_{\mathbb{L}^2}\,\mathrm{d}t.
$
Moreover, by using inequality \eqref{1cor41}, we have
\begin{equation*}
\begin{split}
&\sup_{n\in\mathbb{N}}\mathbb{E}'\left|\int_0^T(\mathbf{u}_n',\phi)_{\mathbb{L}^2}\,\mathrm{d}t\right|^2\leq\sup_{n\in\mathbb{N}}\mathbb{E}'\left(\|\mathbf{u}_n'\|_{L^{\infty}(0,T;\mathbb{L}^2)}^2\|\phi\|_{L^1(0,T;\mathbb{L}^2)}^2\right)\\
&\leq\sup_{n\in\mathbb{N}}\|\mathbf{u}_n'\|_{L^{4}(\Omega';L^{\infty}(0,T;\mathbb{L}^2))}^2\|\phi\|_{L^{4}(\Omega';L^{1}(0,T;\mathbb{L}^2))}^2<\infty.
\end{split}
\end{equation*}
Thus by using the Vitali convergence theorem we have
$
\mathbb{E}'\int_0^T(\mathbf{u}_n',\phi)_{\mathbb{L}^2}\,\mathrm{d}t\rightarrow\mathbb{E}'\int_0^T(\mathbf{u}',\phi)_{\mathbb{L}^2}\,\mathrm{d}t,
$
which means the result \eqref{2lem5-5}. On the other hand, by using the Banach-Alaoglu theorem we infer from \eqref{1cor41} that there exists a subsequence of $\{\mathbf{u}_n'\}$ (still denoted by $\{\mathbf{u}_n'\}$) and $\mathbf{v}\in L^{2p}(\Omega';L^{\infty}(0,T;\mathbb{H}^1)\cap L^{2}(0,T;\mathbb{H}^3))$ such that
$
\mathbf{u}_n'\rightarrow\mathbf{v}~\textrm{weakly in}~L^{2p}(\Omega';L^{\infty}(0,T;\mathbb{H}^1)\cap L^{2}(0,T;\mathbb{H}^3))\subset L^{\frac{4}{3}}(\Omega';L^{4}(0,T;\mathbb{L}^2)).
$
By the uniqueness of weak limit, we infer that
$
\mathbf{u}'=\mathbf{v}~\textrm{in}~L^{2p}(\Omega';L^{\infty}(0,T;\mathbb{H}^1)\cap L^{2}(0,T;\mathbb{H}^3)).
$
Moreover by Theorem 2.1 of \cite{40}, for any $\beta\geq0$ the following imbedding is continuous
\begin{equation*}\label{3lem54}
\begin{split}
L^{\infty}(0,T;\mathbb{H}^1)\cap C([0,T];(\mathbb{H}^{\beta})^*)\subset C([0,T];\mathbb{H}_w^1).
\end{split}
\end{equation*}
This together with the weakly convergence of $\mathbf{u}_n$ to $\mathbf{u}$ in $L^{2p}(\Omega';L^{\infty}(0,T;\mathbb{H}^1))$ and the completeness of $C([0,T];\mathbb{H}_w^1)$ imply that $\mathbf{u}'\in L^{2p}(\Omega;C([0,T];\mathbb{H}^1_w))$. The proof is thus completed.~~$\Box$

\section{Local pathwise weak solution}\label{sec4}
This section is devoted to show that $\mathbf{u}'$ provided by Proposition \ref{pro41} is a local martingale weak solution of \eqref{sys1}.

Thanks to the conclusions provided by Proposition \ref{pro41}, Corollary \ref{cor41} and Lemma \ref{lem54}, we have the following convergence results.
\begin{lemma}\label{lem53}
For any $\phi\in \mathbb{H}^1$,
\begin{equation*}\label{3lem55}
\begin{split}
\lim_{n\rightarrow\infty}\mathbb{E}\int_0^t\left(\Delta \mathbf{u}_n',\phi\right)_{\mathbb{L}^2}\,\mathrm{d}s=-\mathbb{E}\int_0^t\left(\nabla\mathbf{u}',\nabla\phi\right)_{\mathbb{L}^2}\,\mathrm{d}s,
\end{split}
\end{equation*}
\begin{equation*}\label{4lem55}
\begin{split}
\lim_{n\rightarrow\infty}\mathbb{E}\int_0^t\left(\Delta^2 \mathbf{u}_n',\phi\right)_{\mathbb{L}^2}\,\mathrm{d}s=-\mathbb{E}\int_0^t\left(\nabla\Delta\mathbf{u}',\nabla\phi\right)_{\mathbb{L}^2}\,\mathrm{d}s,
\end{split}
\end{equation*}
\begin{equation*}\label{1lem53}
\begin{split}
\lim_{n\rightarrow\infty}\mathbb{E}\int_0^t\left(\Pi_n\left(\left(1-|\mathbf{u}_n'|^2\right)\mathbf{u}_n'\right),\phi\right)_{\mathbb{L}^2}\,\mathrm{d}s=\mathbb{E}\int_0^t\left(\left(1-|\mathbf{u}'|^2\right)\mathbf{u}',\phi\right)_{\mathbb{L}^2}\,\mathrm{d}s,
\end{split}
\end{equation*}
\begin{equation*}\label{2lem53}
\begin{split}
\lim_{n\rightarrow\infty}\mathbb{E}\int_0^t\left(\Pi_n\left(\left(\mathbf{u}_n'\times\mathbf{h}_j\right)\times\mathbf{h}_j\right),\phi\right)_{\mathbb{L}^2}\,\mathrm{d}s=\mathbb{E}\int_0^t\left(\left(\left(\mathbf{u}'\times\mathbf{h}_j\right)\times\mathbf{h}_j\right),\phi\right)_{\mathbb{L}^2}\,\mathrm{d}s,
\end{split}
\end{equation*}
\begin{equation}\label{1lem55}
\begin{split}
\lim_{n\rightarrow\infty}\mathbb{E}\int_0^t\left(\Pi_n\left(\mathbf{u}_n'\times\Delta \mathbf{u}_n'\right),\phi\right)_{\mathbb{L}^2}\,\mathrm{d}s=-\mathbb{E}\int_0^t\left(\mathbf{u}'\times\nabla\mathbf{u}',\nabla\phi\right)_{\mathbb{L}^2}\,\mathrm{d}s,
\end{split}
\end{equation}
\begin{equation}\label{2lem55}
\begin{split}
\lim_{n\rightarrow\infty}\mathbb{E}\int_0^t\theta_R(\|\nabla \mathbf{u}_n'\|_{\mathbb{L}^2})\left(\Pi_n\left(\Delta\left(|\mathbf{u}_n'|^2\mathbf{u}_n'\right)\right),\phi\right)_{\mathbb{L}^2}\,\mathrm{d}s=-\mathbb{E}\int_0^t\theta_R(\|\nabla \mathbf{u}'\|_{\mathbb{L}^2})\left(\nabla\left(|\mathbf{u}'|^2\mathbf{u}'\right),\nabla\phi\right)_{\mathbb{L}^2}\,\mathrm{d}s
\end{split}
\end{equation}
\end{lemma}
\noindent\textbf{Proof.} We shall present the proofs of \eqref{1lem55} and \eqref{2lem55}, as the others are similar. We first prove that for $t\in[0,T]$ and $\mathbb{P}'$-a.s.,
\begin{equation*}
\begin{split}
\lim_{n\rightarrow\infty}\int_0^t\left(\Pi_n\left(\mathbf{u}_n'\times\Delta \mathbf{u}_n'\right),\phi\right)_{\mathbb{L}^2}\,\mathrm{d}s=-\int_0^t\left(\mathbf{u}'\times\nabla\mathbf{u}',\nabla\phi\right)_{\mathbb{L}^2}\,\mathrm{d}s.
\end{split}
\end{equation*}
By using integration by parts, it is sufficient to prove that $\mathbb{P}'$-a.s.,
\begin{equation*}
\begin{split}
\lim_{n\rightarrow\infty}\int_0^T\left(\mathbf{u}'_n\times\nabla\mathbf{u}'_n,\nabla\phi\right)_{\mathbb{L}^2}\,\mathrm{d}s=\int_0^T\left(\mathbf{u}'\times\nabla\mathbf{u}',\nabla\phi\right)_{\mathbb{L}^2}\,\mathrm{d}s.
\end{split}
\end{equation*}
By using Corollary \ref{cor41} and Proposition \ref{pro41}, it follows that $\mathbb{P}'$-a.s.,
\begin{equation*}
\begin{split}
&\left|\int_0^T\left(\mathbf{u}'_n\times\nabla\mathbf{u}'_n,\nabla\phi\right)_{\mathbb{L}^2}\,\mathrm{d}s-\int_0^T\left(\mathbf{u}'\times\nabla\mathbf{u}',\nabla\phi\right)_{\mathbb{L}^2}\,\mathrm{d}s\right|\\
&\leq\left|\int_0^T\left(\mathbf{u}'_n\times\nabla\left(\mathbf{u}'_n-\mathbf{u}'\right),\nabla\phi\right)_{\mathbb{L}^2}\,\mathrm{d}s\right|+\left|\int_0^T\left(\left(\mathbf{u}'_n-\mathbf{u}'\right)\times\nabla\mathbf{u}',\nabla\phi\right)_{\mathbb{L}^2}\,\mathrm{d}s\right|\\
&\leq\|\nabla\phi\|_{\mathbb{L}^4}\|\mathbf{u}_n'\|_{L^2(0,T;\mathbb{L}^4)}\|\mathbf{u}_n'-\mathbf{u}'\|_{L^2(0,T;\mathbb{H}^1)}+\|\nabla\phi\|_{\mathbb{L}^4}\|\nabla\mathbf{u}'\|_{L^2(0,T;\mathbb{L}^4)}\|\mathbf{u}_n'-\mathbf{u}'\|_{L^2(0,T;\mathbb{L}^2)}\\
&\rightarrow0~\textrm{as}~n\rightarrow\infty.
\end{split}
\end{equation*}
Moreover, since $\mathbb{H}^1\hookrightarrow\mathbb{L}^4$,
\begin{equation*}
\begin{split}
&\sup_{n\in\mathbb{N}}\mathbb{E}\left|\int_0^T\left(\mathbf{u}'_n\times\nabla\mathbf{u}'_n,\nabla\phi\right)_{\mathbb{L}^2}\,\mathrm{d}s\right|^2\leq \sup_{n\in\mathbb{N}}C\|\nabla\phi\|_{\mathbb{L}^4}^2\mathbb{E}\int_0^T\|\nabla\mathbf{u}'_n\|_{\mathbb{L}^2}^2\|\mathbf{u}'_n\|_{\mathbb{L}^4}^2\,\mathrm{d}s\\
&\leq\sup_{n\in\mathbb{N}}C\|\nabla\phi\|_{\mathbb{L}^4}^2\mathbb{E}\sup_{s\in[0,T]}\|\mathbf{u}'_n(s)\|_{\mathbb{H}^1}^4<\infty.
\end{split}
\end{equation*}
Thus the Vitali convergence theorem yields \eqref{1lem55}.

Similarly, to show \eqref{2lem55}, it is sufficient to prove that
\begin{equation*}
\begin{split}
&\lim_{n\rightarrow\infty}\mathbb{E}\int_0^t\theta_R(\|\nabla \mathbf{u}_n'\|_{\mathbb{L}^2})\left(\nabla\left(|\mathbf{u}_n'|^2\mathbf{u}_n'\right),\nabla\phi\right)_{\mathbb{L}^2}\,\mathrm{d}s=\mathbb{E}\int_0^t\theta_R(\|\nabla \mathbf{u}'\|_{\mathbb{L}^2})\left(\nabla\left(|\mathbf{u}'|^2\mathbf{u}'\right),\nabla\phi\right)_{\mathbb{L}^2}\,\mathrm{d}s.
\end{split}
\end{equation*}
By using the triangle inequality, Corollary \ref{cor41}, Proposition \ref{pro41} and Lemma \ref{lem54}, we infer that
\begin{equation*}
\begin{split}
&\left|\int_0^T\theta_R(\|\nabla \mathbf{u}_n'\|_{\mathbb{L}^2})\left(\nabla\left(|\mathbf{u}_n'|^2\mathbf{u}_n'\right),\nabla\phi\right)_{\mathbb{L}^2}\,\mathrm{d}s-\int_0^T\theta_R(\|\nabla \mathbf{u}'\|_{\mathbb{L}^2})\left(\nabla\left(|\mathbf{u}'|^2\mathbf{u}'\right),\nabla\phi\right)_{\mathbb{L}^2}\,\mathrm{d}s\right|\\
&\leq\left|\sup_{s\in[0,T]}\theta_R(\|\nabla \mathbf{u}_n'(s)\|_{\mathbb{L}^2})\right|\left|\int_0^T\left(\nabla\left(|\mathbf{u}_n'|^2\mathbf{u}_n'\right),\nabla\phi\right)_{\mathbb{L}^2}\,\mathrm{d}s-\int_0^T\left(\nabla\left(|\mathbf{u}'|^2\mathbf{u}'\right),\nabla\phi\right)_{\mathbb{L}^2}\,\mathrm{d}s\right|\\
&+\int_0^T\left|\theta_R(\|\nabla \mathbf{u}_n'\|_{\mathbb{L}^2})-\theta_R(\|\nabla \mathbf{u}'\|_{\mathbb{L}^2})\right|\left|\left(\nabla\left(|\mathbf{u}'|^2\mathbf{u}'\right),\nabla\phi\right)_{\mathbb{L}^2}\right|\,\mathrm{d}s\\
&\leq\left|\int_0^T\left(\nabla\left(|\mathbf{u}_n'|^2\mathbf{u}_n'\right)-\nabla\left(|\mathbf{u}_n'|^2\mathbf{u}'\right),\nabla\phi\right)_{\mathbb{L}^2}\,\mathrm{d}s\right|+\left|\int_0^T\left(\nabla\left(|\mathbf{u}_n'|^2\mathbf{u}'\right)-\nabla\left(|\mathbf{u}'|^2\mathbf{u}'\right),\nabla\phi\right)_{\mathbb{L}^2}\,\mathrm{d}s\right|\\
&+C\int_0^T\|\nabla\left(\mathbf{u}_n'-\mathbf{u}'\right)\|_{\mathbb{L}^2}\left|\left(\nabla\left(|\mathbf{u}'|^2\mathbf{u}'\right),\nabla\phi\right)_{\mathbb{L}^2}\right|\,\mathrm{d}s\\
&\leq C\|\nabla\phi\|_{\mathbb{L}^2}\int_0^T\|\mathbf{u}_n'\|_{\mathbb{L}^6}\|\nabla\mathbf{u}_n'\|_{\mathbb{L}^6}\|\mathbf{u}_n'-\mathbf{u}'\|_{\mathbb{L}^6}\,\mathrm{d}s+C\|\nabla\phi\|_{\mathbb{L}^2}\int_0^T\|\mathbf{u}_n'\|_{\mathbb{L}^6}^2\|\nabla\left(\mathbf{u}_n'-\mathbf{u}'\right)\|_{\mathbb{L}^6}\,\mathrm{d}s\\
&+C\|\nabla\phi\|_{\mathbb{L}^2}\int_0^T\|\mathbf{u}_n'\|_{\mathbb{L}^6}\|\mathbf{u}'\|_{\mathbb{L}^6}\|\nabla(\mathbf{u}_n'-\mathbf{u}')\|_{\mathbb{L}^6}\,\mathrm{d}s+C\|\nabla\phi\|_{\mathbb{L}^2}\int_0^T\|\mathbf{u}'\|_{\mathbb{L}^6}\|\nabla\mathbf{u}'\|_{\mathbb{L}^6}\|\mathbf{u}_n'-\mathbf{u}'\|_{\mathbb{L}^6}\,\mathrm{d}s\\
&+C\|\nabla\phi\|_{\mathbb{L}^2}\int_0^T\|\nabla\mathbf{u}'\|_{\mathbb{L}^6}\||\mathbf{u}_n'|+|\mathbf{u}'|\|_{\mathbb{L}^6}\|\mathbf{u}_n'-\mathbf{u}'\|_{\mathbb{L}^6}\,\mathrm{d}s+C\|\nabla\phi\|_{\mathbb{L}^2}\int_0^T\|\nabla\mathbf{u}'\|_{\mathbb{L}^6}\|\mathbf{u}'\|_{\mathbb{L}^6}^2\|\nabla(\mathbf{u}_n'-\mathbf{u}')\|_{\mathbb{L}^2}\,\mathrm{d}s\\
&\leq C\|\phi\|_{\mathbb{H}^1}\|\mathbf{u}_n'\|_{L^{\infty}(0,T;\mathbb{H}^1)}\|\mathbf{u}_n'\|_{L^{2}(0,T;\mathbb{H}^2)}\|\mathbf{u}_n'-\mathbf{u}'\|_{L^{2}(0,T;\mathbb{H}^2)}+C\|\phi\|_{\mathbb{H}^1}\|\mathbf{u}_n'\|_{L^{\infty}(0,T;\mathbb{H}^1)}^2\|\mathbf{u}_n'-\mathbf{u}'\|_{L^{2}(0,T;\mathbb{H}^2)}\\
&+C\|\phi\|_{\mathbb{H}^1}\|\mathbf{u}_n'\|_{L^{\infty}(0,T;\mathbb{H}^1)}\|\mathbf{u}'\|_{L^{\infty}(0,T;\mathbb{H}^1)}\|\mathbf{u}_n'-\mathbf{u}'\|_{L^{2}(0,T;\mathbb{H}^2)}\\
&+C\|\phi\|_{\mathbb{H}^1}\|\mathbf{u}_n'\|_{L^{\infty}(0,T;\mathbb{H}^1)}\|\mathbf{u}'\|_{L^{2}(0,T;\mathbb{H}^2)}\|\mathbf{u}_n'-\mathbf{u}'\|_{L^{2}(0,T;\mathbb{H}^1)}\\
&+C\|\phi\|_{\mathbb{H}^1}\left(\|\mathbf{u}_n'\|_{L^{\infty}(0,T;\mathbb{H}^1)}+\|\mathbf{u}'\|_{L^{\infty}(0,T;\mathbb{H}^1)}\right)\|\mathbf{u}'\|_{L^{\infty}(0,T;\mathbb{H}^2)}\|\mathbf{u}_n'-\mathbf{u}'\|_{L^{2}(0,T;\mathbb{H}^1)}\\
&+C\|\phi\|_{\mathbb{H}^1}\|\mathbf{u}'\|_{L^{\infty}(0,T;\mathbb{H}^1)}^2\|\mathbf{u}'\|_{L^{2}(0,T;\mathbb{H}^2)}\|\mathbf{u}_n'-\mathbf{u}'\|_{L^{2}(0,T;\mathbb{H}^1)}\rightarrow0~\textrm{as}~n\rightarrow\infty.
\end{split}
\end{equation*}
Moreover, the sequence $\int_0^T\theta_R(\|\nabla \mathbf{u}_n'\|_{\mathbb{L}^2})\left(\nabla\left(|\mathbf{u}_n'|^2\mathbf{u}_n'\right),\nabla\Pi_n\phi\right)_{\mathbb{L}^2}\,\mathrm{d}s$ is uniformly integrable on $\Omega'$. Indeed,
\begin{equation*}
\begin{split}
&\sup_{n\in\mathbb{N}}\mathbb{E}\left|\int_0^T\theta_R(\|\nabla \mathbf{u}_n'\|_{\mathbb{L}^2})\left(\nabla\left(|\mathbf{u}_n'|^2\mathbf{u}_n'\right),\nabla\phi\right)_{\mathbb{L}^2}\,\mathrm{d}s\right|^2\\
&\leq \sup_{n\in\mathbb{N}}C\|\nabla\phi\|_{\mathbb{L}^2}^2\mathbb{E}\left|\int_0^T\|\nabla \mathbf{u}_n'\|_{\mathbb{L}^6}\|\mathbf{u}_n'\|_{\mathbb{L}^6}^2\,\mathrm{d}s\right|^2\\
&\leq \sup_{n\in\mathbb{N}}C\|\phi\|_{\mathbb{H}^1}^2\left(\mathbb{E}\sup_{s\in[0,T]}\|\mathbf{u}_n'(s)\|_{\mathbb{H}^1}^4\right)^{\frac{1}{2}}\left(\mathbb{E}\int_0^T\|\mathbf{u}_n'\|_{\mathbb{H}^2}^2\,\mathrm{d}s\right)^{\frac{1}{2}}<\infty.
\end{split}
\end{equation*}
This together with Vitali convergence theorem yields \eqref{2lem55}.~~$\Box$

Let $(\mathbf{M}_n'(t))$ be a sequence of stochastic process on $(\Omega',\mathcal{F}',\mathbb{P}')$ defined by
$$
\mathbf{M}_n'(t):=\mathbf{u}_n'(t)-\mathbf{u}_n'(0)-\sum_{k=1}^7\mathbf{B}_{n,k}(\mathbf{u}_n')(t),
$$
where $\mathbf{B}_{n,k}$ is defined in \eqref{3cor32}. Let us also define
\begin{equation*}
\begin{split}
\mathbf{M}'(t)&:=\mathbf{u}'(t)-\mathbf{u}'(0)-\beta_1\int_0^t\Delta \mathbf{u}'(s)\,\mathrm{d}s+\beta_2\int_0^t\Delta^2 \mathbf{u}'(s)\,\mathrm{d}s-\beta_3\int_0^t\left(1-|\mathbf{u}'(s)|^2\mathbf{u}'(s)\right)\,\mathrm{d}s\\
&+\beta_4\int_0^t\mathbf{u}'(s)\times\Delta \mathbf{u}'(s)\,\mathrm{d}s-\beta_5\int_0^t\theta_R(\|\nabla \mathbf{u}'(s)\|_{\mathbb{L}^2})\Delta \left(|\mathbf{u}'(s)|^2\mathbf{u}'(s)\right)\,\mathrm{d}s\\
&-\frac{1}{2}\sum_{j=1}^{\infty}\int_0^t\left(-\mathbf{u}'(s)\times\mathbf{h}_j+\mathbf{h}_j-\Delta\mathbf{h}_j\right)\times\mathbf{h}_j\,\mathrm{d}s.
\end{split}
\end{equation*}
According to Lemma \ref{lem53}, it is clear that for each $t\in[0,T]$ and $\phi\in\mathbb{H}^1$,
\begin{equation}\label{1cor555-1}
\begin{split}
\lim_{n\rightarrow\infty}\mathbb{E}'\left|\left\langle \mathbf{M}_n'(t)-\mathbf{M}'(t),\phi\right\rangle_{(\mathbb{H}^{1})^*,\mathbb{H}^1}\right|=0.
\end{split}
\end{equation}
Moreover by a standard argument as shown in Lemma 5.2 of \cite{8}, we have
\begin{equation}\label{2cor555-1}
\begin{split}
\lim_{n\rightarrow\infty}\mathbb{E}'\left|\left\langle \mathbf{M}_n'(t)-\sum_{j=1}^n\int_0^t-\mathbf{u}'\times\mathbf{h}_j+\mathbf{h}_j-\Delta\mathbf{h}_j\,\mathrm{d}W_{j}'(s),\phi\right\rangle_{(\mathbb{H}^{1})^*,\mathbb{H}^1}\right|=0.
\end{split}
\end{equation}
Therefore by \eqref{1cor555-1} and \eqref{2cor555-1}, it follows that for $t\in[0,T]$ and $\mathbb{P}'$-a.s.,
\begin{equation*}\label{equ51}
\begin{split}
\mathbf{M}'(t)=\sum_{j=1}^{\infty}\int_0^{t}\left(-\mathbf{u}'\times \mathbf{h}_j+\mathbf{h}_j-\Delta\mathbf{h}_j,\phi\right)_{\mathbb{L}^2}\mathrm{d}W_j'(s)~\textrm{in}~(\mathbb{H}^1)^*.
 \end{split}
\end{equation*}
Thus $\left(\Omega',\mathcal{F}',\mathbb{P}',W',\mathbf{u}'\right)$ is a global martingale weak solution for the following equation
\begin{equation}\label{sys1-2}
\left\{
\begin{aligned}
&\mathrm{d}\mathbf{u}=\bigl[\beta_1\Delta \mathbf{u}-\beta_2\Delta^2\mathbf{u}+\beta_3(1-|\mathbf{u}|^2)\mathbf{u}-\beta_4\mathbf{u}\times\Delta\mathbf{u}\\
&+\beta_5\theta_R(\|\nabla \mathbf{u}\|_{\mathbb{L}^2})\Delta(|\mathbf{u}|^2\mathbf{u})\bigl]\,\mathrm{d}t+\sum_{j=1}^{\infty}(-\mathbf{u}\times \mathbf{h}_j+\mathbf{h}_j-\Delta\mathbf{h}_j)\circ\mathrm{d}W_j(t),&&\textrm{in}~(0,\infty)\times\mathcal{O},\\
&\mathbf{u}(0,\cdot)=\mathbf{u}_0,&&\textrm{in}~\mathcal{O},\\
&\frac{\partial\mathbf{u}}{\partial\mathbf{n}}=\frac{\partial\Delta\mathbf{u}}{\partial\mathbf{n}}=0,&&\textrm{on}~(0,\infty)\times\partial\mathcal{O}.
\end{aligned}
\right.
\end{equation}

Next we shall prove that this martingale solution is actually the pathwise solution of \eqref{sys1-2}.

\begin{proposition}\label{pro51} Let $\mathcal{O}\subset\mathbb{R}^d,~d=1,2,3$ and let $\mathbf{u}_0\in \mathbb{H}^1$ be fixed. Assume that $\left(\Omega,\mathcal{F},\mathbb{P},W,\mathbf{u}_1\right)$ and $\left(\Omega,\mathcal{F},\mathbb{P},W,\mathbf{u}_2\right)$ are two martingale weak solution of \eqref{sys1-2} such that for $i=1,2$,
$$\mathbf{u}_i(0)=\mathbf{u}(0);~\mathbf{u}_i\in L^{\infty}(0,T;\mathbb{H}^1)\cap L^2(0,T;\mathbb{H}^3);~\mathbf{u}_i~\textrm{satisfies equation}~ \eqref{sys1-2}.
$$
Then $\mathbf{u}_1=\mathbf{u}_2$, $\mathbb{P}$-a.s.
\end{proposition}
\noindent\textbf{Proof.} Let $\mathbf{u}^*:=\mathbf{u}_1-\mathbf{u}_2$. Then $\mathbf{u}^*$ satisfies the following equation
\begin{equation*}
\begin{split}
\mathrm{d}\mathbf{u}^*&=\Bigl[\beta_1\Delta \mathbf{u}^*-\beta_2\Delta^2\mathbf{u}^*+\beta_3\left(\mathbf{u}^*-|\mathbf{u}_1|^2\mathbf{u}^*+\mathbf{u}_2\left(|\mathbf{u}_2|^2-|\mathbf{u}_1|^2\right)\right)-\beta_4\left(\mathbf{u}_1\times\Delta\mathbf{u}^*+\mathbf{u}^*\times\Delta\mathbf{u}_2\right)\\
&+\beta_5\theta_R(\|\nabla \mathbf{u}_1\|_{\mathbb{L}^2})\Delta\left(|\mathbf{u}_1|^2\mathbf{u}^*+\left(|\mathbf{u}_1|^2-|\mathbf{u}_2|^2\right)\mathbf{u}_2\right)+\beta_5\left(\theta_R(\|\nabla \mathbf{u}_1\|_{\mathbb{L}^2})-\theta_R(\|\nabla \mathbf{u}_2\|_{\mathbb{L}^2})\right)\Delta\left(|\mathbf{u}_2|^2\mathbf{u}_2\right)\\
&-\frac{1}{2}\sum_{j=1}^{\infty}(\mathbf{u}^*\times \mathbf{h}_j)\times\mathbf{h}_j\Bigl]\,\mathrm{d}t+\sum_{j=1}^{\infty}-\mathbf{u}^*\times \mathbf{h}_j\,\mathrm{d}W_j(t),
\end{split}
\end{equation*}
in $(H^1)^*$ with $\mathbf{u}^*_0=0$. Let
\begin{equation*}
\begin{split}
\xi^K:=\inf\left\{t\geq0:\|\mathbf{u}_1(t)\|_{\mathbb{H}^1}^2+\|\mathbf{u}_2(t)\|_{\mathbb{H}^1}^2+\int_0^t\|\mathbf{u}_1\|_{\mathbb{H}^3}^2\,\mathrm{d}s+\int_0^t\|\mathbf{u}_2\|_{\mathbb{H}^3}^2\,\mathrm{d}s>K\right\}\wedge T,~K>0.
\end{split}
\end{equation*}
Due to Lemma \ref{lem54}, it is easy to check that
$
\xi^K\nearrow T~\textrm{as}~K\rightarrow\infty,~\mathbb{P}\textrm{-a.s.}
$
By using the It\^{o} formula to $\|\mathbf{u}^*(t\wedge\xi^K)\|_{\mathbb{L}^2}^2$, we have
\begin{equation}\label{1pro51}
\begin{split}
&\frac{1}{2}\|\mathbf{u}^*(t\wedge\xi^K)\|_{\mathbb{L}^2}^2+\beta_2\int_0^{t\wedge\xi^K}\|\Delta\mathbf{u}^*\|_{\mathbb{L}^2}^2\,\mathrm{d}s\\
&=-\beta_1\int_0^{t\wedge\xi^K}\|\nabla\mathbf{u}^*\|_{\mathbb{L}^2}^2\,\mathrm{d}s+\beta_3\int_0^{t\wedge\xi^K}\left(\mathbf{u}^*-|\mathbf{u}_1|^2\mathbf{u}^*+\mathbf{u}_2\left(|\mathbf{u}_2|^2-|\mathbf{u}_1|^2\right),\mathbf{u}^*\right)_{\mathbb{L}^2}\,\mathrm{d}s\\
&+\beta_4\int_0^{t\wedge\xi^K}\left(\mathbf{u}_1\times\nabla\mathbf{u}_1-\mathbf{u}_2\times\nabla\mathbf{u}_2,\nabla\mathbf{u}^*\right)_{\mathbb{L}^2}\,\mathrm{d}s\\
&-\beta_5\int_0^{t\wedge\xi^K}\theta_R(\|\nabla \mathbf{u}_1\|_{\mathbb{L}^2})\left(\nabla\left(|\mathbf{u}_1|^2\mathbf{u}^*+\left(|\mathbf{u}_1|^2-|\mathbf{u}_2|^2\right)\mathbf{u}_2\right),\nabla\mathbf{u}^*\right)_{\mathbb{L}^2}\,\mathrm{d}s\\
&-\beta_5\int_0^{t\wedge\xi^K}\left(\theta_R(\|\nabla \mathbf{u}_1\|_{\mathbb{L}^2})-\theta_R(\|\nabla \mathbf{u}_2\|_{\mathbb{L}^2})\right)\left(\nabla\left(|\mathbf{u}_2|^2\mathbf{u}_2\right),\nabla\mathbf{u}^*\right)_{\mathbb{L}^2}\,\mathrm{d}s\\
&+\frac{1}{2}\sum_{j=1}^{\infty}\int_0^{t\wedge\xi^K}\|\mathbf{u}^*\times \mathbf{h}_j\|_{\mathbb{L}^2}^2\,\mathrm{d}s:=\sum_{i=1}^6I_i.
\end{split}
\end{equation}
By using integration by parts, H\"{o}lder's inequality and Young's inequality, we have
\begin{equation}\label{2pro51}
\begin{split}
|I_1|\leq\varepsilon\int_0^{t\wedge\xi^K}\|\Delta\mathbf{u}^*\|_{\mathbb{L}^2}^2\,\mathrm{d}s+C_{\varepsilon}\int_0^{t\wedge\xi^K}\|\mathbf{u}^*\|_{\mathbb{L}^2}^2\,\mathrm{d}s.
\end{split}
\end{equation}
By using the triangle inequality, we see that
\begin{equation}\label{3pro51}
\begin{split}
I_2\leq C\int_0^{t\wedge\xi^K}\|\mathbf{u}^*\|_{\mathbb{L}^2}^2\,\mathrm{d}s+C\int_0^{t\wedge\xi^K}\|\mathbf{u}_2\|_{\mathbb{L}^{\infty}}\left(\|\mathbf{u}_1\|_{\mathbb{L}^{\infty}}+\|\mathbf{u}_2\|_{\mathbb{L}^{\infty}}\right)\|\mathbf{u}^*\|_{\mathbb{L}^2}^2\,\mathrm{d}s.
\end{split}
\end{equation}
Similarly, by the H\"{o}lder inequality and Young¡¯s inequality, it follows that
\begin{equation}\label{4pro51}
\begin{split}
|I_3|&\leq\beta_4\left|\int_0^{t\wedge\xi^K}\left(\mathbf{u}_1\times\mathbf{u}^*,\Delta\mathbf{u}^*\right)_{\mathbb{L}^2}\,\mathrm{d}s\right|\\
&\leq C\int_0^{t\wedge\xi^K}\|\mathbf{u}_1\|_{\mathbb{L}^{\infty}}\|\mathbf{u}^*\|_{\mathbb{L}^{2}}\|\Delta\mathbf{u}^*\|_{\mathbb{L}^{2}}\,\mathrm{d}s\\
&\leq\varepsilon\int_0^{t\wedge\xi^K}\|\Delta\mathbf{u}^*\|_{\mathbb{L}^{2}}^{2}\,\mathrm{d}s+C_{\varepsilon}\int_0^{t\wedge\xi^K}\|\mathbf{u}_1\|_{\mathbb{L}^{\infty}}^2\|\mathbf{u}^*\|_{\mathbb{L}^{2}}^2\,\mathrm{d}s.
\end{split}
\end{equation}
For the fourth term,  by using integration by parts, we see that
\begin{equation}\label{5pro51}
\begin{split}
&|I_4|=\left|\beta_5\int_0^{t\wedge\xi^K}\theta_R(\|\nabla \mathbf{u}_1\|_{\mathbb{L}^2})\left(|\mathbf{u}_1|^2\mathbf{u}^*+\left(|\mathbf{u}_1|^2-|\mathbf{u}_2|^2\right)\mathbf{u}_2,\Delta\mathbf{u}^*\right)_{\mathbb{L}^2}\,\mathrm{d}s\right|\\
&\leq C\int_0^{t\wedge\xi^K}\|\mathbf{u}_1\|_{\mathbb{L}^{\infty}}^2\|\mathbf{u}^*\|_{\mathbb{L}^2}\|\Delta\mathbf{u}^*\|_{\mathbb{L}^2}\,\mathrm{d}s\\
&+C\int_0^{t\wedge\xi^K}(\|\mathbf{u}_1\|_{\mathbb{L}^{\infty}}\|\mathbf{u}_2\|_{\mathbb{L}^{\infty}}+\|\mathbf{u}_2\|_{\mathbb{L}^{\infty}}^2)\|\mathbf{u}^*\|_{\mathbb{L}^2}\|\Delta\mathbf{u}^*\|_{\mathbb{L}^2}\,\mathrm{d}s\\
&\leq \varepsilon\int_0^{t\wedge\xi^K}\|\Delta\mathbf{u}^*\|_{\mathbb{L}^{2}}^{2}\,\mathrm{d}s+C_{\varepsilon}\int_0^{t\wedge\xi^K}(\|\mathbf{u}_1\|_{\mathbb{L}^{\infty}}^4+\|\mathbf{u}_2\|_{\mathbb{L}^{\infty}}^4)\|\mathbf{u}^*\|_{\mathbb{L}^2}^2\,\mathrm{d}s.
\end{split}
\end{equation}
For the fifth term,  noting the property of $\theta_R(\cdot)$, we see that
\begin{equation}\label{6pro51}
\begin{split}
&|I_5|\leq\beta_5\int_0^{t\wedge\xi^K}\left|\theta_R(\|\nabla \mathbf{u}_1\|_{\mathbb{L}^2})-\theta_R(\|\nabla \mathbf{u}_2\|_{\mathbb{L}^2})\right|\left|\left(|\mathbf{u}_2|^2\mathbf{u}_2,\Delta\mathbf{u}^*\right)_{\mathbb{L}^2}\right|\,\mathrm{d}s\\
&\leq C\int_0^{t\wedge\xi^K}\|\nabla \mathbf{u}^*\|_{\mathbb{L}^{2}}\|\mathbf{u}_2\|_{\mathbb{L}^{6}}^3\|\Delta\mathbf{u}^*\|_{\mathbb{L}^{2}}\,\mathrm{d}s\\
&\leq \varepsilon\int_0^{t\wedge\xi^K}\|\Delta\mathbf{u}^*\|_{\mathbb{L}^{2}}^2\,\mathrm{d}s+ C_{\varepsilon}\int_0^{t\wedge\xi^K}\|\nabla \mathbf{u}^*\|_{\mathbb{L}^{2}}^2\|\mathbf{u}_2\|_{\mathbb{H}^{1}}^6\,\mathrm{d}s\\
&\leq \varepsilon\int_0^{t\wedge\xi^K}\|\Delta\mathbf{u}^*\|_{\mathbb{L}^{2}}^2\,\mathrm{d}s+C_{\varepsilon}K^3\int_0^{t\wedge\xi^K}\|\nabla \mathbf{u}^*\|_{\mathbb{L}^{2}}^2\,\mathrm{d}s\\
&\leq \varepsilon\int_0^{t\wedge\xi^K}\|\Delta\mathbf{u}^*\|_{\mathbb{L}^{2}}^2\,\mathrm{d}s+C_{\varepsilon,K}\int_0^{t\wedge\xi^K}\|\mathbf{u}^*\|_{\mathbb{L}^{2}}^2\,\mathrm{d}s.
\end{split}
\end{equation}
For the last term, by using the condition \eqref{condh}, it follows that
\begin{equation}\label{1-6pro51}
\begin{split}
&|I_6|\leq C_{\mathbf{h}}\int_0^{t\wedge\xi^K}\|\mathbf{u}^*\|_{\mathbb{L}^2}^2\,\mathrm{d}s
\end{split}
\end{equation}
Thus plugging \eqref{2pro51}-\eqref{1-6pro51} into \eqref{1pro51} and choosing $\varepsilon$ small enough, we infer that
\begin{equation*}
\begin{split}
\|\mathbf{u}^*(t\wedge\xi^K)\|_{\mathbb{L}^2}^2\leq C_{K}\int_0^{t\wedge\xi^K}\mathbf{F}(s)\|\mathbf{u}^*(s)\|_{\mathbb{L}^{2}}^2\,\mathrm{d}s,
\end{split}
\end{equation*}
where
$
\mathbf{F}:=1+\|\mathbf{u}_1\|_{\mathbb{L}^{\infty}}^4+\|\mathbf{u}_2\|_{\mathbb{L}^{\infty}}^4.
$
According to the GN inequality, $
\|f\|_{\mathbb{L}^{\infty}}\leq C\|f\|_{\mathbb{H}^{1}}^{\frac{1}{2}}\|f\|_{\mathbb{H}^{2}}^{\frac{1}{2}},~f\in\mathbb{H}^2.
$ is valid for $d=1,2,3$. Thus
\begin{equation*}
\begin{split}
\int_0^{t\wedge\xi^K}\mathbf{F}(s)\,\mathrm{d}s&\leq t+C\sup_{s\in[0,t\wedge\xi^K]}\|\mathbf{u}_1(s)\|_{\mathbb{H}^{1}}^2\int_0^{t\wedge\xi^K}\|\mathbf{u}_1(s)\|_{\mathbb{H}^{2}}^2\,\mathrm{d}s\\
&+C\sup_{s\in[0,t\wedge\xi^K]}\|\mathbf{u}_2(s)\|_{\mathbb{H}^{1}}^2\int_0^{t\wedge\xi^K}\|\mathbf{u}_2(s)\|_{\mathbb{H}^{2}}^2\,\mathrm{d}s<\infty.
\end{split}
\end{equation*}
Thus by the Gronwall lemma we have $\sup_{s\in[0,t\wedge\xi^K]}\|\mathbf{u}^*(s)\|_{\mathbb{L}^2}^2=0$, a.s. By using the monotone convergence theorem and the fact that $\xi^K\rightarrow T$ as $K\rightarrow\infty$, it follows that $\mathbb{P}$-a.s.,
$
\sup_{s\in[0,T]}\|\mathbf{u}^*(s)\|_{\mathbb{L}^2}^2=0,
$
which implies the uniqueness.~~$\Box$

\begin{corollary}\label{cor51} Let $\mathcal{O}\subset\mathbb{R}^d,~d=1,2,3$. Then for $\mathbf{u}_0\in \mathbb{H}^1$,
\begin{enumerate}
\item [(1)] there exists a unique pathwise weak solution of equation \eqref{sys1-2};
\item [(2)] the martingale solution of \eqref{sys1-2} is unique in law.
\end{enumerate}
\end{corollary}
\noindent\textbf{Proof.} The corollary follows from Theorem 2.2 and 12.1 in \cite{35}.~~$\Box$\\
\textbf{Proof of Theorem \ref{the1}.}
Now we aim is to remove the truncation function $\theta_R(\cdot)$ in order to construct the local solution in the sense of Definition \ref{def1} for \eqref{sys1}. To this end, we define the stopping time
\begin{equation*}\label{equ52}
\begin{split}
\tau^K:=\inf\left\{t\geq0,~\|\mathbf{u}(t)\|_{\mathbb{H}^1}>K\right\},~\forall K>0.
 \end{split}
\end{equation*}
We assume first that $\|\mathbf{u}_0\|_{\mathbb{H}^1}\leq K$ for some deterministic $K>0$, and choose large enough $R>K$. Then $\tau^K$ is strictly positive $\mathbb{P}$-a.s., and $\|\nabla \mathbf{u}(t)\|_{\mathbb{L}^2}\leq R$ for any $t\in[0,\tau^K]$, which means that
\begin{equation*}
\begin{split}
\theta_R(\|\nabla \mathbf{u}(t)\|_{\mathbb{L}^2})\equiv1~\textrm{on}~[0,\tau^K].
\end{split}
\end{equation*}
This together with Corollary \ref{cor51} implies that the pair $\left(\mathbf{u},\tau^K\right)$ is a unique local pathwise solution of system \eqref{sys1}. To pass to the general case $\|\mathbf{u}_0\|_{\mathbb{H}^1}<\infty$, we decompose $\mathbf{u}_0$ as
\begin{equation*}
\begin{split}
\mathbf{u}_0=\sum_{m\in\mathbb{N}}\mathbf{u}_0^m,~\textrm{where~}\mathbf{u}_0^m:=\mathbf{u}_0I_{\{m\leq\|\mathbf{u}_0\|_{\mathbb{H}^1}<m+1\}}.
\end{split}
\end{equation*}
Then for each $m\in\mathbb{N}$, one can choose $R>m+1$ to construct a unique local pathwise solution $\left(\mathbf{u}^m,\tau^m\right)$ to \eqref{sys1} with initial profile $\mathbf{u}_0^m$. Then we define a pair $\left(\mathbf{u},\tau\right)$ by
\begin{equation*}
\begin{split}
\mathbf{u}=\sum_{m\in\mathbb{N}}\mathbf{u}^mI_{\{m\leq\|\mathbf{u}_0\|_{\mathbb{H}^1}<m+1\}},~\tau=\sum_{m\in\mathbb{N}}\tau^mI_{\{m\leq\|\mathbf{u}_0\|_{\mathbb{H}^1}<m+1\}}.
\end{split}
\end{equation*}
Since $\mathbf{u}^m\in L^{\infty}([0,T];\mathbb{H}^1)\cap L^2(0,T;\mathbb{H}^3)$ a.s., it follows that $\mathbf{u}\in L^{\infty}([0,T];\mathbb{H}^1)\cap L^2(0,T;\mathbb{H}^3)$, $\mathbb{P}$-a.s. Moreover, since $\|\mathbf{u}_0\|_{\mathbb{H}^1}<\infty$, it follows that $\sum_{m\in\mathbb{N}}I_{\{m\leq\|\mathbf{u}_0\|_{\mathbb{H}^1}<m+1\}}\equiv1$. Thus by using the uniform bounds provided by Lemma \ref{lem54}, we infer that for all $p\geq1$
\begin{equation*}
\begin{split}
\mathbb{E}\sup_{t\in[0,\tau]}\|\mathbf{u}(t)\|_{\mathbb{H}^1}^p=\sum_{m\in\mathbb{N}}I_{\{m\leq\|\mathbf{u}_0\|_{\mathbb{H}^1}<m+1\}}\mathbb{E}\sup_{t\in[0,\tau^m]}\|\mathbf{u}^m(t)\|_{\mathbb{H}^1}^p<\infty,
\end{split}
\end{equation*}
and
\begin{equation*}
\begin{split}
\mathbb{E}\left(\int_0^{\tau}\|\mathbf{u}(t)\|_{\mathbb{H}^3}^2\,\mathrm{d}t\right)^p=\sum_{m\in\mathbb{N}}I_{\{m\leq\|\mathbf{u}_0\|_{\mathbb{H}^1}<m+1\}}\mathbb{E}\left(\int_0^{\tau^m}\|\mathbf{u}^m(t)\|_{\mathbb{H}^3}^2\,\mathrm{d}t\right)^p<\infty.
\end{split}
\end{equation*}
Moreover, it follows that
\begin{equation*}
\begin{split}
&\mathbf{u}(t\wedge\tau)=\sum_{m\in\mathbb{N}}I_{\{m\leq\|\mathbf{u}_0\|_{\mathbb{H}^1}<m+1\}}\mathbf{u}^m(t\wedge\tau^m)\\
&=\sum_{m\in\mathbb{N}}I_{\{m\leq\|\mathbf{u}_0\|_{\mathbb{H}^1}<m+1\}}\Bigl(\mathbf{u}_0^m+\int_0^{t\wedge\tau^m}\beta_1\Delta \mathbf{u}^m-\beta_2\Delta^2\mathbf{u}^m+\beta_3(1-|\mathbf{u}^m|^2)\mathbf{u}^m-\beta_4\mathbf{u}^m\times\Delta\mathbf{u}^m\\
&+\beta_5\Delta(|\mathbf{u}^m|^2\mathbf{u}^m)\,\mathrm{d}s\Bigl)+\sum_{m\in\mathbb{N}}I_{\{m\leq\|\mathbf{u}_0\|_{\mathbb{H}^1}<m+1\}}\int_0^{t\wedge\tau^m}\sum_{j=1}^{\infty}(-\mathbf{u}^m\times \mathbf{h}_j+\mathbf{h}_j-\Delta\mathbf{h}_j)\circ\mathrm{d}W_j(s)\\
&=\mathbf{u}_0+\int_0^{t\wedge\tau}\beta_1\Delta \mathbf{u}-\beta_2\Delta^2\mathbf{u}+\beta_3(1-|\mathbf{u}|^2)\mathbf{u}-\beta_4\mathbf{u}\times\Delta\mathbf{u}+\beta_5\Delta(|\mathbf{u}|^2\mathbf{u})\,\mathrm{d}s\\
&+\int_0^{t\wedge\tau}\sum_{j=1}^{\infty}(-\mathbf{u}\times \mathbf{h}_j+\mathbf{h}_j-\Delta\mathbf{h}_j)\circ\mathrm{d}W_j(s).
\end{split}
\end{equation*}
Thus the pair $\left(\mathbf{u},\tau\right)$ constructed above is a unique local pathwise solution to \eqref{sys1}. Moreover, by using a standard argument provided by \cite{24}, we can extend the solution $\left(\mathbf{u},\tau\right)$ to a maximal time of existence. The proof of Theorem \ref{the1} is thus completed.~~$\Box$

\section{Global pathwise weak solutions and invariant measures for $d=1$}\label{sec5}
The first goal in this section is devoted to prove \eqref{sys1} admits a unique global pathwise weak solution for $d=1$. To this end, we consider the following Galerkin approximation system which does not require the introduction of the truncation function $\theta_R(\cdot)$.
\begin{equation}\label{sys3}
\left\{
\begin{aligned}
&\mathrm{d}\mathbf{u}_n=\bigl[\beta_1\Delta \mathbf{u}_n-\beta_2\Delta^2\mathbf{u}_n+\beta_3\Pi_n\left((1-|\mathbf{u}_n|^2)\mathbf{u}_n\right)-\beta_4\Pi_n(\mathbf{u}_n\times\Delta\mathbf{u}_n)\\
&~~+\beta_5\Pi_n\Delta(|\mathbf{u}_n|^2\mathbf{u}_n)\bigl]\,\mathrm{d}t+\sum_{j=1}^{n}\Pi_n(-\mathbf{u}_n\times \mathbf{h}_j+\mathbf{h}_j-\Delta\mathbf{h}_j)\circ\mathrm{d}W_j(t),&&\textrm{in}~(0,\infty)\times\mathcal{O},\\
&\mathbf{u}_n(0)=\Pi_n\mathbf{u}_0,&&\textrm{in}~\mathcal{O}.
\end{aligned}
\right.
\end{equation}
To obtain the global solution, we need the following uniformly bound.
\begin{lemma}\label{lem62} Let $\mathcal{O}\subset\mathbb{R}$ be a bounded domain with $C^{2,1}$-boundary. Then for any $p\geq1$, $n\in\mathbb{N}$ and $t\in[0,T]$, there exists a positive constant $C=C(\|\mathbf{u}_0\|_{\mathbb{H}^1},p,\mathbf{h},T)$ independent of $n$ such that
\begin{equation}\label{111lem62}
\begin{split}
\mathbb{E}\sup_{s\in[0,t]}\|\mathbf{u}_n(s)\|_{\mathbb{H}^1}^{2p}+\mathbb{E}\left(\int_0^t\|\mathbf{u}_n(s)\|_{\mathbb{H}^3}^2\,\mathrm{d}s\right)^p\leq C.
\end{split}
\end{equation}
\end{lemma}
\noindent\textbf{Proof.} Similar to the proof of Lemma \ref{lem32}, we have
\begin{equation}\label{111lem61}
\begin{split}
&\mathbb{E}\sup_{s\in[0,t]}\|\mathbf{u}_n(s)\|_{\mathbb{L}^2}^{2p}+\mathbb{E}\left(\int_0^t\| \mathbf{u}_n(s)\|_{\mathbb{H}^1}^2\,\mathrm{d}s\right)^p\leq C.
\end{split}
\end{equation}
By a similar argument as shown in Lemma \ref{lem33}, we deduce that
\begin{equation}\label{1lem63}
\begin{split}
&\|\nabla\mathbf{u}_n(t)\|_{\mathbb{L}^2}^2+\int_0^t\|\nabla\Delta \mathbf{u}_n\|_{\mathbb{L}^2}^2\,\mathrm{d}s+\int_0^t\||\mathbf{u}_n||\nabla\mathbf{u}_n|\|_{\mathbb{L}^2}^2\,\mathrm{d}s\\\
&\leq C\|\nabla\mathbf{u}_n(0)\|_{\mathbb{L}^2}^2+C\int_0^t\|\nabla\mathbf{u}_n\|_{\mathbb{L}^2}^2\,\mathrm{d}s+C\int_0^t\|\mathbf{u}_n\|_{\mathbb{L}^2}^2\,\mathrm{d}s+C\int_0^t\|\nabla\mathbf{u}_n\|_{\mathbb{L}^4}^4\,\mathrm{d}s\\
&+\sum_{j=1}^n\int_0^t\left(\nabla G_{nj}(\mathbf{u}_n),\nabla\mathbf{u}_n\right)_{\mathbb{L}^2}\,\mathrm{d}W_j(s).
\end{split}
\end{equation}
Plugging \eqref{ad2lem33} and \eqref{ad3lem33} into \eqref{1lem63} and choosing $\varepsilon$ small enough, we infer that
\begin{equation}\label{4lem63}
\begin{split}
&\|\nabla\mathbf{u}_n(t)\|_{\mathbb{L}^2}^2+\int_0^t\|\nabla\Delta \mathbf{u}_n\|_{\mathbb{L}^2}^2\,\mathrm{d}s+\int_0^t\||\mathbf{u}_n||\nabla\mathbf{u}_n|\|_{\mathbb{L}^2}^2\,\mathrm{d}s\\\
&\leq C\|\nabla\mathbf{u}_n(0)\|_{\mathbb{L}^2}^2+C\int_0^t\|\mathbf{u}_n\|_{\mathbb{L}^2}^2\,\mathrm{d}s+C\int_0^t\|\mathbf{u}_n\|_{\mathbb{L}^2}^{14}\,\mathrm{d}s+C\int_0^t\|\nabla\mathbf{u}_n\|_{\mathbb{L}^2}^2\,\mathrm{d}s\\
&+\sum_{j=1}^n\int_0^t\left(\nabla G_{nj}(\mathbf{u}_n),\nabla\mathbf{u}_n\right)_{\mathbb{L}^2}\,\mathrm{d}W_j(s).
\end{split}
\end{equation}
Recalling the estimate \eqref{15lem33}, it follows from \eqref{4lem63} and \eqref{111lem61} that
\begin{equation*}
\begin{split}
&\mathbb{E}\sup_{s\in[0,t]}\|\nabla\mathbf{u}_n(s)\|_{\mathbb{L}^2}^{2p}+\mathbb{E}\left(\int_0^t\|\nabla\Delta\mathbf{u}_n(s)\|_{\mathbb{L}^2}^2\,\mathrm{d}s\right)^p+\mathbb{E}\left(\int_0^t\||\mathbf{u}_n||\nabla\mathbf{u}_n|\|_{\mathbb{L}^2}^2\,\mathrm{d}s\right)^p\\
&\leq C+C\int_0^t\mathbb{E}\|\nabla\mathbf{u}_n(s)\|_{\mathbb{L}^2}^{2p}\,\mathrm{d}s,
\end{split}
\end{equation*}
which together with Gronwall lemma and \eqref{111lem61} implies the result \eqref{111lem62}.~~$\Box$

\begin{proposition}\label{pro61} Under the same assumptions as in Lemma \ref{lem62}, the system \eqref{sys1} admits a unique global pathwise weak solution in the sense of Definition \ref{def1}.
\end{proposition}
\noindent\textbf{Proof.} With the uniform boundedness provided by Lemma \ref{lem62}, we can once again obtain the tightness result from section \ref{sec3}, as well as the pointwise convergence result for the approximate solutions provided by section \ref{sec4}. Therefore, \eqref{sys1} possesses at least one global martingale weak solutions. Additionally, similar to the discussion in Proposition \ref{pro51}, it is not difficult to deduce that the pathwise solutions of \eqref{sys1} are unique. By Yamada-Watanabe theorem \cite{35,44}, we conclude that \eqref{sys1} possesses a unique global pathwise weak solution.~~$\Box$

The second goal in this section is devoted to show that the existence of invariant measure for equation \eqref{sys1}. Our proof mainly draws inspiration from the ideas presented in \cite{10,12,13}.
\begin{lemma}\label{lem71} Let $\mathbf{u}$ be a weak solution to \eqref{sys1} with properties listed in Theorem \ref{the2}. Then there exists a positive constant $C$ depending on $C_0$ and $\mathbf{h}$ such that for all $t\geq0$
\begin{equation}\label{1lem71}
\begin{split}
\int_0^t\mathbb{E}\|\mathbf{u}(s)\|_{\mathbb{H}^2}^2\,\mathrm{d}s\leq C(1+t).
\end{split}
\end{equation}
\end{lemma}
\noindent\textbf{Proof.} We will use It\^{o}'s formula in the form presented in Pardoux's fundamental work \cite{36}. According to Proposition \ref{pro61} and Lemma \ref{lem62}, it follows that
\begin{equation*}
\begin{split}
&\mathbb{E}\int_0^T\|\mathbf{u}\|_{\mathbb{H}^1}^2\,\mathrm{d}s+\sum_{j=1}^{\infty}\mathbb{E}\int_0^T\|-\mathbf{u}\times \mathbf{h}_j+\mathbf{h}_j-\Delta\mathbf{h}_j\|_{\mathbb{L}^2}^2\,\mathrm{d}s<\infty,\\
&\mathbb{E}\int_0^T\|\Delta\mathbf{u}\|_{(\mathbb{H}^{1})^*}^2\,\mathrm{d}s+\mathbb{E}\int_0^T\|\Delta^2\mathbf{u}\|_{(\mathbb{H}^{1})^*}^2\,\mathrm{d}s+\mathbb{E}\int_0^T\|\mathbf{u}\times\Delta\mathbf{u}\|_{(\mathbb{H}^{1})^*}^2\,\mathrm{d}s+\mathbb{E}\int_0^T\|\Delta(\left|\mathbf{u}\right|^2\mathbf{u})\|_{(\mathbb{H}^{1})^*}^2\,\mathrm{d}s<\infty.
\end{split}
\end{equation*}
Choosing $H=\mathbb{L}^2$ and $V=\mathbb{H}^1$, then the process satisfies the assumptions of Pardoux's theorem, and thus we have
\begin{equation}\label{2lem71}
\begin{split}
&\frac{1}{2}\|\mathbf{u}(t)\|_{\mathbb{L}^2}^2+\beta_1\int_0^t\|\nabla \mathbf{u}(s)\|_{\mathbb{L}^2}^2\,\mathrm{d}s+\beta_2\int_0^t\|\Delta \mathbf{u}(s)\|_{\mathbb{L}^2}^2\,\mathrm{d}s+\beta_3\int_0^t\| \mathbf{u}(s)\|_{\mathbb{L}^4}^4\,\mathrm{d}s\\
&=\frac{1}{2}\|\mathbf{u}_0\|_{\mathbb{L}^2}^2+\beta_3\int_0^t\| \mathbf{u}(s)\|_{\mathbb{L}^2}^2\,\mathrm{d}s-\beta_5\int_0^t\left(\nabla\left(|\mathbf{u}(s)|\mathbf{u}(s)\right),\nabla\mathbf{u}(s)\right)_{\mathbb{L}^2}\,\mathrm{d}s\\
&-\frac{1}{2}\sum_{j=1}^{\infty}\int_0^t\left(G_{j}(\mathbf{u}(s))\times\mathbf{h}_j,\mathbf{u}(s)\right)_{\mathbb{L}^2}\,\mathrm{d}s+\frac{1}{2}\sum_{j=1}^{\infty}\int_0^t\|G_{j}(\mathbf{u}(s))\|_{\mathbb{L}^2}^2\,\mathrm{d}s\\
&+\sum_{j=1}^{\infty}\int_0^t\left(G_{j}(\mathbf{u}(s)),\mathbf{u}(s)\right)_{\mathbb{L}^2}\,\mathrm{d}W_j(s),
\end{split}
\end{equation}
where $G_{j}(\mathbf{u}):=-\mathbf{u}\times \mathbf{h}_j+\mathbf{h}_j-\Delta\mathbf{h}_j$. By a similar argument as shown in \eqref{4lem32}-\eqref{6lem32}, we have
\begin{equation*}
\begin{split}
\sum_{j=1}^{\infty}\left|\left(G_{j}(\mathbf{u})\times\mathbf{h}_j,\mathbf{u}\right)_{\mathbb{L}^2}\right|+\sum_{j=1}^{\infty}\|G_{j}(\mathbf{u})\|_{\mathbb{L}^2}^2\leq C_{\mathbf{h}}+C_{\mathbf{h}}\|\mathbf{u}\|_{\mathbb{L}^2}^2,
\end{split}
\end{equation*}
and
$
\|\nabla \mathbf{u}\|_{\mathbb{L}^2}^2\leq\varepsilon\|\Delta\mathbf{u}\|_{\mathbb{L}^2}^2+C_{\varepsilon}\|\mathbf{u}\|_{\mathbb{L}^2}^2.
$
Choosing $\varepsilon$ small enough, it follows from \eqref{2lem71} that
\begin{equation}\label{3lem71}
\begin{split}
&\|\mathbf{u}(t)\|_{\mathbb{L}^2}^2+\int_0^t\|\Delta \mathbf{u}\|_{\mathbb{L}^2}^2\,\mathrm{d}s+\int_0^t\| \mathbf{u}\|_{\mathbb{L}^4}^4\,\mathrm{d}s+\int_0^t\|\mathbf{u}\cdot\nabla\mathbf{u}\|_{\mathbb{L}^2}^2\,\mathrm{d}s+\int_0^t\| |\mathbf{u}||\nabla \mathbf{u}|\|_{\mathbb{L}^2}^2\,\mathrm{d}s\\
&\leq C\|\mathbf{u}_0\|_{\mathbb{L}^2}^2+C_{\mathbf{h}}+C_{\mathbf{h}}\int_0^t\|\mathbf{u}\|_{\mathbb{L}^2}^2\,\mathrm{d}s+C\sum_{j=1}^{\infty}\int_0^t\left(G_{j}(\mathbf{u}),\mathbf{u}\right)_{\mathbb{L}^2}\,\mathrm{d}W_j(s).
\end{split}
\end{equation}
By using \eqref{condh} and Proposition \ref{pro61} and noting $\left(\mathbf{u}\times\mathbf{h}_j,\mathbf{u}\right)_{\mathbb{L}^2}=0$, we see that
\begin{equation*}
\begin{split}
&\mathbb{E}\int_0^t\left(G_{j}(\mathbf{u}),\mathbf{u}\right)_{\mathbb{L}^2}^2\,\mathrm{d}s=\mathbb{E}\int_0^t\left(\mathbf{h}_j-\Delta\mathbf{h}_j,\mathbf{u}\right)_{\mathbb{L}^2}^2\,\mathrm{d}s\\
&\leq2\left(\|\mathbf{h}_j\|_{\mathbb{L}^2}^2+\|\Delta\mathbf{h}_j\|_{\mathbb{L}^2}^2\right)\mathbb{E}\int_0^t\|\mathbf{u}\|_{\mathbb{L}^2}^2\,\mathrm{d}s<\infty,
\end{split}
\end{equation*}
which means that the process $t\rightarrow\int_0^t\left(G_{j}(\mathbf{u}),\mathbf{u}\right)_{\mathbb{L}^2}\,\mathrm{d}s$ is a martingale on $[0,T]$. In particular, it follows that
$
\mathbb{E}\int_0^t\left(G_{j}(\mathbf{u}),\mathbf{u}\right)_{\mathbb{L}^2}\,\mathrm{d}W_j(s)=0.
$
Hence, by Young's inequality, we infer from \eqref{3lem71} that
\begin{equation}\label{4lem71}
\begin{split}
&\mathbb{E}\|\mathbf{u}(t)\|_{\mathbb{L}^2}^2+\mathbb{E}\int_0^t\|\Delta \mathbf{u}\|_{\mathbb{L}^2}^2\,\mathrm{d}s+\mathbb{E}\int_0^t\| \mathbf{u}\|_{\mathbb{L}^4}^4\,\mathrm{d}s+\mathbb{E}\int_0^t\|\mathbf{u}\cdot\nabla\mathbf{u}\|_{\mathbb{L}^2}^2\,\mathrm{d}s+\mathbb{E}\int_0^t\| |\mathbf{u}||\nabla \mathbf{u}|\|_{\mathbb{L}^2}^2\,\mathrm{d}s\\
&\leq C\mathbb{E}\|\mathbf{u}_0\|_{\mathbb{L}^2}^2+C_{\mathbf{h}}+C_{\mathbf{h}}\mathbb{E}\int_0^t\|\mathbf{u}\|_{\mathbb{L}^2}^2\,\mathrm{d}s\\
&\leq C\mathbb{E}\|\mathbf{u}_0\|_{\mathbb{L}^2}^2+C_{\mathbf{h}}+C_{\mathbf{h}}t+\frac{1}{2}\mathbb{E}\int_0^t\|\mathbf{u}\|_{\mathbb{L}^4}^4\,\mathrm{d}s,
\end{split}
\end{equation}
which implies that
$
\mathbb{E}\|\mathbf{u}(t)\|_{\mathbb{L}^2}^2+\mathbb{E}\int_0^t\|\mathbf{u}\|_{\mathbb{H}^2}^2\,\mathrm{d}s+\mathbb{E}\int_0^t\| |\mathbf{u}||\nabla \mathbf{u}|\|_{\mathbb{L}^2}^2\,\mathrm{d}s\leq C(1+t).
$
~~$\Box$\\
Let $\mathcal{B}_b(\mathbb{H}^1)$ denote the set of all bounded and Borel measurable functions on $\mathbb{H}^1$. According to Proposition \ref{pro61}, the unique solution $\mathbf{u}$ is a $\mathbb{H}^1$-valued Markov process. Thus for any $\psi\in\mathcal{B}_b(\mathbb{H}^1)$, $t\geq0$, we can define the transition semigroup $P_t\psi:\mathbb{H}^1\rightarrow\mathbb{R}$ by
\begin{equation}\label{ad71}
\begin{split}
P_t\psi(\mathbf{u}_0):=\mathbb{E}\psi(\mathbf{u}(t;\mathbf{u}_0)),~\mathbf{u}_0\in\mathbb{H}^1,
\end{split}
\end{equation}
where $\mathbf{u}(t;\mathbf{u}_0)$ stands for the process $\mathbf{u}$ starting at time $t=0$ and $\mathbf{u}(0)=\mathbf{u}_0$. Next we shall show the sequentially weak Feller property of $\{P_t\}$. Before doing so, we need establish the following convergence result. For $p\geq4,~\beta>1$ and $T>0$ let
\begin{equation*}
\begin{split}
\tilde{\mathcal{X}}&:=C([0,T];(\mathbb{H}^{\beta})^*)\cap L^p(0,T;\mathbb{L}^4)\cap L^2(0,T;\mathbb{H}^2)\cap C([0,T];\mathbb{H}^1_w)\\
&=\mathcal{X}\cap C([0,T];\mathbb{H}^1_w),
\end{split}
\end{equation*}
and let $\mathcal{T}$ denote the supremum of the corresponding four topologies.
\begin{lemma}\label{lem72} Assume that a $\mathbb{H}^1$-valued sequence $\{\mathbf{u}_{0,k}\}_{k\in\mathbb{N}}$ is convergent weakly in $\mathbb{H}^1$ to $\mathbf{u}_0\in\mathbb{H}^1$. Let $C_0$ be a positive constant such that $\sup_{k\in\mathbb{N}}\|\mathbf{u}_{0,k}\|_{\mathbb{H}^1}\leq C_0$. Let $\left(\Omega,\mathcal{F},\mathbb{F},\mathbb{P},W,\mathbf{u}^k\right)$ be a unique solution of \eqref{sys1} with the initial data $\mathbf{u}_{0,k}$. Then there exist
\begin{enumerate}
\item [(1)] a subsequence $\{k_m\}_m$,
\item [(2)] a stochastic basis $\left(\tilde{\Omega},\tilde{\mathcal{F}},\tilde{\mathbb{F}},\tilde{\mathbb{P}}\right)$,
\item [(3)] a standard $\tilde{\mathbb{F}}$-Wiener process $\tilde{W}=(\tilde{W}_j)_{j=1}^{\infty}$ defined on this basis,
\item [(4)] progressively measurable process $\tilde{\mathbf{u}}$, $\{\tilde{\mathbf{u}}_{k_m}\}_{m\in\mathbb{N}}$ (defined on this basis) with laws supported in $(\tilde{\mathcal{X}},\mathcal{T})$ such that
\begin{equation*}
\begin{split}
&\tilde{\mathbf{u}}_{k_m}~\textrm{has the same law as}~\tilde{\mathbf{u}}_{k_m}~\textrm{on}~\tilde{\mathcal{X}},\\
&\tilde{\mathbf{u}}_{k_m}\rightarrow\tilde{\mathbf{u}}~\textrm{in}~\tilde{\mathcal{X}}~\textrm{as}~m\rightarrow\infty,~\tilde{\mathbb{P}}\textrm{-a.s.},\\
&\tilde{\mathbf{u}}~\textrm{is a solution of equation}~\eqref{sys1}~\textrm{with the initial data}~\mathbf{u}_0.
\end{split}
\end{equation*}
\end{enumerate}
\end{lemma}
\noindent\textbf{Proof.} The proof is relied on the Jakubowski's version of the Skorokhod theorem \cite{27}. According to Proposition \ref{pro61}, for given $\mathbf{u}_{0,k}\in\mathbb{H}^1$ there exists a unique solution $\mathbf{u}^k$ to \eqref{sys1} defined on the stochastic basis $\left(\Omega,\mathcal{F},\mathbb{F},\mathbb{P}\right)$. Let
$
\mathcal{Y}(\beta):= L^{\infty}(0,T;\mathbb{H}^1)\cap L^{2}(0,T;\mathbb{H}^3)\cap W^{\alpha,p}(0,T;(\mathbb{H}^{\beta})^*)
$
be a Banach space endowed with the norm
$
\|\mathbf{u}\|_{\mathcal{Y}(\beta)}=\|\mathbf{u}\|_{L^{\infty}(0,T;\mathbb{H}^1)}+\|\mathbf{u}\|_{ L^{2}(0,T;\mathbb{H}^3)}+\|\mathbf{u}\|_{W^{\alpha,p}(0,T;(\mathbb{H}^{\beta})^*)}.
$
By using \eqref{1the2} and noting that $\{\mathbf{u}_{0,k}\}_{k\in\mathbb{N}}$ is uniformly bounded in $\mathbb{H}^1$, it follows that for $p\geq4$, $\alpha\in(0,\frac{1}{2})$ and for all $k\in\mathbb{N}$,
$
\mathbf{u}^k\in L^p(\Omega;C([0,T];\mathbb{H}^1_w))
$ and $\mathbb{E}\|\mathbf{u}^k\|_{\mathcal{Y}(\beta)}\leq \hat{C}$, where the positive constant $\hat{C}$ only depends on $C_0,~p,~\mathbf{h}$ and $T$. Let
$
\mathcal{B}_R(\beta):=\{\mathbf{f}\in\mathcal{Y}(\beta):\|\mathbf{f}\|_{\mathcal{Y}(\beta)}\leq R\}.
$
Then by using the Chebyshev inequality, it follows that
\begin{equation}\label{1lem72}
\begin{split}
\sup_{k\in\mathbb{N}}\mathbb{P}\left(\mathbf{u}^k\in\mathcal{B}_R(\beta)\right)\geq1-\frac{\hat{C}}{R^2}.
\end{split}
\end{equation}
Moreover, for any $\beta_1\in(0,\beta)$, the following compact embedding is valid
\begin{equation*}
\begin{split}
W^{\alpha,p}(0,T;(\mathbb{H}^{\beta_1})^*)\cap L^p(0,T;\mathbb{H}^1)\cap L^2(0,T;\mathbb{H}^2)\hookrightarrow \mathcal{X}.
\end{split}
\end{equation*}
Thus $\mathcal{B}_R(\beta_1)$ is a compact subset in $\mathcal{X}$. As a consequence for a certain $R_1>0$, it follows that $\mathcal{B}_R(\beta_1)\subset\mathcal{X}\cap C([0,T];\mathbb{B}^1_w(R_1))\subset\tilde{\mathcal{X}}$ where $\mathbb{B}^1_w(R_1)$ was defined in section \ref{sec1}. Let $\{\mathbf{f}_n\}$ be sequence in $\mathcal{B}_R(\beta_1)$. Then there exists a subsequence of $\{\mathbf{f}_n\}$ (still denoted by $\{\mathbf{f}_n\}$) and $\mathbf{f}\in\mathcal{B}_R(\beta_1)$ such that
\begin{equation*}
\begin{split}
\mathbf{f}_n\rightarrow\mathbf{f}~\textrm{in}~\mathcal{X}~\textrm{and}~\mathbf{f}_n\rightarrow\mathbf{f}~\textrm{weak*}~\textrm{in}~L^{\infty}(0,T;\mathbb{H}^1).
\end{split}
\end{equation*}
Hence
$
\mathbf{f}_n\rightarrow\mathbf{f}~\textrm{in}~\mathcal{X}\cap C([0,T];\mathbb{B}^1_w(R_1))\subset\tilde{\mathcal{X}},
$
which implies that $\mathcal{B}_R(\beta_1)$ is a compact subset in $\tilde{\mathcal{X}}$. This together with \eqref{1lem72} and choosing $R>\sqrt{\frac{\hat{C}}{\varepsilon}}$, we see that the sequence $\{\mathbf{u}^k\}$ of $\tilde{\mathcal{X}}$-valued Borel random variables defined on $\left(\Omega^k,\mathcal{F}^k,\mathbb{F}^k,\mathbb{P}^k\right)$ satisfies the condition of Jakubowski-Skorokhod theorem. The proof is thus completed.~~$\Box$

\begin{proposition}\label{pro71} Let $\psi:\mathbb{H}^1\rightarrow\mathbb{R}$ be a bounded and sequentially weakly continuous function and let $\mathbf{u}_{0,k}\rightarrow\mathbf{u}_{0}$ weakly in $\mathbb{H}^1$ as $k\rightarrow\infty$. Then for every $t\geq0$,
$
P_t\psi(\mathbf{u}_{0,k})\rightarrow P_t\psi(\mathbf{u}_{0})~\textrm{as}~k\rightarrow\infty.
$
\end{proposition}
\noindent\textbf{Proof.} Thanks to Lemma \ref{lem72}, there exists a subsequence of $\mathbf{u}^k$ (still denoted by $\mathbf{u}^k$), a stochastic basis $\left(\tilde{\Omega},\tilde{\mathcal{F}},\tilde{\mathbb{F}},\tilde{\mathbb{P}}\right)$, an $\mathbb{R}^{\infty}$-valued standard $\tilde{\mathbb{F}}$-Wiener process $\tilde{W}=(\tilde{W}_j)_{j=1}^{\infty}$ defined on this basis, progressively measurable processes $\tilde{\mathbf{u}}$ and $\{\tilde{\mathbf{u}}^k\}$ (defined on this basis) with laws supported in $(\tilde{\mathcal{X}},\mathcal{T})$ such that
$
\tilde{\mathbf{u}}^k~\textrm{has the same law as}~\mathbf{u}^k~\textrm{on}~\tilde{\mathcal{X}}
$
and
$
\tilde{\mathbf{u}}^k\rightarrow\tilde{\mathbf{u}}~\textrm{in}~\tilde{\mathcal{X}}~\textrm{as}~k\rightarrow\infty,~\tilde{\mathbb{P}}\textrm{-a.s.}
$
Thus it follows that
\begin{equation}\label{1pro71}
\begin{split}
P_t\psi(\mathbf{u}_{0})=\mathbb{E}\psi(\mathbf{u}(t;\mathbf{u}_0))=\tilde{\mathbb{E}}\psi(\tilde{\mathbf{u}}(t)),
\end{split}
\end{equation}
and $\tilde{\mathbf{u}}^k\rightarrow\tilde{\mathbf{u}}$ in $C([0,T];\mathbb{H}^1_w)$, $\tilde{\mathbb{P}}$-a.s. Moreover, by using the sequential weak continuity of $\psi$, we see that
$
\psi(\tilde{\mathbf{u}}^k(t))\rightarrow\psi(\tilde{\mathbf{u}}(t))~\textrm{in}~\mathbb{R}~\textrm{as}~k\rightarrow\infty.
$
Since the function $\psi$ is bounded, by the Lebesgue dominated convergence theorem we infer that
\begin{equation}\label{2pro71}
\begin{split}
\lim_{k\rightarrow\infty}\tilde{\mathbb{E}}\psi(\tilde{\mathbf{u}}^k(t))=\tilde{\mathbb{E}}\psi(\tilde{\mathbf{u}}(t)).
\end{split}
\end{equation}
Hence, combining \eqref{1pro71} with \eqref{2pro71} and noting the laws of $\mathbf{u}^k$ and $\tilde{\mathbf{u}}^k$ are equivalent, we derive that
\begin{equation*}
\begin{split}
\lim_{k\rightarrow\infty}P_t\psi(\mathbf{u}_{0,k})=\lim_{k\rightarrow\infty}\mathbb{E}\psi(\mathbf{u}^k)=\lim_{k\rightarrow\infty}\tilde{\mathbb{E}}\psi(\tilde{\mathbf{u}}^k)=\tilde{\mathbb{E}}\psi(\tilde{\mathbf{u}})=P_t\psi(\mathbf{u}_{0}).
\end{split}
\end{equation*}
The proof is completed.~~$\Box$

\begin{theorem}\label{the71}  Under the same assumptions as in Proposition \ref{pro61}, there exists an invariant measure $\mu$ of the semigroup $\{P_t\}$ defined by \eqref{ad71}, such that for any $t\geq0$ and $\psi\in C_b(\mathbb{H}^1_w)$,
$
\int_{\mathbb{H}^1}P_t\psi(\mathbf{u})\mu(\mathrm{d}\mathbf{u})=\int_{\mathbb{H}^1}\psi(\mathbf{u})\mu(\mathrm{d}\mathbf{u}).
$
\end{theorem}
\noindent\textbf{Proof.} The proof is based on the Maslowski-Seidler theorem \cite{33}, see also Lemma \ref{lem96}. Proposition \ref{pro71} means that the semigroup $\{P_t\}$ is sequentially weakly Feller in $\mathbb{H}^1$. Additionally, by using the Chebyshev inequality and Lemma \ref{lem71}, we derive that for every $T>0$ and $R>0$,
$
\frac{1}{T}\int_0^T\mathbb{P}\left(\|\mathbf{u}(s;\mathbf{u}_0)\|_{\mathbb{H}^1}>R\right)\,\mathrm{d}s\leq\frac{1}{TR^2}\int_0^T\mathbb{E}\|\mathbf{u}(s;\mathbf{u}_0)\|_{\mathbb{H}^1}^2\,\mathrm{d}s\leq\frac{C+CT}{TR^2}.
$
Let $T\geq1$ be fixed. For any $\varepsilon>0$  there exists $R>\sqrt{\frac{2C}{\varepsilon}}$ such that
\begin{equation*}
\begin{split}
\frac{1}{T}\int_0^T\mathbb{P}\left(\|\mathbf{u}(s;\mathbf{u}_0)\|_{\mathbb{H}^1}>R\right)\,\mathrm{d}s\leq\varepsilon.
\end{split}
\end{equation*}
Thus by Lemma \ref{lem96}, we infer that there exists as least one invariant measure for equation \eqref{sys1}.~~$\Box$\\
\textbf{Proof of Theorem \ref{the2}.}
Theorem \ref{the2} is a direct combination of Proposition \ref{pro61} and Theorem \ref{the71}.~~$\Box$

\section{Proof of Theorem \ref{the3}}\label{sec7}
This section is dedicated to proving the final theorem of this paper. The proof framework of Theorem \ref{the3} is similar to that of Theorems \ref{the1} and Theorem \ref{the2}. To avoid repetition, the proof processes of certain lemmas in this section will be omitted.

\subsection{Pathwise very weak solution}\label{1subsec8}
We now consider the following Galerkin approximation system which does not require the introduction of the truncation function $\theta_R(\cdot)$.
\begin{equation}\label{sys4}
\left\{
\begin{aligned}
&\mathrm{d}\mathbf{u}_n=\bigl[\beta_1\Delta \mathbf{u}_n-\beta_2\Delta^2\mathbf{u}_n+\beta_3\Pi_n\left((1-|\mathbf{u}_n|^2)\mathbf{u}_n\right)-\beta_4\Pi_n(\mathbf{u}_n\times\Delta\mathbf{u}_n)\\
&~~+\beta_5\Pi_n\Delta(|\mathbf{u}_n|^2\mathbf{u}_n)\bigl]\,\mathrm{d}t+\sum_{j=1}^{n}\Pi_n(-\mathbf{u}_n\times \mathbf{h}_j+\mathbf{h}_j-\Delta\mathbf{h}_j)\circ\mathrm{d}W_j(t),&&\textrm{in}~(0,\infty)\times\mathcal{O},\\
&\mathbf{u}_n(0)=\Pi_n\mathbf{u}_0,&&\textrm{in}~\mathcal{O}.
\end{aligned}
\right.
\end{equation}
\begin{lemma}\label{lem81} Let $\mathcal{O}\subset\mathbb{R}^d,~d=1,2,3$, be a bounded domain with $C^{1,1}$-boundary. Then for any $p\geq1$, $n\in\mathbb{N}$ and every $t\in[0,T]$, there exists a positive constant $C=C(\|\mathbf{u}_0\|_{\mathbb{L}^2},p,\mathbf{h},T)$ independent of $n$ such that
\begin{equation}\label{1lem61}
\begin{split}
&\mathbb{E}\sup_{s\in[0,t]}\|\mathbf{u}_n(s)\|_{\mathbb{L}^2}^{2p}+\mathbb{E}\left(\int_0^t\| \mathbf{u}_n(s)\|_{\mathbb{H}^2}^2\,\mathrm{d}s\right)^p+\mathbb{E}\left(\int_0^t\|\mathbf{u}_n(s)\|_{\mathbb{L}^4}^4\,\mathrm{d}s\right)^p\leq C.
\end{split}
\end{equation}
\end{lemma}
\noindent\textbf{Proof.} The proof process is standard, one can see Lemma \ref{lem32}.~~$\Box$

\begin{lemma}\label{lem82} Let $\mathcal{O}\subset\mathbb{R}^d,~d=1,2$, be a bounded domain with $C^{1,1}$-boundary. Let $q\geq1$, $p>2$ and $\alpha\in(0,\frac{1}{2})$ with $p\alpha>1$. Then for any $n\in\mathbb{N}$ and every $t\in[0,T]$, there exists a positive constant $C$ independent of $n$ such that
\begin{equation}\label{0lem82}
\begin{split}
\mathbb{E}\|\mathbf{u}_n\|_{W^{\alpha,p}(0,T;(\mathbb{H}^{2})^*)}^{q}\leq C.
\end{split}
\end{equation}
\end{lemma}
\noindent\textbf{Proof.} Equation \eqref{sys4} can be written in the following way:
\begin{equation}\label{1lem82}
\begin{split}
\mathbf{u}_n(t)&=\mathbf{u}_n(0)+\beta_1\int_0^tF_n^1(\mathbf{u}_n)\,\mathrm{d}s-\beta_2\int_0^tF_n^2(\mathbf{u}_n)\,\mathrm{d}s+\beta_3\int_0^t\Pi_n\mathbf{u}_n\,\mathrm{d}s-\beta_3\int_0^tF_n^3(\mathbf{u}_n)\,\mathrm{d}s\\
&-\beta_4\int_0^tF_n^4(\mathbf{u}_n)\,\mathrm{d}s+\beta_5\int_0^t\Pi_n\Delta\left(|\mathbf{u}_n|^2\mathbf{u}_n\right)\,\mathrm{d}s-\frac{1}{2}\sum_{j=1}^n\int_0^t\Pi_n\left(G_{nj}(\mathbf{u}_n)\times\mathbf{h}_j\right)\,\mathrm{d}s\\
&+\sum_{j=1}^n\int_0^tG_{nj}(\mathbf{u}_n)\mathrm{d}W_j(t)\\
&:=\mathbf{u}_n(0)+\sum_{k=1}^7\mathbf{B}_{n,k}(\mathbf{u}_n)(t)+\mathbf{B}_{n,8}(\mathbf{u}_n,W_n)(t),~t\in[0,T].
\end{split}
\end{equation}
Let $\phi\in\mathbb{H}^2$. By the Sobolev embeddings $\mathbb{H}^2\hookrightarrow \mathbb{L}^{\infty}$ and $\mathbb{H}^1\hookrightarrow \mathbb{L}^{6}$, it follows that
\begin{equation}\label{2lem82}
\begin{split}
&|\left(F_n^1(\mathbf{u}_n),\phi\right)_{\mathbb{L}^2}|=|\left(\mathbf{u}_n,\Delta\phi\right)_{\mathbb{L}^2}|\leq\|\mathbf{u}_n\|_{\mathbb{L}^2}\|\phi\|_{\mathbb{H}^2},\\
&|\left(F_n^2(\mathbf{u}_n),\phi\right)_{\mathbb{L}^2}|=|\left(\Delta\mathbf{u}_n,\Delta\phi\right)_{\mathbb{L}^2}|\leq\|\Delta\mathbf{u}_n\|_{\mathbb{L}^2}\|\phi\|_{\mathbb{H}^2},\\
&|\left(F_n^3(\mathbf{u}_n),\phi\right)_{\mathbb{L}^2}|\leq\|\mathbf{u}_n\|_{\mathbb{L}^3}^3\|\phi\|_{\mathbb{L}^{\infty}}\leq C\|\mathbf{u}_n\|_{\mathbb{L}^3}^3\|\phi\|_{\mathbb{H}^{2}},\\
&|\left(F_n^4(\mathbf{u}_n),\phi\right)_{\mathbb{L}^2}|\leq|\left(\mathbf{u}_n\times\nabla\mathbf{u}_n,\nabla\phi\right)_{\mathbb{L}^2}|\leq C\|\mathbf{u}_n\|_{\mathbb{L}^2}\|\nabla\mathbf{u}_n\|_{\mathbb{L}^3}\|\nabla\phi\|_{\mathbb{L}^6}\leq C\|\mathbf{u}_n\|_{\mathbb{L}^2}\|\mathbf{u}_n\|_{\mathbb{H}^2}\|\phi\|_{\mathbb{H}^2},\\
&|\left(\Pi_n\Delta\left(|\mathbf{u}_n|^2\mathbf{u}_n\right),\phi\right)_{\mathbb{L}^2}|\leq|\left(|\mathbf{u}_n|^2\mathbf{u}_n,\Delta\phi\right)_{\mathbb{L}^2}|\leq C\|\mathbf{u}_n\|_{\mathbb{L}^6}^3\|\phi\|_{\mathbb{H}^2}.
\end{split}
\end{equation}
According to the GN inequality, it follows that
\begin{equation}\label{3lem82}
\begin{split}
&\|\mathbf{u}_n\|_{\mathbb{L}^3}\leq C\|\mathbf{u}_n\|_{\mathbb{L}^2}^{\frac{11}{12}}\|\mathbf{u}_n\|_{\mathbb{H}^2}^{\frac{1}{12}},~\textrm{for}~d=1,\\
&\|\mathbf{u}_n\|_{\mathbb{L}^3}\leq C\|\mathbf{u}_n\|_{\mathbb{L}^2}^{\frac{5}{6}}\|\mathbf{u}_n\|_{\mathbb{H}^2}^{\frac{1}{6}},~\textrm{for}~d=2,\\
&\|\mathbf{u}_n\|_{\mathbb{L}^6}\leq C\|\mathbf{u}_n\|_{\mathbb{L}^2}^{\frac{5}{6}}\|\mathbf{u}_n\|_{\mathbb{H}^2}^{\frac{1}{6}},~\textrm{for}~d=1,\\
&\|\mathbf{u}_n\|_{\mathbb{L}^6}\leq C\|\mathbf{u}_n\|_{\mathbb{L}^2}^{\frac{2}{3}}\|\mathbf{u}_n\|_{\mathbb{H}^2}^{\frac{1}{3}},~\textrm{for}~d=2,\\
\end{split}
\end{equation}
which implies that the following estimates
\begin{equation}\label{4lem82}
\begin{split}
&\|\mathbf{u}_n\|_{\mathbb{L}^3}^6\leq C+C\|\mathbf{u}_n\|_{\mathbb{L}^2}^{10}+C\|\mathbf{u}_n\|_{\mathbb{H}^2}^{2},\\
&\|\mathbf{u}_n\|_{\mathbb{L}^6}^6\leq C\|\mathbf{u}_n\|_{\mathbb{L}^2}^{5}+C\|\mathbf{u}_n\|_{\mathbb{H}^2}^{2}+C\|\mathbf{u}_n\|_{\mathbb{L}^2}^{5}\|\mathbf{u}_n\|_{\mathbb{H}^2}^{2},\\
\end{split}
\end{equation}
are valid for $d=1,2$. Thus by \eqref{2lem82}-\eqref{4lem82}, we infer that
\begin{equation}\label{5lem82}
\begin{split}
&\sum_{k=1}^7\mathbb{E}\|\mathbf{B}_{n,k}(\mathbf{u}_n)\|_{W^{1,2}(0,T;(\mathbb{H}^{2})^*)}^q\\
&\leq C+C\mathbb{E}\left(\int_0^T\|\mathbf{u}_n\|_{\mathbb{L}^{2}}^{10}\,\mathrm{d}s\right)^{\frac{q}{2}}+C\mathbb{E}\left(\int_0^T\|\mathbf{u}_n\|_{\mathbb{H}^{2}}^2\,\mathrm{d}s\right)^{\frac{q}{2}}+C\mathbb{E}\left(\int_0^T\|\mathbf{u}_n\|_{\mathbb{L}^{2}}^2\|\mathbf{u}_n\|_{\mathbb{H}^{1}}^2\,\mathrm{d}s\right)^{\frac{q}{2}}\\
&+C\mathbb{E}\left(\int_0^T\|\mathbf{u}_n\|_{\mathbb{L}^{2}}^5\|\mathbf{u}_n\|_{\mathbb{H}^{2}}^2\,\mathrm{d}s\right)^{\frac{q}{2}}\\
&\leq C+C\mathbb{E}\left[\sup_{t\in[0,T]}\|\mathbf{u}_n(t)\|_{\mathbb{L}^{2}}^{\frac{5q}{2}}\left(\int_0^T\|\mathbf{u}_n\|_{\mathbb{H}^{2}}^2\,\mathrm{d}s\right)^{\frac{q}{2}}\right]\\
&\leq C+C\mathbb{E}\sup_{t\in[0,T]}\|\mathbf{u}_n(t)\|_{\mathbb{L}^{2}}^{5q}+C\mathbb{E}\left(\int_0^T\|\mathbf{u}_n\|_{\mathbb{H}^{2}}^2\,\mathrm{d}s\right)^{q}\leq C.
\end{split}
\end{equation}
Moreover, by \eqref{11cor32} and the embedding $(\mathbb{H}^{1})^*\hookrightarrow(\mathbb{H}^{2})^*$, it follows that
\begin{equation}\label{6lem82}
\begin{split}
&\mathbb{E}\|\mathbf{B}_{n,8}(\mathbf{u}_n,W)\|_{W^{\alpha,p}(0,T;(\mathbb{H}^{2})^*)}^q\leq C_{\mathbf{h}}+C_{\mathbf{h}}\mathbb{E}\left(\int_0^T\|\mathbf{u}_n\|_{\mathbb{L}^2}^q\,\mathrm{d}s\right)\leq C.
\end{split}
\end{equation}
Since
$
W^{1,2}(0,T;(\mathbb{H}^{2})^*)\hookrightarrow W^{\alpha,p}(0,T;(\mathbb{H}^{2})^*),~\textrm{if}~\frac{1}{2}+\frac{1}{p}>\alpha,
$
the inequality \eqref{0lem82} can be directly obtained by using the estimates \eqref{5lem82} and \eqref{6lem82}.~~$\Box$

\begin{lemma}\label{lem83} Let $\mathcal{O}\subset\mathbb{R}^3$, be a bounded domain with $C^{1,1}$-boundary. Let $r\in[1,\frac{4}{3})$, $p>2$ and $\alpha\in(0,\frac{1}{2})$ with $p\alpha>1$. Then for any $n\in\mathbb{N}$ and every $t\in[0,T]$, there exists a positive constant $C$ independent of $n$ such that
\begin{equation}\label{0lem83}
\begin{split}
\mathbb{E}\|\mathbf{u}_n\|_{W^{\alpha,p}(0,T;(\mathbb{H}^{2})^*)}^{q}\leq C,~\textrm{if}~\frac{1}{p}-\alpha>\frac{1}{r}-1.
\end{split}
\end{equation}
\end{lemma}
\noindent\textbf{Proof.} According to the GN inequality, it follows that
\begin{equation}\label{1lem83}
\begin{split}
&\|\mathbf{u}_n\|_{\mathbb{L}^3}\leq C\|\mathbf{u}_n\|_{\mathbb{L}^2}^{\frac{3}{4}}\|\mathbf{u}_n\|_{\mathbb{H}^2}^{\frac{1}{4}},~\textrm{for}~d=3,\\
&\|\mathbf{u}_n\|_{\mathbb{L}^6}\leq C\|\mathbf{u}_n\|_{\mathbb{L}^2}^{\frac{1}{2}}\|\mathbf{u}_n\|_{\mathbb{H}^2}^{\frac{1}{2}},~\textrm{for}~d=3.
\end{split}
\end{equation}
Thus by using the H\"{o}lder inequality and Young's inequality, it follows from \eqref{2lem82} and \eqref{1lem83} that
\begin{equation*}
\begin{split}
&\sum_{k=1}^4\mathbb{E}\|F_n^k(\mathbf{u}_n)\|_{L^2(0,T;(\mathbb{H}^{2})^*)}^q\leq C
\end{split}
\end{equation*}
and
\begin{equation*}
\begin{split}
&\mathbb{E}\|\Pi_n\Delta\left(|\mathbf{u}_n|^2\mathbf{u}_n\right)\|_{L^r(0,T;(\mathbb{H}^{2})^*)}^q\leq C\mathbb{E}\left(\int_0^T\|\mathbf{u}_n\|_{\mathbb{L}^2}^{\frac{3r}{2}}\|\mathbf{u}_n\|_{\mathbb{H}^2}^{\frac{3r}{2}}\,\mathrm{d}s\right)^{\frac{q}{r}}\\
&\leq C\mathbb{E}\left[\left(\int_0^T\|\mathbf{u}_n\|_{\mathbb{L}^2}^{\frac{6r}{4-3r}}\,\mathrm{d}s\right)^{\frac{4-3r}{4}}\left(\int_0^T\|\mathbf{u}_n\|_{\mathbb{H}^2}^{2}\,\mathrm{d}s\right)^{\frac{3r}{4}}\right]^{\frac{q}{r}}\\
&\leq C\mathbb{E}\sup_{t\in[0,T]}\|\mathbf{u}_n(t)\|_{\mathbb{L}^{2}}^{3q}+C\mathbb{E}\left(\int_0^T\|\mathbf{u}_n\|_{\mathbb{H}^{2}}^2\,\mathrm{d}s\right)^{\frac{3q}{2}}\leq C,
\end{split}
\end{equation*}
which together with the Sobolev embedding
$
W^{1,r}(0,T;(\mathbb{H}^{2})^*)\hookrightarrow W^{\alpha,p}(0,T;(\mathbb{H}^{2})^*),~\textrm{if}~\frac{1}{p}-\alpha>\frac{1}{r}-1,
$
yields the result \eqref{0lem83}.~~$\Box$

Thanks to Lemmas \ref{lem81} and \ref{lem82}, it is not hard to obtain the following tightness result.
\begin{lemma}\label{lem84}
Let $\mathcal{O}\subset\mathbb{R}^d,~d=1,2,3$, be a bounded domain with $C^{1,1}$-boundary. If $\beta'>2$ and $p>1$, then the measures $\{\mathcal{L}(\mathbf{u}_n)\}_{n\in\mathbb{N}}$ on $\mathcal{X}_1:=C([0,T];(\mathbb{H}^{\beta'})^*)\cap L^2(0,T;\mathbb{H}^1)\cap L^4_w(0,T;\mathbb{L}^4)$ are tight.
\end{lemma}

By a similar argument as shown in section \ref{sec3}, we have the following conclusions.

\begin{proposition}\label{pro81}
There exist
\begin{enumerate}
\item [(1)] a probability space $(\Omega',\mathcal{F}',\mathbb{P}')$,
\item [(2)] a sequence $\{(\mathbf{u}_n',W_n')\}$ of random variables defined on $(\Omega',\mathcal{F}',\mathbb{P}')$ and taking values in the space $\mathcal{X}_1\times C([0,T];\mathbb{R}^{\infty})$,
\item [(3)] a random variable $(\mathbf{u}',W')$ defined on $(\Omega',\mathcal{F}',\mathbb{P}')$ and taking values in the space $\mathcal{X}_1\times C([0,T];\mathbb{R}^{\infty})$,
\end{enumerate}
such that in the space $\mathcal{X}_1\times C([0,T];\mathbb{R}^{\infty})$ there hold
\begin{enumerate}
\item [(a)] $\mathcal{L}(\mathbf{u}_n,W_n)=\mathcal{L}(\mathbf{u}_n',W_n')$,
\item [(b)] $(\mathbf{u}_n',W_n')$ converges $\mathbb{P}'$-almost surely to $(\mathbf{u}',W')$ in the topology of $\mathcal{X}_1$.
\end{enumerate}
\end{proposition}

\begin{corollary}\label{cor81}
For every $p\geq1$ and $T>0$,
\begin{equation*}\label{1cor81}
\begin{split}
\sup_{n\in\mathbb{N}}\left[\mathbb{E}\sup_{t\in[0,T]}\|\mathbf{u}_n'(t)\|_{\mathbb{L}^2}^{2p}+\mathbb{E}\left(\int_0^T\|\mathbf{u}_n'(t)\|_{\mathbb{H}^2}^2\,\mathrm{d}t\right)^p+\mathbb{E}\left(\int_0^T\|\mathbf{u}_n'(t)\|_{\mathbb{L}^4}^4\,\mathrm{d}t\right)^p\right]<\infty.
\end{split}
\end{equation*}
\end{corollary}

Similar to the proof of Lemma \ref{lem54}, we have the following results.
\begin{lemma}\label{lem85}
For any $p\geq1$ there holds
\begin{equation*}
\begin{split}
\mathbf{u}_n'\rightarrow\mathbf{u}'~\textrm{weakly~in}~L^{2p}(\Omega';L^{\infty}(0,T;\mathbb{L}^2)\cap L^2(0,T;\mathbb{H}^2))\cap L^{4p}(\Omega';L^4(0,T;\mathbb{L}^4)).
\end{split}
\end{equation*}
Moreover, the process $\mathbf{u}'\in L^{2p}(\Omega';C([0,T];\mathbb{L}^2_w))$.
\end{lemma}

Furthermore, we have the following convergence results.
\begin{lemma}\label{lem86}
Let $\mathcal{O}\subset\mathbb{R}^d,~d=1,2,3$. For any $\phi\in \mathbb{H}^2$, there hold
\begin{equation*}\label{5lem86}
\begin{split}
\lim_{n\rightarrow\infty}\mathbb{E}\int_0^t\left(\Delta \mathbf{u}_n',\phi\right)_{\mathbb{L}^2}\,\mathrm{d}s=-\mathbb{E}\int_0^t\left(\nabla\mathbf{u}',\nabla\phi\right)_{\mathbb{L}^2}\,\mathrm{d}s,
\end{split}
\end{equation*}
\begin{equation*}\label{6lem86}
\begin{split}
\lim_{n\rightarrow\infty}\mathbb{E}\int_0^t\left(\Delta^2 \mathbf{u}_n',\phi\right)_{\mathbb{L}^2}\,\mathrm{d}s=\mathbb{E}\int_0^t\left(\Delta\mathbf{u}',\Delta\phi\right)_{\mathbb{L}^2}\,\mathrm{d}s,
\end{split}
\end{equation*}
\begin{equation*}\label{1lem86}
\begin{split}
\lim_{n\rightarrow\infty}\mathbb{E}\int_0^t\left(\Pi_n\left(\left(1-|\mathbf{u}_n'|^2\right)\mathbf{u}_n'\right),\phi\right)_{\mathbb{L}^2}\,\mathrm{d}s=\mathbb{E}\int_0^t\left(\left(1-|\mathbf{u}'|^2\right)\mathbf{u}',\phi\right)_{\mathbb{L}^2}\,\mathrm{d}s,
\end{split}
\end{equation*}
\begin{equation*}\label{2lem86}
\begin{split}
\lim_{n\rightarrow\infty}\mathbb{E}\int_0^t\left(\Pi_n\left(\left(\mathbf{u}_n'\times\mathbf{h}_j\right)\times\mathbf{h}_j\right),\phi\right)_{\mathbb{L}^2}\,\mathrm{d}s=\mathbb{E}\int_0^t\left(\left(\left(\mathbf{u}'\times\mathbf{h}_j\right)\times\mathbf{h}_j\right),\phi\right)_{\mathbb{L}^2}\,\mathrm{d}s
\end{split}
\end{equation*}
\begin{equation*}\label{3lem86}
\begin{split}
\lim_{n\rightarrow\infty}\mathbb{E}\int_0^t\left(\Pi_n\left(\mathbf{u}_n'\times\Delta \mathbf{u}_n'\right),\phi\right)_{\mathbb{L}^2}\,\mathrm{d}s=-\mathbb{E}\int_0^t\left(\mathbf{u}'\times\nabla\mathbf{u}',\nabla\phi\right)_{\mathbb{L}^2}\,\mathrm{d}s,
\end{split}
\end{equation*}
\begin{equation}\label{4lem86}
\begin{split}
\lim_{n\rightarrow\infty}\mathbb{E}\int_0^t\left(\Pi_n\left(\Delta\left(|\mathbf{u}_n'|^2\mathbf{u}_n'\right)\right),\phi\right)_{\mathbb{L}^2}\,\mathrm{d}s=\mathbb{E}\int_0^t\left(|\mathbf{u}'|^2\mathbf{u}',\Delta\phi\right)_{\mathbb{L}^2}\,\mathrm{d}s
\end{split}
\end{equation}
\end{lemma}
\noindent\textbf{Proof.} We only provide a proof of \eqref{4lem86}, as the others are similar. To show \eqref{4lem86}, it is sufficient to prove that
\begin{equation*}
\begin{split}
\lim_{n\rightarrow\infty}\mathbb{E}\int_0^t\left(|\mathbf{u}'_n|^2\mathbf{u}'_n,\Delta\phi\right)_{\mathbb{L}^2}\,\mathrm{d}s=\mathbb{E}\int_0^t\left(|\mathbf{u}'|^2\mathbf{u}',\Delta\phi\right)_{\mathbb{L}^2}\,\mathrm{d}s.
\end{split}
\end{equation*}
By the H\"{o}lder inequality and Sobolev embedding $\mathbb{H}^1\hookrightarrow \mathbb{L}^{6}$, it follows that
\begin{equation}\label{7lem86}
\begin{split}
&\left|\int_0^T\left(|\mathbf{u}'_n|^2\mathbf{u}'_n,\Delta\phi\right)_{\mathbb{L}^2}\,\mathrm{d}s-\int_0^T\left(|\mathbf{u}'|^2\mathbf{u}',\Delta\phi\right)_{\mathbb{L}^2}\,\mathrm{d}s\right|\\
&\leq\left|\int_0^T\left(|\mathbf{u}'_n|^2\left(\mathbf{u}'_n-\mathbf{u}'\right),\Delta\phi\right)_{\mathbb{L}^2}\,\mathrm{d}s\right|+\left|\int_0^T\left(\left(|\mathbf{u}'_n|^2-|\mathbf{u}'|^2\right)\mathbf{u}',\Delta\phi\right)_{\mathbb{L}^2}\,\mathrm{d}s\right|\\
&\leq\|\Delta\phi\|_{\mathbb{L}^2}\|\mathbf{u}_n'\|_{L^4(0,T;\mathbb{L}^6)}^2\|\mathbf{u}_n'-\mathbf{u}'\|_{L^2(0,T;\mathbb{L}^6)}\\
&+\|\Delta\phi\|_{\mathbb{L}^2}\|(|\mathbf{u}'_n|+|\mathbf{u}'|)\mathbf{u}'\|_{L^2(0,T;\mathbb{L}^3)}\|\mathbf{u}_n'-\mathbf{u}'\|_{L^2(0,T;\mathbb{L}^6)}\\
&\leq C\|\phi\|_{\mathbb{H}^2}\left(\|\mathbf{u}_n'\|_{L^4(0,T;\mathbb{L}^6)}^2+\|\mathbf{u}'\|_{L^4(0,T;\mathbb{L}^6)}^2\right)\|\mathbf{u}_n'-\mathbf{u}'\|_{L^2(0,T;\mathbb{H}^1)}.
\end{split}
\end{equation}
Let us now recall the estimate \eqref{3lem82} and \eqref{1lem83}, then we have
\begin{equation*}
\begin{split}
\|f\|_{L^4(0,T;\mathbb{L}^6)}^4\leq\int_0^T\|f\|_{\mathbb{L}^2}^{\frac{10}{3}}\|f\|_{\mathbb{H}^2}^{\frac{2}{3}}\,\mathrm{d}s~\textrm{for}~d=1,\\
\|f\|_{L^4(0,T;\mathbb{L}^6)}^4\leq\int_0^T\|f\|_{\mathbb{L}^2}^{\frac{8}{3}}\|f\|_{\mathbb{H}^2}^{\frac{4}{3}}\,\mathrm{d}s~\textrm{for}~d=2,\\
\|f\|_{L^4(0,T;\mathbb{L}^6)}^4\leq\int_0^T\|f\|_{\mathbb{L}^2}^{2}\|f\|_{\mathbb{H}^2}^{2}\,\mathrm{d}s~\textrm{for}~d=3,
\end{split}
\end{equation*}
which implies that
\begin{equation}\label{8lem86}
\begin{split}
\|f\|_{L^4(0,T;\mathbb{L}^6)}^4\leq C+C\sup_{s\in[0,T]}\|f(s)\|_{\mathbb{L}^2}^4\int_0^T\|f\|_{\mathbb{H}^2}^{2}\,\mathrm{d}s+C\sup_{s\in[0,T]}\|f(s)\|_{\mathbb{L}^2}^4+C\int_0^T\|f\|_{\mathbb{H}^2}^{2}\,\mathrm{d}s
\end{split}
\end{equation}
is valid for $d=1,2,3$. This together with \eqref{7lem86}, Proposition \ref{pro81} and Lemma \ref{lem85} implies that
\begin{equation*}
\begin{split}
&\left|\int_0^T\left(|\mathbf{u}'_n|^2\mathbf{u}'_n,\Delta\phi\right)_{\mathbb{L}^2}\,\mathrm{d}s-\int_0^T\left(|\mathbf{u}'|^2\mathbf{u}',\Delta\phi\right)_{\mathbb{L}^2}\,\mathrm{d}s\right|\\
&\leq C\|\phi\|_{\mathbb{H}^2}\Bigl(1+\|\mathbf{u}_n'\|_{L^{\infty}(0,T;\mathbb{L}^2)}^2\|\mathbf{u}_n'\|_{L^{2}(0,T;\mathbb{H}^2)}+\|\mathbf{u}_n'\|_{L^{\infty}(0,T;\mathbb{L}^2)}^2+\|\mathbf{u}_n'\|_{L^{2}(0,T;\mathbb{H}^2)}\\
&+\|\mathbf{u}'\|_{L^{\infty}(0,T;\mathbb{L}^2)}^2\|\mathbf{u}'\|_{L^{2}(0,T;\mathbb{H}^2)}+\|\mathbf{u}'\|_{L^{\infty}(0,T;\mathbb{L}^2)}^2+\|\mathbf{u}'\|_{L^{2}(0,T;\mathbb{H}^2)}\Bigl)\|\mathbf{u}_n'-\mathbf{u}'\|_{L^2(0,T;\mathbb{H}^1)}\\
&\rightarrow0~\textrm{as}~n\rightarrow\infty.
\end{split}
\end{equation*}
Moreover, thanks to the H\"{o}lder inequality, the inequality \eqref{8lem86} and Corollary \ref{cor81}, we have
\begin{equation*}
\begin{split}
&\sup_{n\in\mathbb{N}}\mathbb{E}\left|\int_0^T\left(|\mathbf{u}'_n|^2\mathbf{u}'_n,\Delta\phi\right)_{\mathbb{L}^2}\,\mathrm{d}s\right|^{\frac{4}{3}}\leq \sup_{n\in\mathbb{N}}C\|\Delta\phi\|_{\mathbb{L}^2}^{\frac{4}{3}}\mathbb{E}\int_0^T\|\mathbf{u}'_n\|_{\mathbb{L}^6}^4\,\mathrm{d}s\\
&\leq C\|\Delta\phi\|_{\mathbb{L}^2}^{\frac{4}{3}}+\sup_{n\in\mathbb{N}}C\|\Delta\phi\|_{\mathbb{L}^2}^{\frac{4}{3}}\Bigl(\mathbb{E}\sup_{s\in[0,T]}\|\mathbf{u}'_n(s)\|_{\mathbb{L}^2}^4\int_0^T\|\mathbf{u}'_n\|_{\mathbb{H}^2}^{2}\,\mathrm{d}s\\
&+\mathbb{E}\sup_{s\in[0,T]}\|\mathbf{u}'_n(s)\|_{\mathbb{L}^2}^4+\mathbb{E}\int_0^T\|\mathbf{u}'_n\|_{\mathbb{H}^2}^{2}\,\mathrm{d}s\Bigl)<\infty.
\end{split}
\end{equation*}
This together with the Vitali theorem yields \eqref{4lem86}. The proof is thus completed.~~$\Box$

By a similar argument to the proof of \eqref{1cor555-1} and \eqref{2cor555-1}, it is not difficult to prove that
\begin{equation*}
\begin{split}
\mathbf{u}'(t)&=\mathbf{u}'(0)+\beta_1\int_0^t\Delta \mathbf{u}'(s)\,\mathrm{d}s-\beta_2\int_0^t\Delta^2 \mathbf{u}'(s)\,\mathrm{d}s+\beta_3\int_0^t\left(1-|\mathbf{u}'(s)|^2\mathbf{u}'(s)\right)\,\mathrm{d}s\\
&-\beta_4\int_0^t\mathbf{u}'(s)\times\Delta \mathbf{u}'(s)\,\mathrm{d}s+\beta_5\int_0^t\Delta \left(|\mathbf{u}'(s)|^2\mathbf{u}'(s)\right)\,\mathrm{d}s\\
&-\frac{1}{2}\sum_{j=1}^{\infty}\int_0^t\left(\left(\mathbf{u}'(s)\times\mathbf{h}_j+\mathbf{h}_j+\Delta\mathbf{h}_j\right)\times\mathbf{h}_j\right)\,\mathrm{d}s+\sum_{j=1}^{\infty}\int_0^t-\mathbf{u}'(s)\times\mathbf{h}_j+\mathbf{h}_j-\Delta\mathbf{h}_j\,\mathrm{d}W_{j}'(s)
\end{split}
\end{equation*}
in $(\mathbb{H}^{2})^*$. Thus $\left(\Omega',\mathcal{F}',\mathbb{P}',W',\mathbf{u}'\right)$ is a global martingale very weak solution of \eqref{sys1}.

Particularly, in the case of $d=1,2$, the global very weak solution is actually pathwise unique.

\begin{proposition}\label{pro82} Let $\mathcal{O}\subset\mathbb{R}^d,~d=1,2$ and let $\mathbf{u}_0\in \mathbb{L}^2$ be fixed. Assume that $\left(\Omega,\mathcal{F},\mathbb{P},W,\mathbf{u}_1\right)$ and $\left(\Omega,\mathcal{F},\mathbb{P},W,\mathbf{u}_2\right)$ are two martingale very weak solution of \eqref{sys1} such that for $i=1,2$,
$$\mathbf{u}_i(0)=\mathbf{u}(0);~\mathbf{u}_i\in L^{\infty}(0,T;\mathbb{L}^2)\cap L^2(0,T;\mathbb{H}^2);~\mathbf{u}_i~\textrm{satisfies equation}~ \eqref{sys1}.
$$
Then $\mathbf{u}_1=\mathbf{u}_2$, $\mathbb{P}$-a.s. Therefore, \eqref{sys1} admits a unique global pathwise very weak solution.
\end{proposition}
\noindent\textbf{Proof.} Let $\mathbf{u}^*:=\mathbf{u}_1-\mathbf{u}_2$. Let
\begin{equation*}
\begin{split}
\chi^K:=\inf\left\{t\geq0:\|\mathbf{u}_1(t)\|_{\mathbb{L}^2}^2+\|\mathbf{u}_2(t)\|_{\mathbb{L}^2}^2+\int_0^t\|\mathbf{u}_1\|_{\mathbb{H}^2}^2\,\mathrm{d}s+\int_0^t\|\mathbf{u}_2\|_{\mathbb{H}^2}^2\,\mathrm{d}s>K\right\}\wedge T,~K>0.
\end{split}
\end{equation*}
Due to Lemma \ref{lem85}, it follows that
$
\chi^K\nearrow T~\textrm{as}~K\rightarrow\infty,~\mathbb{P}\textrm{-a.s.}
$
By a similar argument as shown in Proposition \ref{pro51}, it follows that
\begin{equation*}
\begin{split}
\|\mathbf{u}^*(t\wedge\chi^K)\|_{\mathbb{L}^2}^2\leq C\int_0^{t\wedge\chi^K}\mathbf{F}(s)\|\mathbf{u}^*(s)\|_{\mathbb{L}^{2}}^2\,\mathrm{d}s,
\end{split}
\end{equation*}
where $\mathbf{F}=1+\|\mathbf{u}_1\|_{\mathbb{L}^{\infty}}^4+\|\mathbf{u}_2\|_{\mathbb{L}^{\infty}}^4$. According to the GN inequality, it follows that
\begin{equation*}
\begin{split}
\int_0^{t\wedge\chi^K}\|\mathbf{u}_i\|_{\mathbb{L}^{\infty}}^4\,\mathrm{d}s\leq C\int_0^{t\wedge\chi^K}\|\mathbf{u}_i\|_{\mathbb{L}^{2}}^2\|\mathbf{u}_i\|_{\mathbb{H}^{1}}^2\,\mathrm{d}s\leq C_K,~\textrm{for}~d=1,\\
\int_0^{t\wedge\chi^K}\|\mathbf{u}_i\|_{\mathbb{L}^{\infty}}^4\,\mathrm{d}s\leq C\int_0^{t\wedge\chi^K}\|\mathbf{u}_i\|_{\mathbb{L}^{2}}^2\|\mathbf{u}_i\|_{\mathbb{H}^{2}}^2\,\mathrm{d}s\leq C_K,~\textrm{for}~d=2.
\end{split}
\end{equation*}
Thus by the Gronwall lemma, $\sup_{s\in[0,t\wedge\chi^K]}\|\mathbf{u}^*(s)\|_{\mathbb{L}^2}^2=0$ $\mathbb{P}$-a.s. By using the monotone convergence theorem and the fact that $\chi^K\nearrow T$ as $K\rightarrow\infty$, it follows that $\mathbb{P}$-a.s.,
$
\sup_{s\in[0,T]}\|\mathbf{u}^*(s)\|_{\mathbb{L}^2}^2=0,
$
which implies the uniqueness. Finally by using the Yamada-Watanabe theorem, we complete the proof.~~$\Box$

\subsection{Invariant measures for $d=1,2$}\label{2subsec8}
\begin{lemma}\label{lem8-1} Let $\mathbf{u}$ be a very weak solution to \eqref{sys1} with properties listed in Theorem \ref{the3}. Then there exists a positive constant $C$ depending on $C_0'$ and $\mathbf{h}$ such that for all $t\geq0$
\begin{equation}\label{1lem8-1}
\begin{split}
\int_0^t\mathbb{E}\|\mathbf{u}(s)\|_{\mathbb{H}^2}^2\,\mathrm{d}s\leq C(1+t).
\end{split}
\end{equation}
\end{lemma}
\noindent\textbf{Proof.} According to Lemma \ref{lem81} and Lemma \ref{lem85}, it follows that
\begin{equation*}
\begin{split}
&\mathbb{E}\int_0^T\|\mathbf{u}\|_{\mathbb{H}^2}^2\,\mathrm{d}s+\sum_{j=1}^{\infty}\mathbb{E}\int_0^T\|-\mathbf{u}\times \mathbf{h}_j+\mathbf{h}_j-\Delta\mathbf{h}_j\|_{\mathbb{L}^2}^2\,\mathrm{d}s<\infty,\\
&\mathbb{E}\int_0^T\|\Delta\mathbf{u}\|_{(\mathbb{H}^{2})^*}^2\,\mathrm{d}s+\mathbb{E}\int_0^T\|\Delta^2\mathbf{u}\|_{(\mathbb{H}^{2})^*}^2\,\mathrm{d}s+\mathbb{E}\int_0^T\|\mathbf{u}\times\Delta\mathbf{u}\|_{(\mathbb{H}^{2})^*}^2\,\mathrm{d}s+\mathbb{E}\int_0^T\|\Delta(\left|\mathbf{u}\right|^2\mathbf{u})\|_{(\mathbb{H}^{2})^*}^2\,\mathrm{d}s<\infty.
\end{split}
\end{equation*}
Choosing $H=\mathbb{L}^2$ and $V=\mathbb{H}^2$, then the process satisfies the assumptions of Pardoux's theorem, which means that the It\^{o}'s formula can be used. Thus by a similar argument as shown in Lemma \ref{lem71}, we easily obtain the estimate \ref{1lem8-1}.~~$\Box$

Let $\mathcal{B}_b(\mathbb{L}^2)$ denote the set of all bounded and Borel measurable functions on $\mathbb{L}^2$. For any $\psi\in\mathcal{B}_b(\mathbb{L}^2)$, $t\geq0$, we can define the transition semigroup $P_t\psi:\mathbb{L}^2\rightarrow\mathbb{R}$ by
\begin{equation}\label{ad8-1}
\begin{split}
P_t\psi(\mathbf{u}_0):=\mathbb{E}\psi(\mathbf{u}(t;\mathbf{u}_0)),~\mathbf{u}_0\in\mathbb{L}^2.
\end{split}
\end{equation}
By a similar argument as shown in Lemma \ref{lem72} and Proposition \ref{pro71}, it is not hard to derive that the the semigroup $\{P_t\}$ is sequentially weakly Feller in $\mathbb{L}^2$. Thus we have the following result.

\begin{theorem}\label{the81}  There exists an invariant measure $\mu$ of the semigroup $\{P_t\}$ defined by \eqref{ad8-1}, such that for any $t\geq0$ and $\psi\in C_b(\mathbb{L}^2_w)$
\begin{equation*}
\begin{split}
\int_{\mathbb{L}^2}P_t\psi(\mathbf{u})\mu(\mathrm{d}\mathbf{u})=\int_{\mathbb{L}^2}\psi(\mathbf{u})\mu(\mathrm{d}\mathbf{u}).
\end{split}
\end{equation*}
\end{theorem}
\noindent\textbf{Proof.} By using the Chebyshev inequality and Lemma \ref{lem8-1}, we derive that for every $T>0$ and $R>0$
\begin{equation*}
\begin{split}
\frac{1}{T}\int_0^T\mathbb{P}\left(\|\mathbf{u}(s;\mathbf{u}_0)\|_{\mathbb{L}^2}>R\right)\,\mathrm{d}s\leq\frac{1}{TR^2}\int_0^T\mathbb{E}\|\mathbf{u}(s;\mathbf{u}_0)\|_{\mathbb{L}^2}^2\,\mathrm{d}s\leq\frac{C+CT}{TR^2}.
\end{split}
\end{equation*}
Since the semigroup $\{P_t\}$ is sequentially weakly Feller in $\mathbb{L}^2$, by using Lemma \ref{lem96}, there exists as least one invariant measure for equation \eqref{sys1}.~~$\Box$\\
\textbf{Proof of Theorem \ref{the3}.}
Theorem \ref{the3} follows from Proposition \ref{pro82} and Theorem \ref{the81}. The entire proof progress is thus completed.~~$\Box$

\section{Appendix}\label{app1}
\begin{lemma}\label{lem92} (\cite{13,21}) Assume that $E$ is a separable Hilbert space, $p\in[2,\infty)$ and $\alpha\in(0,\frac{1}{2})$. Then there exists a constants $C$ depending on $T$ and $\alpha$ such that for any progressively measurable process $\xi=(\xi)_{j=1}^{\infty}$ there holds
\begin{equation*}
\begin{split}
\mathbb{E}\left\|\sum_{j=1}^{\infty}I(\xi_j)\right\|_{W^{\alpha,p}(0,T;E)}^p\leq C\mathbb{E}\int_0^T\left(\sum_{j=1}^{\infty}\|\xi_j(t)\|_{E}^2\right)^{\frac{p}{2}}\,\mathrm{d}t,
\end{split}
\end{equation*}
where $I(\xi_j)$ is defined by
$
I(\xi_j):=\int_0^t\xi_j(s)\,\mathrm{d}W_j(s),~t\geq0.
$
In particular, the trajectories of the process $I(\xi_j)$ belong to $W^{\alpha,2}(0,T;E)$, $\mathbb{P}$-a.s.
\end{lemma}

\begin{lemma}\label{lem93} (\cite{21}) Assume that $B_0\subset B_1\subset B_2$ are Banach spaces, $B_0$ and $B_1$ being reflexive. Assume that the embedding $B_0\subset B_1$ is compact, $q\in(1,\infty)$ and $\alpha\in(0,1)$. Then the embedding
\begin{equation*}
\begin{split}
L^p(0,T;B_0)\cap W^{\alpha,q}(0,T;B_2)\hookrightarrow L^p(0,T;B_1)
\end{split}
\end{equation*}
is compact.
\end{lemma}

\begin{lemma}\label{lem94} (\cite{21}) Assume that $X_0\subset X_1$ are Banach spaces such that the embedding $X_0\subset X_1$ is compact. Assume that $p\in(1,\infty)$, $\alpha\in(0,1)$ and $\alpha p>1$. Then the embedding
\begin{equation*}
\begin{split}
W^{\alpha,p}(0,T;X_0)\subset C([0,T];X_1)
\end{split}
\end{equation*}
is compact.
\end{lemma}

\begin{lemma}\label{lem96} (\cite{33}) Assume that
\begin{enumerate}
\item [(1)] the semigroup $\{P_t\}$ is sequentially weakly Feller in $\mathbb{H}^r,~r=0,1,...$;
\item [(2)] there exists $T_0\geq0$ such that for any $\varepsilon>0$ there exists $R>0$ satisfying
\begin{equation*}
\begin{split}
\sup_{T>T_0}\frac{1}{T}\int_0^T\mathbb{P}\left(\|\mathbf{u}(s;\mathbf{u}_0)\|_{\mathbb{H}^r}>R\right)\,\mathrm{d}s\leq\varepsilon.
\end{split}
\end{equation*}
\end{enumerate}
Then there exists at least one invariant measure for equation \eqref{sys1}.
\end{lemma}

\section*{Conflict of interest statement}

The authors declared that they have no conflicts of interest to this work.

\section*{Data availability}

No data was used for the research described in the article.

\section*{Acknowledgements}
This work was partially supported by the National Natural Science Foundation of China (Grant No. 12231008), and the National Key Research and Development Program of China (Grant No.  2023YFC2206100).

\bibliographystyle{plain}%
\bibliography{SLLBar}

\end{document}